\documentclass[11pt]{article}

\usepackage{amsmath,amssymb,amsfonts}
\usepackage{pstricks,pst-node}
\usepackage{xr}

\numberwithin{equation}{section}
\newtheorem{thm}{Theorem}[section] 
\newtheorem{prp}[thm]{Proposition}
\newtheorem{lmm}[thm]{Lemma}   
\newtheorem{crl}[thm]{Corollary} 
\newtheorem{dfn}[thm]{Definition} 

\renewcommand{\frak}{\mathfrak}
\renewcommand{\Bbb}{\mathbb}

\def\C{\mathbb C}
\def\cC{\mathcal C}
\def\d{\mathfrak d}
\def\D{\mathfrak D}

\def\E{\mathbb E}
\def\F{\mathfrak F}

\def\I{\mathfrak i}
\def\J{\mathcal J}
\def\L{\mathfrak L}
\def\M{\mathfrak M}
\def\cM{\mathcal M}
\def\O{\mathcal O}
\def\Q{\mathbb Q}
\def\R{\mathbb R}
\def\S{\mathcal S}
\def\T{\mathcal T}
\def\U{\mathcal U}
\def\fU{\mathfrak U}
\def\V{\mathcal V}
\def\W{\mathcal W}
\def\X{\mathfrak X}
\def\Z{\mathbb Z}
\def\cZ{\mathcal Z}

\def\e_ref#1{(\ref{#1})}
\def\under#1{\underline{#1}}
\def\ov#1{\overline{#1}}
\def\ti#1{\tilde{#1}}
\def\wt#1{\widetilde{#1}}
\def\lra{\longrightarrow}
\def\Lra{\Longrightarrow}

\def\lan{\langle}
\def\ran{\rangle}

\def\blr#1{\big\lan{#1}\big\ran}
\def\llrr#1{\lan\!\lan{#1}\ran\!\ran}
\def\bigllrr#1{\big\lan\!\big\lan{#1}\big\ran\!\big\ran}

\def\al{\alpha}

\def\de{\delta}
\def\ep{\epsilon}
\def\ga{\gamma}
\def\io{\iota}
\def\ka{\kappa}
\def\la{\lambda}
\def\na{\nabla}

\def\si{\sigma}

\def\ve{\varepsilon}
\def\ups{\upsilon}
\def\vph{\varphi}

\def\vr{\varrho}
\def\vt{\vartheta}
\def\ze{\zeta}

\def\Ga{\Gamma}
\def\La{\Lambda}
\def\Om{\Omega}
\def\Si{\Sigma}

\def\rk{\textnormal{rk}\,}

\def\Im{\textnormal{Im}\,}
\def\ev{\textnormal{ev}}
\def\reg{\textnormal{reg}}
\def\rig{\textnormal{rig}}

\def\Hom{\textnormal{Hom}}

\def\P{\Bbb{P}^n}
\def\PP{\Bbb{P}^2}
\def\PPP{\Bbb{P}^3}
\def\Pf{\Bbb{P}^4}
\def\i{\infty}
\def\eset{\emptyset}
\def\bpar{\bar\partial}

\begin{document}

\title{On the Genus-One Gromov-Witten Invariants\linebreak 
of a Quintic Threefold}

\author{Jun Li\thanks{Partially supported by an NSF grant}~~and 
Aleksey Zinger\thanks{Partially supported by an NSF Postdoctoral Fellowship}}

\date{\today}
\maketitle

\begin{abstract}
\noindent
We rederive a relation between the genus-one GW-invariants
of a quintic threefold in $\Pf$ and the genus-zero and genus-one GW-invariants of~$\Pf$.
In contrast to the more general derivation in our previous paper,
the present derivation relies on a widely believed, but still unproven,
statement concerning rigidity of holomorphic curves in Calabi-Yau threefolds.
On the other hand, this paper's derivation is more direct and geometric.
It requires a bit more effort, but relies on less outside work.
\end{abstract}

\thispagestyle{empty}

\tableofcontents

\section{Introduction}
\label{intro_sec}

\subsection{Summary}
\label{statement_subs}

\noindent
Suppose $\ga\!\lra\!\P$ is the tautological line bundle, $a\!\in\!\Z^+$, and
$$\L=\ga^{*\,\otimes a}.$$
If $s\!\in\!H^0(\P;\L)$ is a generic holomorphic section,
$$Y\equiv s^{-1}(0)$$
is a smooth hypersurface in $\P$.
It has long been known how to express the genus-zero Gromov-Witten invariants
of $Y$ in terms of the genus-zero GW-invariants of~$\P$; 
see~\e_ref{genus0_e} below for a special case.
The latter can be computed using the classical localization theorem of~\cite{ABo}.
In~\cite{LZ}, we prove a genus-one analogue of~\e_ref{genus0_e} for an arbitrary
hypersurface~$Y$.
The proof itself is rather simple.
However, it relies on the constructions of {\tt reduced} genus-one GW-invariants
in~\cite{g1comp2} and of euler classes of certain natural cones in a setting more
general than in~\cite{g1cone}.
The latter in fact constitutes most of~\cite{LZ}.\\

\noindent
In this paper, we rederive a genus-one analogue of~\e_ref{genus0_e} for a quintic threefold
$Y$ in $\Pf$, i.e.~for $a\!=\!5$ in the above notation, in a direct, 
albeit more laborious, way.
In order to do this, we will assume a certain rigidity property
for genus-zero and genus-one $J$-holomorphic curves in a quintic threefold; 
see the next subsection.
While it is not known whether the rigidity property is satisfied,
it is widely believed to be the case, at least for curves up genus two or three.
Our derivation generalizes to arbitrary Calabi-Yau complete-intersection threefolds 
in projective spaces.
It can be used for Fano complete-intersection threefolds as well,
but in such cases it can be obtained by taking $\nu\!=\!0$ in 
Subsection~\ref{g1gw-hyperppf_subs} of~\cite{LZ}.
In the Calabi-Yau cases, this cannot be done and this paper's derivation is
different from that in Subsection~\ref{g1gw-hyperppf_subs} of~\cite{LZ}.\\

\noindent
Quintic threefolds, as well as other Calabi-Yau manifolds,
play a prominent role in theoretical physics.
As a result physicists have made a number of important predictions concerning CY-manifolds.
Some of these predictions have been verified mathematically; others have~not.
This paper indicates that one of them fits in nicely with known mathematical facts.\\

\noindent
If $X$ is a Kahler manifold, $g$ and $k$ are nonnegative integers, and $A\!\in\!H_2(X;\Z)$,
we denote by $\ov\M_{g,k}(X,A)$ the moduli space of (equivalence classes of) stable holomorphic maps from genus-$g$ curves with $k$~marked points in the homology class~$A$.
Let
$$\ov\M_g(X,A)=\ov\M_{g,0}(X,A).$$
If $\io\!:Y\!\lra\!\P$ is an inclusion and $\ell$ is the homology class of a line
in $\P$, let
$$\ov\M_{g,k}(Y,d)=\bigcup_{\io_*A=d\ell}\!\!\ov\M_{g,k}(Y,A).$$
If $Y$ is a Calabi-Yau threefold, the virtual, or expected dimension, of $\ov\M_g(Y,d)$
is zero.
The virtual degree of $\ov\M_g(Y,d)$ is the genus-$g$ degree-$d$ GW-invariant of~$Y$.
If $Y$ is a quintic threefold, we denote this invariant by~$N_g(d)$.\\

\noindent
Let 
$$\pi_g^d\!\!: \fU_g(\P,d)\lra\ov\M_g(\P,d)   \qquad\hbox{and}\qquad
\ev_g^d\!:\fU_g(\P,d)\lra\P$$ 
be the semi-universal family and the natural evaluation~map.
In other words, the fiber of $\pi_g^d$ over $[\cC,u]$
is the curve~$\cC$, while 
$$\ev_g^d\big([\cC,u;z]\big)=u(z) \qquad\hbox{if}\quad z\!\in\!\cC.$$
We define a section~$s_g^d$ of the sheaf 
$\pi_{g*}^d\ev_g^{d*}\L\!\lra\!\ov\M_g(\P,d)$ by
$$s_g^d\big([{\cal C},u]\big)=[s\circ u].$$
If $Y\!=\!s^{-1}(0)$, $\ov\M_g(Y,d)$ is the zero set of this section.\\

\noindent
If $a\!=\!5$, it has long been know that
\begin{equation}\label{genus0_e}
N_0(d) =\blr{ e\big(\pi_{0*}^d\ev_0^{d*}\L\big),\big[\ov\M_0(\Pf,d)\big]}.
\end{equation}
The moduli space $\ov\M_0(\Pf,d)$ is a smooth orbivariety  and
\begin{equation}\label{g0sheaf_e}
\pi_{0*}^d\ev_0^{d*}\L\lra\ov\M_0(\Pf,d)
\end{equation}
is a locally free sheaf, i.e.~a vector bundle.
Furthermore,
$$\dim_{\C}\ov\M_0(\Pf,d)=5d+1 \qquad\hbox{and}\qquad
\rk\!_{\C}\,\pi_{0*}^d\ev_0^{d*}\L=5d+1.$$
Thus, the right-hand side of~\e_ref{genus0_e} is well-defined.
It can be computed via the classical localization theorem of~\cite{ABo}.
The complexity of this computation increases quickly with the degree~$d$,
but a closed formula has been obtained in 
\cite{Ber}, \cite{Ga}, \cite{Gi}, \cite{Le}, and~\cite{LLY}.\\

\noindent
If $g\!>\!0$, the sheaf 
$\pi_{g*}^d\ev_g^{d*}\L\!\lra\!\ov\M_g(\Pf,d)$ 
is not locally free and does not define an euler class.
Thus, the right-hand side of~\e_ref{genus0_e} does not even
make sense if $0$ is replaced by $g\!>\!0$.
Instead one might try to generalize~\e_ref{genus0_e} as
\begin{equation}\label{allgenus_e}
N_g(d) \stackrel{?}{=}\blr{ 
e\big(R^0\pi_{g*}^d\ev_{g}^{d*}\L-R^1\pi_{g*}^d\ev_g^{d*}\L\big),
\big[\ov\M_g(\Pf,d)\big]^{vir}},
\end{equation}
where $R^i\pi_{g*}^d\ev_g^{d*}\L\!\lra\!\ov\M_g(\Pf,d)$
is the $i$th direct image sheaf.
The right-hand side of~\e_ref{allgenus_e} can be computed via 
the virtual localization theorem of~\cite{GrP1}.
However, 
$$N_1(d)\neq \blr{e\big(R^0\pi_{1*}^d\ev_1^{d*}\L\!-\!
R^1\pi_{1*}^d\ev_1^{d*}\L\big),
\big[\ov\M_1(\Pf,d)\big]^{vir}},$$
according to a low-degree check of~\cite{GrP2} and~\cite{K}.\\

\noindent
Let 
$$\M_g^0(\Pf,d)=\big\{[\cC,u]\!\in\!\ov\M_g(\Pf,d)\!:
\cC~\hbox{is smooth}\big\}.$$
We denote by $\ov\M_g^0(\Pf,d)$ the closure of $\M_g^0(\Pf,d)$ in~$\ov\M_g(\Pf,d)$.
If $g\!>\!0$, $\ov\M_g^0(\Pf,d)$ is one of the many irreducible components of 
the moduli space~$\ov\M_g(\Pf,d)$.

\begin{thm}
\label{main_thm}
If $d$ is a positive integer, $\L\!=\!\ga^{*\otimes 5}\!\lra\!\Pf$, 
$$\pi_1^d\!:{\frak U}_1(\Pf,d)\lra\ov\M_1^0(\Pf,d)
\quad\hbox{and}\quad
\ev_1^d\!:{\frak U}_1(\Pf,d)\lra\Pf$$ 
are the semi-universal family and the natural evaluation~map, respectively,
then the euler class of the~sheaf
$$\pi_{1*}^d\ev_1^{d*}\L\lra\ov\M_1^0(\Pf,d)$$
is well-defined.
Furthermore, 
\begin{equation}\label{main_thm_e}
N_1(d)=\frac{1}{12}N_0(d)+
\blr{e\big(\pi_{1*}^d\ev_1^{d*}\L\big),\big[\ov\M_1^0(\Pf,d)\big]}.
\end{equation}\\
\end{thm}

\noindent
The moduli space $\ov\M_1^0(\Pf,d)$ is not a smooth orbifold.
Nevertheless, it determines a fundamental class in $H_{10d}(\ov\M_1^0(\Pf,d);\Q)$,
as its singularities are fairly simple.
The~sheaf
\begin{equation}\label{sheaf_e}
\pi_{1*}^d\ev_1^{d*}\L\lra\ov\M_1^0(\Pf,d)
\end{equation}
is not locally free.
Nevertheless, its euler class is well-defined.
In other words, the euler class of every desingularization of this sheaf is the~same,
in the sense described in Subsection~\ref{g1cone-appr_subs} of~\cite{g1cone}.
The last expression in~\e_ref{main_thm_e} can be computed
via the classical localization theorem.
Of course, the singularities of the space $\ov\M_1^0(\Pf,d)$
cause additional complications.
However, since these singularities can be understood,
these complications can be handled.
A desingularization of $\ov\M_1^0(\Pf,d)$, 
i.e.~a smooth orbifold $\wt\M_1^0(\Pf,d)$ and a~map
$$\ti\pi\!:\wt\M_1^0(\Pf,d)\lra\ov\M_1^0(\Pf,d),$$
which is biholomorphic onto $\M_1^0(\Pf,d)$, is constructed in~\cite{VZ}.
This desingularization of $\ov\M_1^0(\Pf,d)$ comes with a desingularization
of the sheaf~\e_ref{sheaf_e},
i.e.~a vector bundle 
$$\ti\V \lra \wt\M_1^0(\Pf,d) \qquad\hbox{s.t.}\qquad
\ti\pi_*\ti\V=\pi_{1*}^d\ev_1^{d*}\L.$$
In particular,
$$\blr{e\big(\pi_{1*}^d\ev_1^{d*}\L),
\big[\ov\M_1^0(\Pf,d)\big]}
=\blr{e(\ti\V),\big[\wt\M_1^0(\Pf,d)\big]}.$$
The localization theorem of~\cite{ABo} is directly applicable
to the right-hand side of this equality.\\

\noindent
Using Theorem~\ref{main_thm} and the desingularization constructed in~\cite{VZ},
we have computed the numbers $N_1(d)$ for $d\!=\!1,2,3,4$.
The results agree with those predicted in~\cite{BCOV};
see Subsection~\ref{g1gw-mainthmcrl_subs} in~\cite{LZ} for more details.\\

\noindent
From the point of view of symplectic topology as described in
\cite{FuO} and~\cite{LT}, the numbers $N_g(d)$ 
can be interpreted as the euler class of a vector bundle, albeit of an infinite-rank
vector bundle over a space of the ``same" dimension.
As in the finite-dimensional case,
this euler class is the number of zeros, counted with appropriate
multiplicities, of a transverse (multivalued, admissible) section.\\

\noindent
In brief, we prove Theorem~\ref{main_thm} by slightly perturbing the complex structure~$J_0$
on~$\Pf$, then expressing each of the three terms appearing in~\e_ref{main_thm_e} 
as the number of zeros of a transverse section of a vector bundle
and comparing the results for the two sides of~\e_ref{main_thm_e}.
There are a vector bundle $\F\!\lra\!\X$, possibly
of infinite rank, and a section $\vph$ of $\F$ associated to each of three terms.
The zero set of~$\vph$ is easy to describe.
However, $\vph$ is not transverse to the zero set.
We determine the number $\cC_{\cZ}(\vph)$ of zeros of $\vph\!+\!\ve$, for 
a small generic multisection~$\ve$, that lie near each stratum~$\cZ$ of~$\vph^{-1}(0)$.
These numbers in turn determine the contribution of each $J$-holomorphic curve in $Y$ to
the three numbers in~\e_ref{main_thm_e}.
We will see that every such curve contributes equally to the 
two sides of~\e_ref{main_thm_e}.\\

\noindent
Theorem~\ref{main_thm} follows immediately from Propositions~\ref{g0_prp1}-\ref{cone_prp1}
and separately from Propositions~\ref{g0_prp2}-\ref{cone_prp2}.
The first three propositions are easier to state and can be deduced from 
the last three propositions.
While the statements of Propositions~\ref{g0_prp2}-\ref{cone_prp2} are more technical,
they are easier to prove.

\subsection{Rigidity Properties}
\label{rigidassump_subs}

\noindent
Throughout the rest of the paper, $Y$ will denote a quintic threefold in $\Pf$.
If $J$ is an almost complex structure on $Y$, $(\Si,j)$ is a Riemann surface, and 
$u\!:\Si\!\lra\!Y$ is a $J$-holomorphic map, let
$$D_{J,u}\!: \Ga(\Si;u^*TY)\lra
\Ga\big(\Si;\La^{0,1}_{J,j}T^*\Si\!\otimes\!u^*TY\big)$$
be the linearization of $\bpar_J$-operator at $u$;
see Subsection~\ref{review_subs}.

\begin{dfn}
\label{rigid_dfn}
An almost complex structure $J$ on $Y$ {\tt satisfies the genus-$g$ rigidity property}
if for every smooth connected genus-$g$ Riemann surface $(\Si,j)$ and 
nonconstant $J$-holomorphic map $u\!:\Si\!\lra\!Y$\\
${}\quad$ ($J_Y1$) $u(\Si)$ is a smooth curve;\\
${}\quad$ ($J_Y2$) $\ker D_{J,u}\!\subset\!\Ga(\Si;u^*Tu(\Si))$.
\end{dfn}

\noindent
If $J$ satisfies the genus-$g$ rigidity property, 
all genus-$g$ $J$-holomorphic curves in~$Y$ are smooth and isolated.
We denote by $\J(Y)$  the space of all $C^1$-smooth almost complex structures on~$Y$,
with the $C^1$-topology, and by $\J_{\rig}^g(Y)\!\subset\!\J(Y)$ the subspace of
almost complex structures that satisfy the genus-$g$ rigidity property.\\

\noindent
{\bf Rigidity Conjecture}~ For all $g$ and all Calabi-Yau threefolds~$Y$, 
$\J_{\rig}^g(Y)$ is dense in $\J(Y)$.\\

\noindent
{\bf Rigidity Assumption}~ If $Y$ is a quintic threefold,
the closure of $\J_{\rig}^0(Y)\!\cap\!\J_{\rig}^1(Y)$ in
$\J(Y)$ contains~$J_0$.\\

\noindent
We note that $\J_{\rig}^g(Y)$ is open in $\J(Y)$.
Thus, the $g\!=\!0,1$ cases of the Rigidity Conjecture imply our Rigidity Assumption.\\

\noindent
Since the expected dimension of the moduli space $\ov\M_g(Y,d;J)$ of 
genus-$g$ degree-$d$ $J$-holomorphic maps into $Y$  is zero,
it is easy to show that the property~($J_Y1$) of Definition~\ref{rigid_dfn}
is satisfied by a generic almost complex structure~$J$.
However, despite years of attempts, this has not been shown to be the case
for~($J_Y2$), even for $g\!=\!0$.
Nevertheless, this is believed to be case, though with some hesitation 
for $g$ above $2$ or~$3$.\\

\noindent
For each $J\!\in\!\J(Y)$, let $\S_g^d(Y;J)$ be the set of 
$J$-holomorphic genus-$g$ degree-$d$ (simple) {\it curves} in~$Y$. 
If $J\!\in\!\J_{\rig}^g(Y)$, this set is finite.
By Propositions~\ref{g0_prp1} and~\ref{g1_prp1} below,
the number of elements in $\S_g^d(Y;J)$, counted with appropriate signs,
is independent of $J\!\in\!\J_{\rig}^g(Y)$ for $g\!=\!0,1$.
We denote this number by~$n_g(d)$.

\begin{prp}
\label{g0_prp1}
For all $d\!\in\!\Z^+$,  
$$N_0(d)=\sum_{\si|d}\frac{n_0(d/\si)}{\si^3}.$$
\end{prp}

\begin{prp}
\label{g1_prp1}
For all $d\!\in\!\Z^+$,
$$N_1(d)= \frac{1}{12}\sum_{\si|d}\frac{n_0(d/\si)}{\si}
+\sum_{\si|d}\frac{n_1(d/\si)}{\si}.$$
\end{prp}

\begin{prp}
\label{cone_prp1}
If $d$, $\L$, $\pi_1^d$, and $\ev_1^d$ are as in the statement of Theorem~\ref{main_thm}, 
then
$$\blr{e(\pi_{1*}^d\ev_1^{d*}\L),\big[\ov\M_1^0(\Pf,d)\big]}
=\frac{1}{12}\sum_{\si|d}\frac{\si^2\!-\!1}{\si^3}n_0(d/\si)+
\sum_{\si|d}\frac{n_1(d/\si)}{\si}.$$\\
\end{prp}

\noindent
We do not prove these three propositions as stated,
since this is not necessary for the proof of Theorem~\ref{main_thm}.
Instead, we prove the less elegant and more notationally involved Propositions~\ref{g0_prp2}-\ref{cone_prp2}, that also imply Theorem~\ref{main_thm}.
Propositions~\ref{g0_prp1}, \ref{g1_prp1}, and \ref{cone_prp1} 
can be derived from Propositions~\ref{g0_prp2}, \ref{g1_prp2}, and \ref{cone_prp2},
respectively; see the end of Subsection~\ref{outline_subs2}.

\section{Preliminaries}
\label{prelim_sec}

\subsection{Review of Key Definitions}
\label{review_subs}

\noindent
In this subsection, we give geometric definitions of the three terms  
that appear in~\e_ref{main_thm_e}.
The construction of the Gromov-Witten invariants described below
is a slight variation on that of~\cite{FuO} and~\cite{LT}, but 
it is easy to see the only difference is in the presentation.
Below we use the term {\it multisection}, or {\it multivalued section},  
of a vector orbi-bundle as defined in Section~3 of~\cite{FuO}.\\

\noindent
If $X$ is a smooth submanifold of $\Bbb{P}^n$, we denote by $\X_g(X,d)$ the space of 
equivalence classes of stable degree-$d$ smooth maps from genus-$g$ Riemann surfaces to~$X$. 
Let $\X_g^0(X,d)$ be the subset of $\X_g(X,d)$
consisting of stable maps with smooth domains. 
The spaces $\X_g(X,d)$ are topologized
using $L^p_1$-convergence on compact subsets of smooth points of the domain
and certain convergence requirements near the nodes.
Here and throughout the rest of the paper, $p$~denotes a real number 
greater than~two.
The spaces $\X_g(X,d)$ can be stratified by the smooth infinite-dimensional orbifolds
$\X_{\T}(X)$
of stable maps from domains of the same geometric type.
The closure of the main stratum, $\X_g^0(X,d)$, is $\X_g(X,d)$.\\

\noindent
If $J$ is an almost complex structure on $\Bbb{P}^n$, let 
$$\Ga_g^{0,1}(X,d;J)\!\lra\!\X_g(X,d)$$
be the bundle of $(TX,J)$-valued $(0,1)$-forms. 
In other words, the fiber of $\Ga_g^{0,1}(X,d;J)$ over a point
$[b]\!=\![\Si,j;u]$ in $\X_g(X,d)$ is the space
$$\Ga_g^{0,1}(X,d;J)\big|_{[b]}=\Ga^{0,1}(b;TX;J)\big/\hbox{Aut}(b),
\quad\hbox{where}\quad
\Ga^{0,1}(b;TX;J)=\Ga\big(\Si;\La_{J,j}^{0,1}T^*\Si\!\otimes\!u^*TX).$$
Here $j$ is the complex structure on $\Si$, the domain of the smooth map~$u$.
The bundle \hbox{$\La_{J,j}^{0,1}T^*\Si\!\otimes\!u^*TX$} over $\Si$
consists of $(J,j)$-antilinear homomorphisms:
$$\La_{J,j}^{0,1}T^*\Si\!\otimes\!u^*TX=\big\{
\al\!\in\!\hbox{Hom}(T\Si,u^*TX)\!:\al\!\circ\!j\!=\!-\!J\!\circ\!\al\big\}.$$
The total space of the bundle $\Ga_g^{0,1}(TX,d;J)\!\lra\!\X_g(X,d)$ 
is topologized using $L^p$-convergence on compact subsets of smooth points of the domain
and certain convergence requirements near the nodes.
The restriction of $\Ga_g^{0,1}(TX,d;J)$ to each stratum
$\X_{\T}(X)$ is a smooth vector orbibundle of infinite rank.\\

\noindent
We define a continuous section of the bundle 
$\Ga_g^{0,1}(TX,d;J)\!\lra\!\X_g(X,d)$  by
$$\bar{\partial}_J\big([\Si,j;u]\big)=\frac{1}{2}\big(
du+J\!\circ\!du\!\circ\!j\big).$$
By definition, the zero set of this section
is the moduli space $\ov\M_g(X,d;J)$ 
of equivalence classes of stable $J$-holomorphic degree-$d$ maps 
from genus-$g$ curves into~$X$.
The restriction of~$\bpar_J$ to each stratum of $\X_g(X,d)$ is smooth.
For each element $[b]\!=\![\Si,j,u]$ of $\X_g(X,d)$,
we~put 
$$D_{J,b}\xi=\frac{1}{2}\big(\na^X\xi+J\circ\na^X\xi\circ j\big)
+\frac{1}{2}(\na_{\xi}^XJ)\circ du\circ j
~~\hbox{if}~\xi\!\in\!\Ga(b;TX)\!\equiv\!\Ga(\Si;u^*TX),$$
where $\na^X$ denotes the Levi-Civita connection of a $J$-compatible metric on~$X$.
The linear operator $D_{J,b}$ describes the restriction of a linearization 
of $\bpar_J$ at $[b]$  to a finite-codimensional subspace of 
the tangent bundle of the stratum $\X_{\T}(X)$ of  $\X_g(X,d)$ containing~$[b]$.\\

\noindent
The section $\bpar_J\!:\X_g(X,d)\!\lra\!\Ga_g^{0,1}(X,d;J)$ 
is Fredholm, i.e.~its linearization at every point of $\bpar_J^{-1}(0)$
has finite-dimensional kernel and cokernel.
The index of~$\bpar_J$ at a point of $\X_g^0(X,d)$
is the expected dimension of the moduli space~$\ov\M_g(X,d;J)$.
If $X\!=\!Y$, this expected dimension is~$0$.
By definition, 
\begin{equation}\label{gw_dfn_e}
N_g(d)=\, ^{\pm}\! \big|\{\bpar_J\!+\!\ve\}^{-1}(0)\big|,
\end{equation}
where $\ve$ is a small multivalued perturbation such that 
$\bpar_J\!+\!\ve$ is transverse to the zero~set
along each stratum $\X_{\T}(Y)$ of  $\X_g(Y,d)$ 
and 
$$^{\pm}\!\big|\{\bpar_J\!+\!\ve\}^{-1}(0)\big|$$
is the number of elements in the finite set $\{\bpar_J\!+\!\ve\}^{-1}(0)$, 
counted with appropriate multiplicities.
By the transversality condition,
$$\{\bpar_J\!+\!\ve\}^{-1}(0) \subset \X_g^{0}(Y,d).$$
The smallness condition implies in particular that the set
$\{\bar{\partial}_J\!+\!\ve\}^{-1}(0)$ is close to~$\bar{\partial}_J^{-1}(0)$.
Since the set $\bar{\partial}_J^{-1}(0)$ is compact, it follows that 
the set $\{\bpar_J\!+\!\ve\}^{-1}(0)$ is also compact.
Let ${\cal A}_g^d(\bpar_J)$ denote the set of all perturbations~$\ve$ 
of $\bpar_J$ that satisfy the two conditions above.
Such perturbations will be called $\bpar_J$-admissible.
Below we will refer to the number in~\e_ref{gw_dfn_e} as the euler
class of the tuple 
$$\V_g^d(\bpar;J)
\equiv \big(\X_g(Y,d),\Ga_g^{0,1}(Y,d;J),\pi;\bpar_J,{\cal A}_g^d(\bpar_J)\big).$$
This euler class depends on the Fredholm homotopy
class of the section~$\bpar_J$.\\

\noindent
We now describe the last term in~\e_ref{main_thm_e} in a similar way.
If $\L\!\lra\!\Pf$ is as in Theorem~\ref{main_thm}, let 
$\Ga_g(\L,d)\!\lra\!\X_g(\Pf,d)$ be the cone 
such that the fiber of $\Ga_g(\L,d)$ over $[b]\!=\![\Si,j;u]$ in 
$\X_g(\Pf,d)$ is the Banach space
$$\Ga_g(\L,d)\big|_{[b]}=\Ga(b;\L)\big/\hbox{Aut}(b),
\qquad\hbox{where}\qquad \Ga(b;\L)=L^p_1(\Si;u^*\L),$$
and the topology on $\Ga_g(\L,d)$ in defined analogously
to the topology on~$\Ga_g(\Pf,d)$.
Let $\na$ denote the hermitian connection in the line bundle $\L\!\lra\!\Pf$ 
induced from the standard connection on the tautological line bundle over~$\Pf$.
If $(\Si,j)$ is a Riemann surface and $u\!:\Si\!\lra\!\Pf$ is a smooth map, 
let
$$\na^u\!:\Ga(\Si;u^*\L)\lra \Ga\big(\Si;T^*\Si\!\otimes\!u^*\L\big)$$
be the pull-back of $\na$ by $u$.
If $b\!=\!(\Si,j;u)$,
we define the corresponding $\bpar$-operator by
\begin{equation}\label{vdfn_e}
\bpar_{\na,b}\!:\Ga(\Si;u^*\L)\lra
\Ga\big(\Si;\La_{\I,j}^{0,1}T^*\Si\!\otimes\!u^*\L\big),
\quad
\bpar_{\na,b}\xi=\frac{1}{2}
\big(\na^u\xi+\I\na^u\xi\circ j\big),
\end{equation}
where $\I$ is the complex multiplication in the bundle $u^*\L$. 
Let
$$\V_g^d=\big\{[b,\xi]\!\in\!\Ga_g(\L,d)\!: [b]\!\in\!\X_g(\Pf,d),~
\xi\!\in\!\ker\bpar_{\na,b}\!\subset\Ga_g(b;\L)\big\}
\subset\Ga_{g,k}(\L,d).$$
The cone $\V_g^d\!\lra\!\X_{g,k}(\Pf,d)$ inherits 
its topology from~$\Ga_g(\L,d)$.\\

\noindent
Let $\ov\M_1^0(\Pf,d;J)\!\subset\!\ov\M_1(\Pf,d;J)$ 
denote the closed subset containing the~set
$$\M_1^0(\Pf,d;J)=\big\{[\cC,u]\!\in\!\ov\M_1(\Pf,d;J)\!:
\cC~\hbox{is smooth}\big\},$$
which is defined in~\cite{g1comp}.
If the almost complex structure $J$ is sufficiently close to~$J_0$, $\ov\M_1^0(\Pf,d;J)$ 
is the closure of $\M_1^0(\Pf,d;J)$ in~$\ov\M_1(\Pf,d;J)$.
Furthermore, in such a case,
$\M_1^0(\Pf,d;J)$ is a smooth orbifold of dimension~$10d$,
while  $\partial\ov\M_1^0(\Pf,d;J)$ is a finite union of smooth
orbifolds of dimension at most $10d\!-\!2$.
On the other hand, $\V_1^d|_{\M_1^0(\Pf,d;J)}$
is a complex vector orbibundle of rank~$5d$.
The last term in~\e_ref{main_thm_e} is the number of zeros, 
counted with appropriate multiplicities,
of any continuous multisection $\vph$ of the cone~$\V_1^d$
over $\ov\M_1^0(\Pf,d;J)$
such that $\vph^{-1}(0)$ is contained in $\M_1^0(\Pf,d;J)$
and $\vph|_{\M_1^0(\Pf,d;J)}$ is smooth and transverse to the zero~set;
see Subsections~\ref{g1cone-appr_subs} and~\ref{g1cone-mainres_subs} in~\cite{g1cone}.
Proposition~\ref{g1cone-g1cone_prp} in~\cite{g1cone}
guarantees that a section~$\vph$ satisfying the two conditions exists.
In our case, it is more convenient to think of $\vph$ as $s_1^d\!+\!\ve$, 
where $\ve$ is a multivalued perturbation of $s_1^d$.
We denote by ${\cal A}_1^d(s;J)$ the set of all perturbations~$\ve$ of~$s_1^d$
such that $s_1^d\!+\!\ve$ satisfies the two conditions above.
Such perturbations~$\ve$ will be called $s_1^d$-admissible.
Let
$$\V_1^d(s;J)\equiv\big(\ov\M_1^0(\Pf,d;J),\V_1^d,\pi;
s_1^d,{\cal A}_1^d(s)\big).$$
This tuple will be the focus of Section~\ref{cone_sec}.\\

\noindent
{\it Remark:} Since $Y$ is a semi-positive symplectic manifold,
one can define the numbers $N_g(d)$ without using 
the infinite-rank orbibundles~$\Ga_g^{0,1}(Y,d;J)$; see~\cite{RT}.
However, there would be no effect on the proofs of 
Propositions~\ref{g0_prp1}-\ref{cone_prp1} and \ref{g0_prp2}-\ref{cone_prp2}, and 
the construction described above appears more natural in 
the present context, even though it involves more complicated objects.

\subsection{Components of the Proof}
\label{outline_subs2}

\noindent
We now set up additional notation that allows us to state 
more notationally involved, but also easier-to-prove, versions of 
Propositions~\ref{g0_prp1}-\ref{cone_prp1}.\\

\noindent
By Theorems~\ref{g1comp-reg_thm} and~\ref{g1comp-str_thm} in~\cite{g1comp},
there exists $\de(d)\!\in\!\Bbb{R}^+$ with the property that
if $J$ is an almost complex structure on $\Pf$ such that
$\|J\!-\!J_0\|_{C^1}\!\le\!\de(d)$, 
then $J$ is genus-one $d\ell$-regular in the sense of Definition~\ref{g1comp-g1reg_dfn}
in~\cite{g1comp}.
This regularity condition
implies that the  moduli spaces $\ov\M_{0,k}(\Pf,d;J)$ and $\ov\M_{1,k}(\Pf,d;J)$ 
have the same stratification structure as the moduli spaces
$$\ov\M_{0,k}(\Pf,d)\equiv\ov\M_{0,k}(\Pf,d;J_0)
\qquad\hbox{and}\qquad
\ov\M_{1,k}(\Pf,d)\equiv\ov\M_{1,k}(\Pf,d;J_0),$$
respectively.
In addition, by Theorem~\ref{g1cone-main_thm} in~\cite{g1cone}, 
$\de(d)\!\in\!\R^+$ can be chosen so that the euler class of the cone
$$\V_1^d\lra \ov\M_1^0(\Pf,d;J)$$
is well-defined and 
\begin{equation}\label{g1cone_e}
\blr{e(\V_1^d),\big[\ov\M_1^0(\Pf,d;J)\big]}
=\blr{e(\V_1^d),\big[\ov\M_1^0(\Pf,d)\big]},
\end{equation}
if $\|J\!-\!J_0\|_{C^1}\!\le\!\de(d)$.\\

\noindent
If $J$ is an almost complex structure on $\Pf$ and
$[\Si,u]\!\in\!\ov\M_1(\Pf,d;J)$, we~put
$$s_1^d\big([\Si,u]\big)=[s\circ u]\in\Ga(\L,d)\big|_{[\Si,u]}.$$
If $J$ is {\tt $\na s$-equivalent} to $J_0$, i.e.
$$\na s\circ J_0=\na s\circ J\in \Ga\big(\Pf;\Hom_{\R}(T\Pf,\L)\big),$$
then $s_1^d([\Si,u])\!\in\!\V_1^d|_{[\Si,u]}$.
Thus, in such a case, we obtain a continuous section of the~cone
$$\V_1^d\lra \ov\M_1^0(\Pf,d;J),$$
which restricts to a smooth section over each stratum of $\ov\M_1^0(\Pf,d;J)$.
Note that
\begin{equation}\label{zeroset_e}
\big\{s_1^d\big|_{\ov\M_1^0(\Pf,d;J)}\big\}^{-1}(0)
=\ov\M_1^0(Y,d;J) =\ov\M_1^0(\Pf,d;J)\cap \ov\M_1(Y,d;J).
\end{equation}\\

\noindent
Since the $(\na,J_0)$-holomorphic section~$s$ of Subsection~\ref{statement_subs}
is transverse to the zero set in~$\L$,
the $(\I,J_0)$-linear map
$$\na s\!\!:T\Pf\lra\L$$
does not vanish along $Y\!=\!s^{-1}(0)$.
Let $U_s$ be a small neighborhood of~$Y$ in~$\Pf$ such that 
$\na s$ does not vanish over~$U_s$.
The kernel of $\na s$ over~$U_s$ is then a rank-three complex subbundle
of $(T\Pf,J_0)|_{U_s}$, which restricts to $TY$ along~$Y$.
We denote this subbundle by $\ti{T}Y$.
If $J$ is an almost complex structure on~$\Pf$ such~that\\
${}\quad$ ($J1$) $J\!=\!J_0$ on $\Pf\!-\!U_s$;\\
${}\quad$ ($J2$) $J(\ti{T}Y)\!=\!\ti{T}Y$ and $J\!=\!J_0$ on $T\Pf|_{U_s}/\ti{T}Y$,\\
then $J_0$ and $J$ are $\na s$-equivalent.
Thus, every almost complex structure $J_Y$ on $Y$ extends to 
an almost complex structure $J$ on $\Pf$ which is $\na s$-equivalent to~$J_0$.
Furthermore, such an extension can be chosen so~that
\begin{equation}\label{jext_e}
\|J-J_0\|_{C^1} \le 2\big\|J_Y-J_0|_{TY}\big\|_{C^1}.
\end{equation}
We denote by $\J_{\rig}(s)$ the set of almost complex structures~$J$
on $\Pf$ such that $J$ is $\na s$-equivalent to $J_0$ and
$J_Y\!\equiv\!J|_{TY}$ is an element of $\J_{\rig}^0(Y)\!\cap\!\J_{\rig}^1(Y)$.
By the above and the Rigidity Assumption in Subsection~\ref{rigidassump_subs},
the $C^1$-closure of~$\J_{\rig}(s)$ in $\J(\Pf)$ contains~$J_0$.\\

\begin{figure}
\begin{pspicture}(-1.1,-1.8)(10,1.25)
\psset{unit=.4cm}
\rput{45}(0,-4){\psellipse(5,-1.5)(2.5,1.5)
\psarc[linewidth=.05](5,-3.3){2}{60}{120}\psarc[linewidth=.05](5,0.3){2}{240}{300}
\pscircle[fillstyle=solid,fillcolor=gray](5,-4){1}\pscircle*(5,-3){.2}
\pscircle[fillstyle=solid,fillcolor=gray](6.83,.65){1}\pscircle*(6.44,-.28){.2}
\pscircle[fillstyle=solid,fillcolor=gray](3.17,.65){1}\pscircle*(3.56,-.28){.2}}
\rput(.2,-.9){$d_1$}\rput(3.1,2.3){$d_2$}\rput(7.8,-2.5){$d_3$}
\psarc(15,-1){3}{-60}{60}\psline(17,-1)(22,-1)\psline(16.8,-2)(21,-3)\psline(16.8,0)(21,1)
\rput(15.2,-3.5){$(1,0)$}\rput(22.4,1){$(0,d_1)$}
\rput(23.4,-1){$(0,d_2)$}\rput(22.4,-3){$(0,d_3)$}
\rput(33,-1){$d_1\!+\!d_2\!+\!d_3\!=\!d$}
\end{pspicture}
\caption{Domain of a Typical Element of $\M_1^3(\ka,d)$}
\label{m3_fig}
\end{figure}
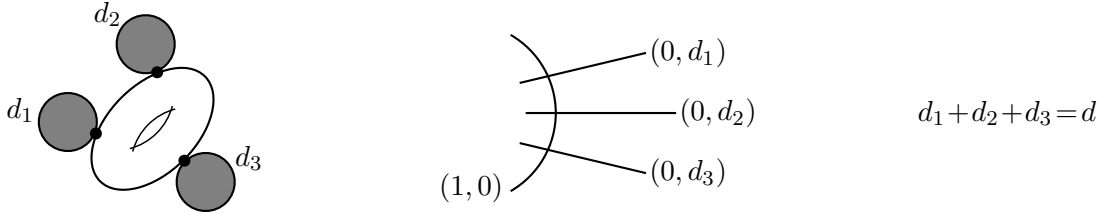

\noindent
From now on, we assume that $J\!\in\!{\cal J}_{\rig}(s)$
is an almost complex structure on $\Pf$ sufficiently close to~$J_0$.
For $g\!=\!0,1$, we~put
$$\S_g^d(Y;J)=\S_g^d(Y;J_Y) \quad\forall\,d\!\in\!\Z^+
\qquad\hbox{and}\qquad
\S_g(Y;J)=\bigsqcup_{d=1}^{\i}\S_g^d(Y;J).$$
If $\ka\!\in\!\S_g(Y;J)$, let $d_{\ka}$ denote the degree of $\ka$ in~$\Pf$.
If $\ka\!\in\!\S_0(Y;J)$  and $q$ is a positive integer, 
let $\M_1^q(\ka,d)$ be the subset of $\ov\M_1(\ka,d)$ consisting of stable maps 
$[\cC,u]$ such that $\cC$ is an elliptic curve~$E$ with $q$~rational 
components attached directly to~$E$ and $u|_E$ is constant.
Figure~\ref{m3_fig} shows the domain of a typical element of $\M_1^3(\ka,d)$,
from the points of view of symplectic topology and of algebraic geometry.
In the first diagram, each shaded disc represents a sphere;
the integer next to each rational component $\cC_i$ indicates the degree  of~$u|_{\cC_i}$.
In the second diagram, the components of $\cC$ are represented by curves,
and the pair of integers next to each component $\cC_i$ shows 
the genus of $\cC_i$ and the degree of~$u|_{\cC_i}$.
For stability reasons, the restriction of $u$ to each rational component
must be non-constant.
We denote by $\ov\M_1^q(\ka,d)$ the closure of $\M_1^q(\ka,d)$
in $\ov\M_1(\ka,d)$.
Note~that
\begin{equation}\label{spacesdim_e1}
\dim_{\C}\ov\M_1^q(\ka,d)=
\begin{cases}
2d,&\hbox{if}~q\!=\!0;\\
2d\!+\!1\!-\!q,&\hbox{if}~q\!\in\!\Z^+.
\end{cases}
\end{equation}
If $q\!\in\!\Z^+$, $\ov\M_1^q(\ka,d)$ is a smooth orbi-variety.
In contrast, $\ov\M_1^0(\ka,d)$ is a singular orbivariety, if $d\!>\!2$;
its structure is described in Subsection~\ref{conestr_subs}.\\

\noindent
For each $q\!\in\!\Z^+$, let $[q]\!=\!\{1,\ldots,q\}$.
If $\under{d}\!=\!(d_1,\ldots,d_q)$ is a $q$-tuple of positive integers
and $\ka\!\in\!\S_0(Y;J)$, we~put
\begin{equation}\label{wedge_e1}
\ov\M_0(\ka,\under{d})=
\big\{(b_1,\ldots,b_q)\!\in\!\prod_{i=1}^{i=q}\ov\M_{0,1}(\ka,d_i)\!:
\ev_0(b_i)\!=\!\ev_0(b_j)~\forall i,j\!\in\![q]\big\},
\end{equation}
where $\ev_0\!:\ov\M_{0,1}(\ka,d_i)\!\lra\!\ka$ is the evaluation
map corresponding to the marked point.
Let 
$$\ov\M_0^q(\ka,d)=\bigsqcup_{d_i>0,\sum d_i=d}\!\!\!\!\!\!\!
\ov\M_0\big(\ka,(d_1,\ldots,d_q)\big).$$
The spaces $\ov\M_0^q(\ka,d)$  are smooth orbi-varieties.
We note~that 
\begin{equation}\label{spacesdim_e0}
\dim_{\C}\ov\M_0^q(\ka,d)=2d\!+\!1\!-\!2q.
\end{equation}
By definition, 
\begin{equation}\label{mkdecomp_e}
\ov\M_1^q(\ka,d)=\big(\ov\cM_{1,q}\!\times\!\ov\M_0^q(\ka,d)\big)\big/S_q,
\end{equation}
where $\ov\cM_{1,q}$ is the moduli space of genus-one curves with $q$ marked points
and $S_q$ is the $q$th symmetric group.
The splitting \e_ref{mkdecomp_e} is illustrated in Figure~\ref{mkdecomp_fig}.
In this figure, we represent an entire space of stable maps
by the domain of a typical element of the space.
We shade the components of the domain on which the maps are non-constant.
The vertical bar in the last diagram indicates that the three marked points
are mapped to the same point in~$\ka$, as specified by~\e_ref{wedge_e1}.
Let 
$$\pi_P,\pi_B\!:\ov\cM_{1,q}\!\times\!\ov\M_0^q(\ka,d)
\lra\ov\cM_{1,q},\ov\M_0^q(\ka,d)$$
be the projection maps.\\

\begin{figure}
\begin{pspicture}(-1.1,-1.8)(10,1.25)
\psset{unit=.4cm}
\rput{45}(5,-4){\psellipse(5,-1.5)(2.5,1.5)
\psarc[linewidth=.05](5,-3.3){2}{60}{120}\psarc[linewidth=.05](5,0.3){2}{240}{300}
\pscircle[fillstyle=solid,fillcolor=gray](5,-4){1}\pscircle*(5,-3){.2}
\pscircle[fillstyle=solid,fillcolor=gray](6.83,.65){1}\pscircle*(6.44,-.28){.2}
\pscircle[fillstyle=solid,fillcolor=gray](3.17,.65){1}\pscircle*(3.56,-.28){.2}}
\rput(16,-1.5){$\approx$}
\rput(22,-1.5){$\ov\cM_{1,3}$}
\rput(24,-1.5){$\times$}
\psline(27,1)(27,-4)
\pscircle[fillstyle=solid,fillcolor=gray](28,1){1}\pscircle*(27,1){.2}
\pscircle[fillstyle=solid,fillcolor=gray](28,-1.5){1}\pscircle*(27,-1.5){.2}
\pscircle[fillstyle=solid,fillcolor=gray](28,-4){1}\pscircle*(27,-4){.2}
\end{pspicture}
\caption{The Decomposition~\e_ref{mkdecomp_e} for $\ov\M_1^3(\ka,d)$}
\label{mkdecomp_fig}
\end{figure}

\noindent
For each $\ka\!\!\in\!\S_0(Y;J)$, we denote by $N_Y\ka$ the normal bundle
of $\ka$ in~$Y$.
If  $q\!\in\!\Z^+$ and $[b]\!=\!([b_i])_{i\in[q]}$ is an element of
$\ov\M_0^q(\ka,d)$, let
$$\Ga(b;TY)\!=\!\big\{\xi\!=\!(\xi_i)_{i\in[q]}\!\in\!
\bigoplus_{i\in\![q]}\Ga(b_i;TY)\!:
\xi_i(y_0(b_i))\!=\!\xi_j(y_0(b_j))~\forall i,j\!\in\![q]\big\},$$
where $y_0(b_i)$ is the marked point of the component map~$b_i$.
Since $J_Y\!\in\!\J_{\rig}^0(Y)$, by the Index Theorem
the cokernel $H^1_J(b;TY)$ of the operator
\begin{equation}\label{wedgedbar_e}
\Ga(b;TY)\lra\Ga^{0,1}(b;TY;J)\!\equiv\!\bigoplus_{i\in\![q]}\Ga^{0,1}(b_i;TY;J),\qquad
D_{J,b}\big((\xi_i)_{i\in[q]}\big)=(D_{J,b_i}\xi_i)_{i\in[q]},
\end{equation}
is a vector space of dimension $2d\!-\!2$.
It is naturally isomorphic to the cokernel $H^1_J(b;N_Y\ka)$
of the operator
$$\Ga(b;N_Y\ka)\lra\Ga^{0,1}(b;N_Y\ka;J),\qquad
D_{J,b}^{\perp}\big((\xi_i)_{i\in[q]}\big)=(D_{J,b_i}^{\perp}\xi_i)_{i\in[q]},$$
induced by the operator $D_{J,b}$.
These cokernels induce a vector orbibundle over $\ov\M_0^q(\ka,d)$,
which will be denoted by~${\cal W}_{\ka,d}^{0,q}$.
If $q\!=\!1$, this bundle is the pullback by the forgetful map
$$\ti\pi\!: \ov\M_0^1(\ka,d)\!\equiv\!\ov\M_{0,1}(\ka,d) \lra\ov\M_0(\ka,d)$$
of the vector bundle defined in a similar way.
We denote this last vector bundle by~${\cal W}_{\ka,d}^0$.
We have
\begin{equation}\label{rankk_e}
\rk\W_{\ka,d}^{0,q}=2d-2 \quad\forall\, q\!\in\!\Z^+
\quad\hbox{and}\quad
\rk\W_{\ka,d}^0=2d-2.
\end{equation}\\

\noindent
From the decomposition~\e_ref{mkdecomp_e}, we see that the cokernel bundle for
the operators~$D_{J,b}$  over $\ov\M_1^q(\ka,d)$, for $q\!\in\!\Z^+$, is given~by
\begin{equation}\label{cokerkdecomp_e}
\W_{\ka,d}^{1,q}\approx
\big(\pi_P^*\E^*\!\otimes\!\pi_B^*\ev_0^*TY\oplus
\pi_B^*\W_{\ka,d}^{0,q}\big)/S_q,
\end{equation}
where $\E\!\lra\!\ov\cM_{1,q}$ is the Hodge line bundle and
$$\ev_0\!:\ov\M_0^q(\ka,d)\lra\ka$$
is the natural evaluation map, corresponding to the marked point common
to all factors.
We note~that
\begin{equation}\label{rankk_e1}
\rk\W_{\ka,d}^{1,q}=2d+1.
\end{equation}
On the other hand, similarly to the genus-zero case,
the cokernel $H^1_J(b;TY)$ of the operator~$D_{J,b}$ for
$$b\in\M_1^0(\ka,d)\subset\ov\M_1(\ka,d)$$
is naturally isomorphic to the cokernel $H^1_J(b;N_Y\ka)$ 
of the operator~$D_{J,b}^{\perp}$ induced by~$D_{J,b}$.
The cokernels $H^1_J(b;N_Y\ka)$  have the expected rank
for all $b\!\in\ov\M_1^0(\ka,d)$ and thus form a vector bundle over~$\ov\M_1^0(\ka,d)$,
which we denote by ${\cal W}_{\ka,d}^{1,0}$.
We have 
\begin{gather}\label{coker0decomp_e1}
\rk\W_{\ka,d}^{1,0}=2d \qquad\hbox{and}\\
\label{coker0decomp_e2}
\W_{\ka,d}^{1,0}\big|_{\ov\M_1^0(\ka,d)\cap\ov\M_1^q(\ka,d)}
\approx \big(
\pi_P^*\E^*\!\otimes\!\pi_B^*\ev_0^*N_Y\ka\oplus
\pi_B^*\W_{\ka,d}^{0,q}\big)/S_q~~\forall\, q\!\in\!\Z^+.
\end{gather}
We are now ready to reformulate Propositions~\ref{g0_prp1}-\ref{cone_prp1}.

\begin{prp}
\label{g0_prp2}
If $d$ and $\L$ are as in Theorem~\ref{main_thm}, 
$s\!\in\!H^0(\Pf;\L)$ is a transverse section, and $Y\!=\!s^{-1}(0)$, 
there exists $\de(d)\!\in\!\R^+$ with the following property.
If $J\!\in\!\J_{\rig}(s)$ and $\|J\!-\!J_0\|_{C^1}\!\le\!\de(d)$, then
$$N_0(d)=\sum_{\ka\in\S_0(Y;J)}\!\!\!
\blr{e(\W_{\ka,d/d_{\ka}}^0),\big[\ov\M_0(\ka,d/d_{\ka})\big]},$$
where $\W_{\ka,d/d_{\ka}}^0\!\lra\!\ov\M_0(\ka,d/d_{\ka})$
is the cokernel bundle corresponding to the almost complex structure~$J$,
as above.
\end{prp}

\begin{prp}
\label{g1_prp2}
If $d$, $\L$, $s$, and $Y$ are as in Proposition~\ref{g0_prp2}, 
there exists $\de(d)\!\in\!\R^+$ with the following property.
If $J\!\in\!\J_{\reg}(s)$ and $\|J\!-\!J_0\|_{C^1}\!\le\!\de(d)$, then
\begin{equation*}\begin{split}
&N_1(d)={\sum_{\ka\in\S_1(Y;J)}} {^{\pm}\big|\M_1^0(\ka,d/d_{\ka})\big|}\\
&\qquad\qquad\qquad
+\sum_{\ka\in\S_0(Y;J)}\!\!\!
\Big(\blr{e(\W_{\ka,d/d_{\ka}}^{1,0}),\big[\ov\M_1^0(\ka,d/d_{\ka})\big]}
+\frac{d/d_{\ka}}{12}\blr{e(\W_{\ka,d/d_{\ka}}^0),\big[\ov\M_0(\ka,d/d_{\ka})\big]}\Big),
\end{split}\end{equation*}
where $\W_{\ka,d/d_{\ka}}^0\!\lra\!\ov\M_0(\ka,d/d_{\ka})$
and $\W_{\ka,d/d_{\ka}}^{1,0}\!\lra\!\ov\M_1^0(\ka,d/d_{\ka})$
are the cokernel bundles corresponding to the almost complex structure~$J$,
as above.
\end{prp}

\begin{prp}
\label{cone_prp2}
If $d$, $\L$, $s$, and $Y$ are as in Proposition~\ref{g0_prp2}
and $\V_1^d\!\lra\X_1(\Pf,d)$ is the cone corresponding
to the line bundle $\L\!\lra\!\Pf$ with its standard connection,
there exists $\de(d)\!\in\!\R^+$ with the following properties.
If $\|J\!-\!J_0\|_{C^1}\!\le\!\de(d)$, then the moduli space 
$\ov\M_1^0(\Pf,d;J)$ carries a rational fundamental class 
of dimension~$10d$, the euler class of the~cone 
$$\V_1^d\lra\ov\M_1^0(\Pf,d;J)$$
is a well-defined element of $H^{10d}(\ov\M_1^0(\Pf,d;J);\Q)$, and
$$\blr{e(\V_1^d),\big[\ov\M_1^0(\Pf,d;J)\big]}
=\blr{e(\V_1^d),\big[\ov\M_1^0(\Pf,d)\big]}.$$
If in addition $J\!\in\!{\cal J}_{\reg}(s)$, 
\begin{equation*}\begin{split}
&\blr{e(\V_1^d),\big[\ov\M_1^0(\Pf,d;J)\big]}
={\sum_{\ka\in\S_1(Y;J)}} {^{\pm}\big|\M_1^0(\ka,d/d_{\ka})\big|}\\
&\qquad\qquad
+\!\!\sum_{\ka\in\S_0(Y;J)}\!\!\!
\Big(\blr{e(\W_{\ka,d/d_{\ka}}^{1,0}),\big[\ov\M_1^0(\ka,d/d_{\ka})\big]}
+\frac{d/d_{\ka}-1}{12}\blr{e(\W_{\ka,d/d_{\ka}}^0),\big[\ov\M_0(\ka,d/d_{\ka})\big]}\Big),
\end{split}\end{equation*}
where $\W_{\ka,d/d_{\ka}}^0\!\lra\!\ov\M_0(\ka,d/d_{\ka})$
and $\W_{\ka,d/d_{\ka}}^{1,0}\!\lra\!\ov\M_1^0(\ka,d/d_{\ka})$
are as in Proposition~\ref{g1_prp2}.
\end{prp}

\noindent
In the last two propositions, 
the moduli space consists $\M_1^0(\ka,d/d_{\ka})$, for $\ka\!\in\!\S_1(Y;J)$,
contains only one element:
the equivalence class of the degree-$d/d_{\ka}$ cover of the elliptic curve $\ka$
by an elliptic curve.
Since the order of the automorphism group of such a cover is~$d/d_{\ka}$,
$$^{\pm}\big|\M_1^0(\ka,d/d_{\ka})\big|=\pm\frac{1}{d/d_{\ka}}.$$
The sign is determined by viewing the zero-dimensional suborbifold
$\M_1^0(\ka,d/d_{\ka})$ of $\X_1(Y,d)$ as a transverse zero 
of the section~$\bar{\partial}_J$.
This sign is the same as the sign of~$\ka$ as an element of 
the set $\S_1^{d_{\ka}}(Y;J)$.
In particular, 
\begin{equation}\label{g1reduc_e}
{\sum_{\ka\in\S_1(Y;J)}} {^{\pm}\big|\M_1^0(\ka,d/d_{\ka})\big|}
=\sum_{\si|d}\frac{n_1(d/\si)}{\si},
\end{equation}
where $n_1(\cdot)$ is as in Subsection~\ref{rigidassump_subs}.\\

\noindent
If $\ka\!\in\!\S_0(Y;J)$, the orientations of the vector bundles
$$\W_{\ka,d/d_{\ka}}^0\!\lra\!\ov\M_0(\ka,d/d_{\ka})
\qquad\hbox{and}\qquad
\W_{\ka,d/d_{\ka}}^{1,0}\!\lra\!\ov\M_1^0(\ka,d/d_{\ka})$$
are determined by the linearizations of the sections $\bpar_J$ 
over $\X_0(Y,d)$ and~$\X_1(Y,d)$. 
According to~\cite{IP}, by a spectral-flow argument it can be shown that
\begin{alignat}{1}
\label{flow_e1}
\blr{e(\W_{\ka,\si}^0),\big[\ov\M_0(\ka,\si)\big]}
&=\pm\blr{e\big(R^1\pi_{0*}^{\si}\ev_0^{\si*}(
\O_{\ka}(-1)\!\oplus\!\O_{\ka}(-1))\big),\big[\ov\M_0(\ka,\si)\big]},\\
\label{flow_e2}
\blr{e(\W_{\ka,\si}^{1,0}),\big[\ov\M_1^0(\ka,\si)\big]}
&=\pm\blr{e\big(R^1\pi_{1*}^{\si}\ev_1^{\si*}
(\O_{\ka}(-1)\!\oplus\!\O_{\ka}(-1))\big),\big[\ov\M_1^0(\ka,\si)\big]},
\end{alignat}
where $\si\!=\!d/d_{\ka}$ and the sign agrees with the sign of $\ka$
as an element of $\S_0^{d_{\ka}}(Y;J)$.
By localization, 
\begin{equation}\label{local_e1}
\blr{e\big(R^1\pi_{0*}^{\si}\ev_0^{\si*}
(\O_{\ka}(-1)\!\oplus\!\O_{\ka}(-1))\big),\big[\ov\M_0(\ka,\si)\big]}
=\frac{1}{\si^3};
\end{equation}
see Section 27.5 of~\cite{H}.
Using the desingularization of $\ov\M_1^0(\ka,\si)$ constructed in~\cite{VZ}, 
it should be possible to show~that
\begin{equation}\label{local_e2}
\blr{e\big(R^1\pi_{1*}^{\si}\ev_1^{\si*}
(\O_{\ka}(-1)\!\oplus\!\O_{\ka}(-1))\big),\big[\ov\M_1^0(\ka,\si)\big]}
=\frac{1}{12}\,\frac{\si-1}{\si^2}.
\end{equation}
Propositions~\ref{g0_prp1}-\ref{cone_prp1} follow from 
Propositions~\ref{g0_prp2}-\ref{cone_prp2} via \e_ref{g1reduc_e}-\e_ref{local_e2}.\\

\noindent
Since Theorem~\ref{main_thm} follows immediately from
Propositions~\ref{g0_prp2}-\ref{cone_prp2},
we do not need to deduce Propositions~\ref{g0_prp1}-\ref{cone_prp1} 
from Propositions~\ref{g0_prp2}-\ref{cone_prp2}.
We prove Propositions~\ref{g1_prp2} and~\ref{cone_prp2} in Sections~\ref{dbar_sec}
and~\ref{cone_sec}, respectively; see also Propositions~\ref{maincontr_prp}
and~\ref{bdcontr_prp}.
The proof of Proposition~\ref{g0_prp2} is very similar to the proof
of Proposition~\ref{g1_prp2}, but simpler, and we omit~it.

\subsection{Summary of the Proof of Proposition~\ref{g1_prp2}}
\label{summary_subs}

\noindent
A key notion in our argument, which is also used in the proof of Proposition~\ref{cone_prp2}, 
is Definition~\ref{contr_dfn} below.
For its purposes, we will call either of the two tuples
$\V_g^d(\bpar;J)$ and $\V_1^d(s;J)$,
defined in Subsection~\ref{review_subs}, a {\tt generalized vector bundle}.
The first tuple involves an infinite-rank bundle over 
an infinite-dimensional space;
the second one involves finite-dimensional objects, albeit non-smooth ones.
Nevertheless, both are generalizations of a rank-$n$ vector bundle~$\F$ 
over an $n$-dimensional complex compact manifold~$\X$, 
with a choice of a section~$\vph$ and of an appropriate subset~${\cal A}(\vph)$
of $\Ga(\X;\F)$ of second category.
Such a collection of data can also be considered to be a generalized vector bundle.

\begin{dfn}
\label{contr_dfn}
Suppose $\V\!=\!\big(\X,\F,\pi;\vph,{\cal A}(\vph)\big)$ 
is a generalized vector bundle. 
Subset $\cZ$ of~$\vph^{-1}(0)$ is a {\tt regular set for $\V$} if there exists
${\cC}_{\cZ}(\V)\!\in\!\Q$ and a dense open subset 
${\cal A}_{\cZ}(\vph)$ of ${\cal A}(\vph)$ with the following properties.
For every $\nu\!\in\!{\cal A}_{\cZ}(\vph)$,\\
${}\quad$ (a) there exists $\ep_{\nu}\!\in\!\R^+$ such that $t\nu\!\in\!{\cal A}(\vph)$
for all $t\!\in\!(0,\ep_{\nu})$;\\
${}\quad$ (b) there exist a compact subset $K_{\nu}\!\subset\!\cZ$,
open neighborhood $U_{\nu}(K)$ of $K$ in $\X$ for each\\ 
${}\qquad~~$ compact subset $K\!\subset\!\cZ$,
and $\ep_{\nu}(U)\!\in\!(0,\ep_{\nu})$ for each open subset $U$ of $\X$
such that
$$ ^{\pm}\!\big|\{\vph\!+\!t\nu\}^{-1}\cap U\big|=\cC_{\cZ}(\V)
\qquad\hbox{if}~~t\!\in\!\big(0,\ep_{\nu}(U)\big)
~\hbox{and}~K_{\nu}\!\subset\!K\!\subset\!U\!\subset\!U_{\nu}(K).$$\\
\end{dfn}

\noindent
Every connected component of $\vph^{-1}(0)$ is regular.
However, a  regular subset of $\vph^{-1}(0)$ need not be closed.
For example, if $\vph$ is a holomorphic section of 
a rank-$k$ algebraic vector bundle~${\frak F}$ over 
a $k$-dimensional compact algebraic variety~${\frak X}$,
every Zariski open subset of~$\vph^{-1}(0)$ is regular.
The sections $\bpar_J$ and $s_1^d$ that play a central role
in this paper are in a sense {\tt generalized holomorphic sections}.\\

\noindent
If $\cZ$ is a regular set for the generalized vector bundle~$\V$,
we will call the number $\cC_{\cZ}(\V)$ 
the {\tt $\vph$-contribution of $\cZ$ to the euler class of~$\V$}.
Note that if $\vph^{-1}(0)\!=\!\sqcup_{i\in I}\cZ_i$ 
is a partition of $\vph^{-1}(0)$ into regular sets,
the euler class of $\V$, or its Poincare dual, is the sum of $\vph$-contributions:
\begin{equation}\label{contr_e}
e(\V)=\sum_{i\in I}\cC_{\cZ_i}(\eta).
\end{equation}
We prove Theorem~\ref{main_thm} by expressing each of the three terms appearing 
in~\e_ref{main_thm_e} in the form~\e_ref{contr_e}
and show that we end up with the same terms on the two sides of~\e_ref{main_thm_e}.\\

\noindent
If $d$, $s$, and $Y$ are as in the previous subsection and $J\!\in\!{\J}_{\rig}(s)$,
\begin{equation}\label{cydecomp_e}
\ov\M_1(Y,d;J)
=\bigsqcup_{\ka\in\S_0(Y;J)}\!\!\!\!\!\ov\M_1(\ka,d/d_{\ka})
~\sqcup  
\bigsqcup_{\ka\in\S_1(Y;J)}\!\!\!\!\!\!\M_1^0(\ka,d/d_{\ka}).
\end{equation}
For any $\ka\!\in\!\S_0(Y;J)$, $\si\!\in\!\Z^+$, and 
subset $\vr$ of $\bar\Z^+\!\equiv\!\Z^+\!\cup\!\{0\}$, let
$$\M_1^{\vr}(\ka,\si)= \bigcap_{q\in\vr}\ov\M_1^q(\ka,\si)
~-\!
\bigcup_{q\in\bar\Z^+\!-\vr}\!\!\!\!\ov\M_1^q(\ka,\si).$$

\begin{prp}
\label{maincontr_prp}
If $d$, $\L$, $s$, and $Y$ are as in Proposition~\ref{g0_prp2}, 
$J\!\in\!\J_{\rig}(s)$ is sufficiently close to~$J_0$,
and $\ka\!\in\!\S_1(Y;J)$, then
$$\cC_{\M_1^0(\ka,d/d_{\ka})}\big(\V_1^d(\bpar;J)\big)
=\, ^{\pm}\big|\M_1^0(\ka,d/d_{\ka})\big|.$$
If $\ka\!\in\!\S_0(Y;J)$, 
$$\cC_{\M_1^{\{0\}}(\ka,d/d_{\ka})}
\big(\V_1^d(\bpar;J)\big)
=\blr{e(\W_{\ka,d/d_{\ka}}^{1,0}),\big[\ov\M_1^0(\ka,d/d_{\ka})\big]}.$$
\end{prp}

\begin{prp}
\label{bdcontr_prp}
If $d$, $\L$, $s$, $Y$, and $J$ are as in Proposition~\ref{maincontr_prp},
$\ka\!\in\!\S_0(Y;J)$, and $\vr$ is a subset $\bar\Z^+$ different from $\{0\}$, then
$$\cC_{\M_1^{\vr}(\ka,d/d_{\ka})}\big(\V_1^d(\bar{\partial};J)\big)
=\begin{cases}
\frac{d/d_{\ka}}{12}
\blr{e(\W_{\ka,d/d_{\ka}}^0),\big[\ov\M_0(\ka,d/d_{\ka})\big]},
&\hbox{if}~\vr\!=\!\{1\};\\
0,&\hbox{if}~\vr\!\neq\!\{1\}.
\end{cases}$$\\\end{prp}

\noindent
One consequence of Propositions~\ref{maincontr_prp} and~\ref{bdcontr_prp}
is that most boundary strata of the moduli space $\ov\M_1(Y,d;J)$
do not contribute to the number~$N_1(d)$.
In fact, we will show that only the strata 
$\M_1^0(Y,d;\ka)$  and $\M_1^1(Y,d;\ka)$ contribute to the number~$N_1(d)$.\\

\noindent
We now outline the proofs of Propositions~\ref{maincontr_prp} and~\ref{bdcontr_prp}.
Let 
$$\nu\in\Ga\big(\X_1(\P,d);\Ga_1^{0,1}(\Pf,d;J)\big)$$ 
be a small generic multisection such that
$$\nu\in \Ga\big(\X_s;\Ga_1^{0,1}(Y,d;J)\big)$$
for a small neighborhood $\X_s$ of
$\ov\M_1(Y,d;J)$ in $\X_1(\Pf,d)$ and vanishes outside of~$\U_s$.
By definition, $N_1(d)$ is the number of elements 
$\exp_u\!\xi\!\in\!\X_1(\Pf,d)$ such that $(u,\xi)$ solves the system
\begin{equation}\label{setup_e}
\begin{cases}
\bpar_J\exp_u\!\xi+\nu(\exp_u\!\xi)=0;\\
s\circ\exp_u\xi=0;
\end{cases}
\qquad u\!\in\!\ov\M_1(\Pf,d;J),~\xi\!\in\!T_u\X_1(\Pf,d).
\end{equation} 
Note that 
$$\bpar_{\na,\exp_u\!\xi}s_1^d(\exp_u\!\xi)=0$$ 
if $(u,\xi)$ solves the first equation, due to our assumptions on~$\nu$.
If $u\!\in\!\M_1^0(\Pf,d;J)$ and $\nu$ is sufficiently small,
the first equation has a unique small solution $\xi_{\nu}(u)$ in~$\Ga_+(u)$,
the orthogonal complement of $T_u\ov\M_1(\Pf,d;J)$ in~$T_u\X_1(\Pf,d)$.
Plugging this solution into the second equation, we obtain 
\begin{equation}\label{main_gen_e}
0=s\circ\exp_u\xi_{\nu}=s_1^d(u)+\pi_{TY}^{\perp}\xi_{\nu}(u)\in\V_1^d,
\end{equation}
where $\pi_{TY}^{\perp}$ is the projection map $T\Pf\!\lra\!T\Pf/TY$, 
defined on a neighborhood of $Y$ in~$\Pf$.
Since all solutions of the system~\e_ref{setup_e} are transverse,
so are the solutions of~\e_ref{main_gen_e}.
Thus, the zeros of a generic perturbation $\nu$ of the section $\bar{\partial}_J$
that lie close to $\M_1^0(Y,d;J)$ correspond to the zeros
of a perturbation of the section $s_1^d$ that lie close to~$\M_1^0(Y,d;J)$.
In Subsection~\ref{maincontr_subs}, we show that the number of these zeros
that lie near each component $\M_1^0(\ka,d/d_{\ka})$ of $\M_1^0(Y,d;J)$
is the euler class
of the bundle ${\cal W}_{\ka,d/d_{\ka}}^{1,0}$ over $\ov\M_1^0(\ka,d/d_{\ka})$.\\

\noindent
We next look for solutions near $\M_1^{\{1\}}(Y,d;J)$, 
i.e.~we assume that $u\!\in\!\M_1^{\{1\}}(Y,d;J)$.
Note that 
\begin{equation}\label{m1decomp_e}
\ov\M_1^1(Y,d;J)\approx\ov\cM_{1,1}\times\ov\M_{0,1}(Y,d;J).
\end{equation}
We denote the projection maps onto $\ov\cM_{1,1}$ and $\ov\M_{0,1}(Y,d;J)$
by $\pi_P$ and $\pi_B$, respectively. 
Let 
$$\ti\pi_B\!:\ov\M_1^1(Y,d;J)\lra\ov\M_0(Y,d;J)$$
be the composition of $\pi_B$
with the forgetful map $\ov\M_{0,1}(Y,d;J)\!\lra\!\ov\M_0(Y,d;J)$.
The bundle $T\X_1(Y,d)$ contains the line subbundle
${\cal L}\!\equiv\!\pi_P^*L_{P,1}\!\otimes\!\pi_B^*L_0$, where 
$$L_{P,1}\lra\ov\cM_{1,1}
\qquad\hbox{and}\qquad L_0\!\lra\!\ov\M_{0,1}(Y,d;J)$$
are the universal tangent line bundles at the marked points.
If $u\!\in\!\M_1^1(Y,d;J)$ and $\ups\!\in\!{\cal L}_u$ is small,
we denote by $u_{\ups}$ the element $\exp_u\!\ups$ of~$\X_1(Y,d)$.
Let 
$$\ev_P\!:\ov\M_1^1(Y,d;J)\lra Y$$ 
be the composition of $\pi_B$ with the evaluation map at the marked point.
This map sends an element $[\cC,u]$ of $\ov\M_1^1(Y,d;J)$ to 
the value of $u$ on the principal component of~$\cC$.\\

\noindent
In this case, we work with the analogue of~\e_ref{setup_e} intrinsic to~$Y$,
i.e.~we look for solutions of the equation
\begin{equation}\label{setup2_e}
\bpar_J\exp_{u_{\ups}}\!\xi+\nu(\exp_{u_{\ups}}\!\xi)=0
\qquad u\!\in\!\ov\M_1^{\{1\}}(Y,d;J),~
\xi\!\in\!\Ga(\ups;TY)\!\equiv\!\Ga(u_{\ups}^*TY).
\end{equation}
This equation usually does not have a small solution in~$\xi$ for a fixed~$u_{\ups}$, 
as there is an obstruction bundle
$$\Ga_-^{0,1}(u;TY;J)= \ti\pi_B^*H^1(u_B^*TY)\oplus
\pi_P^*\E^*\!\otimes\!\ev_P^*TY\subset\Ga^{0,1}_1(Y,d;J),$$
where $u_B$ is the restriction of~$u$ to the bubble components.
Taking the projections $(\pi^{0,1}_{-,B}\!\oplus\!\pi^{0,1}_{-,P})$ and 
$\pi^{0,1}_+$ of~\e_ref{setup2_e} onto $\Ga^{0,1}_-(u;TY;J)$ and 
its complement $\Ga^{0,1}_+(u;TY;J)$ in $\Ga^{0,1}_1(Y,d;J)$, respectively,
we obtain
\begin{equation}\label{setup2_e2}
\begin{cases}
\pi^{0,1}_+\bpar\exp_{u_{\ups}}\!\xi+\pi^{0,1}_+\nu(\exp_{u_{\ups}}\!\xi)=0
\in\Ga^{0,1}_+(u;TY;J);\\
\pi^{0,1}_{-,B}\bpar\exp_{u_{\ups}}\!\xi+\pi^{0,1}_{-,1}\nu(\exp_{u_{\ups}}\!\xi)=0
\in\ti\pi_B^*H^1(u_B^*TY);\\
\pi^{0,1}_{-,P}\bpar\exp_{u_{\ups}}\!\xi+\pi^{0,1}_{-,0}\nu(\exp_{u_{\ups}}\!\xi)=0
\in\pi_P^*\E^*\!\otimes\!\ev_P^*TY.
\end{cases}
\end{equation}
If $\nu$ and $\ups$ are sufficiently small,
the first equation has a unique small solution $\xi_{\nu}(u,\ups)$ in $\Ga_+(u;TY)$.
With appropriate choice of neighborhood charts and of the perturbation~$\nu$,
$\xi_{\nu}(u,\ups)$ depends only on~$u_B$, and the system~\e_ref{setup2_e2} is equivalent~to
\begin{equation}\label{setup2_e3}
\pi^{0,1}_{-,P}\bpar\exp_{u_{\ups}}\!\xi+\pi^{0,1}_{-,P}\nu(\exp_{u_{\ups}}\!\xi)=0
\in\pi_P^*\E^*\!\otimes\!\ev_P^*TY\subset\Ga^{0,1}_1(Y,d;J),\quad
u_B\!\in\!\cZ_0,~\ups\!\in\!{\cal L},
\end{equation}
where $\cZ_0$ is the zero set of a section of the first component of
the bundle~$\Ga_-^{0,1}(\cdot,TY)$ over $\ov\M_0(Y,d;J)$.
In particular, $^{\pm}\!|\cZ_0|\!=\!N_0(d)$.\\

\noindent
Equation \e_ref{setup2_e3} is equivalent to
\begin{equation}\label{setup2_e4}
{\cal D}_u\ups+\pi^{0,1}_{-,P}\nu(u)=0\in \pi_P^*\E^*\!\otimes\!T_{\ev_P(u)}Y,
\quad u_B\!\in\!{\cal Z}_0,~\ups\!\in\!{\cal L},
\end{equation}
where ${\cal D}_u\!\in\!\hbox{Hom}(L_0,T_{\ev_0(u_B)}Y)$.
The image of ${\cal D}_u$ in $T_{\ev_0(u_B)}Y$ is precisely the tangent line
at $\ev_0(u_B)$ to the rational curve~$\Im{u_B}$,
as long as the differential of the map~$u_B$ does not vanish
at the marked point.
Thus, for each $u_B\!\in\!{\cal Z}_0$, the number of solutions of~\e_ref{setup2_e4}
is the number of times $\pi^{0,1}_{-,P}\nu(u)$ lies in $\E^*\!\otimes\!T_{\ev_P(u)}\Im u_B$.
We conclude that
\begin{equation}\label{bdcontrcomp_e}\begin{split}
\cC_{\M_1^{\{1\}}(Y,d;J)}\big(\V_1^d(\bpar;J)\big)
&=\sum_{u_B\in\cZ_0}\!\blr{c\big(\pi_P^*\E^*\!\otimes\!\ev_P^*TY\big)
c\big(\pi_P^*\E^*\!\otimes\!\ev_P^*T\Im u_B\big)^{-1},
\big[\ti\pi_B^{-1}(u_B)\big]}\\
&=\sum_{u_B\in\cZ_0}\!
\blr{c_1(\E^*)\big(c_1(TY)\!-\!c_1(T\Im u_B)\big),[\ov\cM_{1,1}]\!\times\![\Bbb{P}^1]}\\
&=-\frac{1}{24}\cdot(0\!-\!2)\cdot 
\sum_{\ka\in\S_0(Y;J)}\!\!(d/d_{\ka})\cdot\,
^{\pm}\!\big|\cZ_0\cap\ov\M_0(\ka,d/d_{\ka})\big|\\
&=\sum_{\ka\in\S_0(Y;J)}\!\!\frac{d/d_{\ka}}{12}
\blr{e(\W_{\ka,d/d_{\ka}}^0),\big[\ov\M_0(\ka,d/d_{\ka})\big]},
\end{split}\end{equation}
as claimed in Proposition~\ref{bdcontr_prp}.\\

\noindent
We analyze the contribution to the number $N_1(d)$ from 
the complement of $\M_1^{\{0\}}(\Pf,d;J)$ and
$\M_1^{\{1\}}(\Pf,d;J)$  in $\ov\M_1(Y,d;J)$
in a similar way, but we encounter one of two key differences.
If $\vr\!=\!\{0,1\}$ and $u\!\in\!\M_1^{\vr}(Y,d;J)$, ${\cal D}_u\!=\!0$.
Thus, equation~\e_ref{setup2_e4} has no solutions near $\M_1^{\vr}(Y,d;J)$
if~$\nu$ is generic.
On the other hand, if $\vr$ is any other subset of $\bar\Z^+$ containing~$0$
and at least one other element,
the analogue of the set~$\cZ_0$ is empty for dimensional reasons.
Thus, 
$$\cC_{\M_1^{\vr}(Y,d;J)}\big(\V_1^d(\bpar;J)\big)=0
\qquad\hbox{if}\quad \{0\}\!\subsetneq\!\vr\!\subset\!\bar\Z^+,$$
as claimed.\\

\noindent
The computation of the contribution from $\M_1^0(Y,d;J)$ 
to the number $N_1(d)$ can also be carried out in~$Y$, instead of~$\Pf$.
However, the presented version of the computation is meant to indicate 
why the cone $\V_1^d$ should enter into the Gromov-Witten theory of~$Y$.\\

\noindent
We supply more details of the proof of 
Propositions~\ref{maincontr_prp} and~\ref{bdcontr_prp} in Section~\ref{dbar_sec}. 
In particular, in order to use the gluing and obstruction-bundle setup
described in~\cite{gluing}, 
we stratify the moduli spaces that appear in the statements
of Propositions~\ref{maincontr_prp} and~\ref{bdcontr_prp} 
according to the bubble type, or the dual graph, of stable maps.
The notion of contribution to the euler class used in this paper is a direct adaptation,
to the orbifold and multisection setting of~\cite{FuO} and~\cite{LT},
of the analogous notion used in~\cite{g2n2and3} and~\cite{genus0pr}.
However, in the present case, we can get by with far less detailed understanding 
of the behavior of the bundle sections involved.

\subsection{Notation: Genus-Zero Maps}
\label{notation0_subs}

\noindent
We now summarize our notation for bubble maps from genus-zero Riemann surfaces,
with one marked point, and for related objects.
For more details on the notation described below, the reader is referred
to Sections~ in~\cite{gluing}.\\

\noindent
In general, moduli spaces of stable maps can stratified by the dual graph.
However, in the present situation, it is more convenient to make use
of {\it linearly ordered sets}:

\begin{dfn}
\label{index_set_dfn1}
(1) A finite nonempty partially ordered set $I$ is a {\tt linearly ordered set} if 
for all \hbox{$i_1,i_2,h\!\in\!I$} such that $i_1,i_2\!<\!h$, 
either $i_1\!\le\!i_2$ \hbox{or $i_2\!\le\!i_1$.}\\
(2) A linearly ordered set $I$ is a {\tt rooted tree} if
$I$ has a unique minimal element, 
i.e.~there exists \hbox{$\hat{0}\!\in\!I$} such that $\hat{0}\!\le\!i$ 
for {all $i\!\in\!I$}.
\end{dfn}

\noindent
We use rooted trees to stratify the moduli space $\ov\M_{0,1}(\Pf,d;J)$
of degree-$d$ $J$-holomorphic maps from genus-zero Riemann surfaces with 
one marked point to~$\Pf$.\\

\noindent
If $I$ is a linearly ordered set, let $\hat{I}$ be 
the subset of the non-minimal elements of~$I$.
For every $h\!\in\!\hat{I}$,  denote by $\io_h\!\in\!I$
the largest element of $I$ which is smaller than~$h$, i.e.
$$\io_h=\max\big\{i\!\in\!I:i\!<\!h\big\}.$$
A {\tt genus-zero $\Pf$-valued bubble map} is a tuple $b\!=\!(I;x,u)$,
where $I$ is a rooted tree, and
$$x\!:\hat{I}\!\lra\!\C\!=\!S^2\!-\!\{\i\} \hbox{~~~and~~~} 
u\!:I\!\lra\!C^{\i}(S^2;\Pf)$$
are maps such that $u_h(\i)\!=\!u_{\io_h}(x_h)$ for all $h\!\in\!\hat{I}$.
Such a tuple describes a Riemann surface $\Si_b$ and 
a continuous map $u_b\!:\Si_b\!\lra\!\Pf$.
The irreducible components $\Si_{b,i}$ of $\Si_b$ are indexed by the set~$I$ 
and $u_b|_{\Si_{b,i}}\!=\!u_i$.
The Riemann surface $\Si_b$ carries a marked point,
i.e.~the point $(\hat{0},\i)\!\in\!\Si_{b,\hat{0}}$,
if $\hat{0}$ is the minimal element of~$I$.
The general structure of genus-zero bubble maps is described
by tuples $\T\!=\!(I;\under{d})$, where 
$\under{d}\!:I\!\lra\!\Z$ is a map specifying the degree 
of $u_b|_{\Si_{b,i}}$, if $b$ is a bubble map of type~$\T$.
We call such tuples {\tt bubble types}.\\

\noindent
If $\T$ is a bubble type as above, 
let $\U_{\T}(\Pf;J)$ be the subset of $\ov\M_{0,1}(\Pf,d;J)$ 
consisting of stable maps $[\cC,y_1,u]$ such that
$$[\cC,y_1,u]=\big[(\Si_b,(\hat{0},\i)),u_b\big],$$
for some bubble map $b$ of type~$\T$.
Subsection~2.5 of~\cite{gluing} describes a space $\U_{\T}^{(0)}(X;J)$ 
of {\it balanced} stable maps, not of equivalence classes of such maps,
such that
$$\U_{\T}(X;J)=\U_{\T}^{(0)}(X;J)\big/ \hbox{Aut}(\T)\!\propto\!(S^1)^I,$$
for a natural action of $\hbox{Aut}(\T)$ on $(S^1)^I$.
This space is convenient to use in gluing constructions.

\subsection{Notation: Genus-One Maps}
\label{notation1_subs}

\noindent
We next set up analogous notation for genus-one stable maps;
see Subsection~\ref{g1comp-notation1_subs} in~\cite{g1comp} for more details.
In this case, we also need to specify the structure of the principal component.
Thus, we index the strata of~$\ov\M_1(\Pf,d;J)$ by
{\it enhanced linearly ordered sets}:

\begin{dfn}
\label{index_set_dfn2}
An {\tt enhanced linearly ordered set} is a pair $(I,\aleph)$,
where $I$ is a linearly ordered set, $\aleph$ is a subset of $I_0\!\times\!I_0$,
and $I_0$ is the subset of minimal elements of~$I$,
such that if $|I_0|\!>\!1$, 
$$\aleph=\big\{(i_1,i_2),(i_2,i_3),\ldots,(i_{n-1},i_n),(i_n,i_1)\big\}$$
for some bijection $i\!:\{1,\ldots,n\}\!\lra\!I_0$.
\end{dfn}

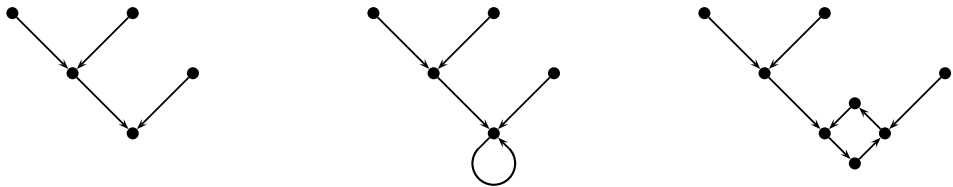
\begin{figure}
\begin{pspicture}(-1.1,-2)(10,1)
\psset{unit=.4cm}
\pscircle*(6,-3){.2}
\pscircle*(4,-1){.2}\psline[linewidth=.06]{->}(4.14,-1.14)(5.86,-2.86)
\pscircle*(8,-1){.2}\psline[linewidth=.06]{->}(7.86,-1.14)(6.14,-2.86)
\pscircle*(2,1){.2}\psline[linewidth=.06]{->}(2.14,.86)(3.86,-.86)
\pscircle*(6,1){.2}\psline[linewidth=.06]{->}(5.86,.86)(4.14,-.86)
\pscircle*(18,-3){.2}\psline[linewidth=.06](17.86,-3.14)(17.5,-3.5)
\psarc(18,-4){.71}{135}{45}\psline[linewidth=.06]{->}(18.5,-3.5)(18.14,-3.14)
\pscircle*(16,-1){.2}\psline[linewidth=.06]{->}(16.14,-1.14)(17.86,-2.86)
\pscircle*(20,-1){.2}\psline[linewidth=.06]{->}(19.86,-1.14)(18.14,-2.86)
\pscircle*(14,1){.2}\psline[linewidth=.06]{->}(14.14,.86)(15.86,-.86)
\pscircle*(18,1){.2}\psline[linewidth=.06]{->}(17.86,.86)(16.14,-.86)
\pscircle*(30,-2){.2}\pscircle*(30,-4){.2}\pscircle*(29,-3){.2}\pscircle*(31,-3){.2}
\psline[linewidth=.06]{->}(29.86,-2.14)(29.14,-2.86)
\psline[linewidth=.06]{->}(29.14,-3.14)(29.86,-3.86)
\psline[linewidth=.06]{->}(30.14,-3.86)(30.86,-3.14)
\psline[linewidth=.06]{->}(30.86,-2.86)(30.14,-2.14)
\pscircle*(27,-1){.2}\psline[linewidth=.06]{->}(27.14,-1.14)(28.86,-2.86)
\pscircle*(33,-1){.2}\psline[linewidth=.06]{->}(32.86,-1.14)(31.14,-2.86)
\pscircle*(25,1){.2}\psline[linewidth=.06]{->}(25.14,.86)(26.86,-.86)
\pscircle*(29,1){.2}\psline[linewidth=.06]{->}(28.86,.86)(27.14,-.86)
\end{pspicture}
\caption{Some Enhanced Linearly Ordered Sets}
\label{index_set_fig}
\end{figure}

\noindent
An enhanced linearly ordered set can be represented by an oriented connected graph.
In Figure~\ref{index_set_fig}, the dots denote the elements of~$I$.
The arrows outside the loop, if there are any, 
specify the partial ordering of the linearly ordered set~$I$.
In fact, every directed edge outside of the loop
connects a non-minimal element $h$ of $I$ with~$\io_h$.
Inside of the loop, there is a directed edge from $i_1$ to $i_2$
if and only if $(i_1,i_2)\!\in\!\aleph$.\\

\noindent
The subset $\aleph$ of $I_0\!\times\!I_0$ will be used to describe
the structure of the principal curve of the domain of stable maps in 
a stratum of the moduli space~$\ov\M_1(\Pf,d;J)$.
If $\aleph\!=\!\eset$, and thus $|I_0|\!=\!1$,
the corresponding principal curve $\Si_{\aleph}$ 
is a smooth torus, with some complex structure.
If $\aleph\!\neq\!\eset$, the principal components form a circle of spheres:
$$\Si_{\aleph}=\Big(\bigsqcup_{i\in I_0}\{i\}\!\times\!S^2\Big)\Big/\sim,
\qquad\hbox{where}\qquad
(i_1,\i)\sim(i_2,0)~~\hbox{if}~~(i_1,i_2)\!\in\!\aleph.$$
A {\tt genus-one $\Pf$-valued bubble map} is a tuple
$b\!=\!\big(I,\aleph;S,x,u\big)$,
where $S$ is a smooth Riemann surface of genus one if $\aleph\!=\!\eset$
and the circle of spheres $\Si_{\aleph}$ otherwise.
The objects $x$, $u$, and $(\Si_b,u_b)$ are as in 
the genus-zero case, except 
the sphere $\Si_{b,\hat{0}}$ is replaced by the genus-one curve $\Si_{b,\aleph}\!\equiv\!S$.
Furthermore, if $\aleph\!=\!\eset$, and thus $I_0\!=\!\{\hat{0}\}$ is a single-element set,
$u_{\hat{0}}\!\in\!C^{\i}(S;\Pf)$.
In the genus-one case, the general structure of bubble maps is encoded by
the tuples of the form $\T\!=\!(I,\aleph;\under{d})$.
Similarly to the genus-zero case, we denote by $\U_{\T}(\Pf;J)$
the subset of $\ov\M_1(\Pf,d;J)$ 
consisting of stable maps $[{\cal C},u]$ such that
$[{\cal C},u]\!=\!\big[\Si_b,u_b]$,
for some bubble map $b$ of type $\T$ as above.\\

\noindent
If $\T\!=\!(I,\aleph;\under{d})$ is a bubble type as above, let
\begin{gather*}
I_1=\big\{h\!\in\!\hat{I}\!:\io_h\!\in\!I_0\big\},
\qquad
\T_0=\big(I_1,I_0,\aleph;\io|_{I_1},\under{d}|_{I_0}\big),\\
\hbox{and}\qquad
\hbox{Aut}^*(\T)=\hbox{Aut}(\T)/
\{g\!\in\!\hbox{Aut}(\T)\!:g\cdot h\!=\!h~\forall h\!\in\!I_1\}.
\end{gather*}
where $I_0$ is the subset of minimal elements of $I$.
For each $h\!\in\!I_1$, we put
$$I_h=\big\{i\!\in\!I\!:h\!\le\!i\big\} \qquad\hbox{and}\qquad
\T_h=\big(I_h;\under{d}|_{I_h}\big).$$
The tuple $\T_0$ describes bubble maps from genus-one Riemann surfaces
with the marked points indexed by the set~$I_1$;
see Subsection~\ref{g1comp-notation1_subs} in~\cite{g1comp}.
We have a natural isomorphism
\begin{equation}\label{g1gendecomp_e3}\begin{split}
\U_{\T}(\Pf;J)\approx \Big(\big\{
\big(b_0,(b_h)_{h\in I_1}\big)\!\in\!\U_{\T_0}(\Pf;J)\!\times\!\!
\prod_{h\in I_1}\!\!\U_{\T_h}(\Pf;J)\!:\qquad&\\
\ev_0(b_h)\!=\!\ev_{\io_h}(b_0)~\forall h\!\in\!I_1&\big\}\Big)
\big/\hbox{Aut}^*(\T).
\end{split}\end{equation}
This decomposition is illustrated in Figure~\ref{g1gendecomp_fig}.
In this figure, we represent an entire stratum of bubble maps
by the domain of the stable maps in that stratum.
The right-hand side of Figure~\ref{g1gendecomp_fig} 
represents the subset of the cartesian product of the three spaces
of bubble maps, corresponding to the three drawings,
on which the appropriate evaluation maps agree pairwise,
as indicated by the dotted lines and defined in~\e_ref{g1gendecomp_e3}.\\

\begin{figure}
\begin{pspicture}(-1.1,-2)(10,1.25)
\psset{unit=.4cm}
\psellipse(8,-1.5)(1.5,2.5)
\psarc[linewidth=.05](6.2,-1.5){2}{-30}{30}\psarc[linewidth=.05](9.8,-1.5){2}{150}{210}
\pscircle[fillstyle=solid,fillcolor=gray](5.5,-1.5){1}\pscircle*(6.5,-1.5){.2}
\pscircle[fillstyle=solid,fillcolor=gray](3.5,-1.5){1}\pscircle*(4.5,-1.5){.2}
\pscircle(10.5,-1.5){1}\pscircle*(9.5,-1.5){.2}
\pscircle[fillstyle=solid,fillcolor=gray](11.91,-.09){1}\pscircle*(11.21,-.79){.2}
\pscircle[fillstyle=solid,fillcolor=gray](11.91,-2.91){1}\pscircle*(11.21,-2.21){.2}
\rput(5.5,0){$h_1$}\rput(3.5,0){$h_2$}\rput(10.3,0){$h_3$}
\rput(13.5,0.1){$h_4$}\rput(13.5,-2.9){$h_5$}
\rput(17.5,-1.5){$\approx$}
\psellipse(23,-1.5)(2.5,1.5)
\psarc[linewidth=.05](23,.3){2}{240}{300}\psarc[linewidth=.05](23,-3.3){2}{60}{120}
\psline[linewidth=.06,linestyle=dotted](28,1)(23,0)
\psline[linewidth=.06,linestyle=dotted](28,-3)(23,-3)
\pscircle*(23,0){.2}\pscircle*(23,-3){.2}
\pscircle[fillstyle=solid,fillcolor=gray](29,1){1}\pscircle*(28,1){.2}
\pscircle[fillstyle=solid,fillcolor=gray](31,1){1}\pscircle*(30,1){.2}
\pscircle(29,-3){1}\pscircle*(28,-3){.2}
\pscircle[fillstyle=solid,fillcolor=gray](30.41,-1.59){1}\pscircle*(29.71,-2.29){.2}
\pscircle[fillstyle=solid,fillcolor=gray](30.41,-4.41){1}\pscircle*(29.71,-3.71){.2}
\rput(29,2.5){$h_1$}\rput(31,2.5){$h_2$}\rput(28.8,-1.5){$h_3$}
\rput(32,-1.4){$h_4$}\rput(32,-4.1){$h_5$}
\end{pspicture}
\caption{An Example of the Decomposition~\e_ref{g1gendecomp_e3}}
\label{g1gendecomp_fig}
\end{figure}

\noindent
Let ${\cal FT}\!\lra\!\U_{\T}(\Pf;J)$
be the bundle of gluing parameters, or of smoothings at the nodes.
This orbi-bundle has the form 
$${\cal FT}=\Big(\!\bigoplus_{(h,i)\in\aleph}\!\!\!L_{h,0}\!\otimes\!L_{i,1}\oplus
\bigoplus_{h\in\hat{I}}L_{h,0}\!\otimes\!L_{h,1}\Big)\big/\hbox{Aut}(\T),$$
for certain line orbi-bundles $L_{h,0}$ and $L_{h,1}$.
Similarly to the genus-zero case,
\begin{gather}
\label{g1gendecomp_e4a}
\U_{\T}(\Pf;J)=\U_{\T}^{(0)}(\Pf;J)\big/
\hbox{Aut}(\T)\!\propto\!(S^1)^{\hat{I}},\qquad\hbox{where}\\
\label{g1gendecomp_e4b}
\U_{\T}^{(0)}(\Pf;J)=
\big\{\big(b_0,(b_h)_{h\in I_1}\big)\!\in\!\U_{\T_0}(\Pf;J)\!\times\!\!
\prod_{h\in I_1}\!\!\U_{\T_h}^{(0)}(\Pf;J)\!:
\ev_0(b_h)\!=\!\ev_{\io_h}(b_0)~\forall h\!\in\!I_1\big\}.
\end{gather}
The line bundles $L_{h,0}$ and $L_{h,1}$ arise from the quotient~\e_ref{g1gendecomp_e4a},
and 
$${\cal FT}=\ti{\cal F}\T\big/\hbox{Aut}(\T)\!\propto\!(S^1)^{\hat{I}},
\qquad\hbox{where}\qquad
\ti{\cal F}\T=\ti{\cal F}_{\aleph}\T\oplus
\bigoplus_{h\in\hat{I}}\ti{\cal F}_h\T,$$
$\ti{\cal F}_{\aleph}\T\!\lra\!\U_{\T}^{(0)}(\Pf;J)$
is the bundle of smoothings for the $|{\cal N}|$ nodes of the circle of spheres $\Si_{\aleph}$
and $\ti{\cal F}_h\T\!\lra\!\U_{\T}^{(0)}(\Pf;J)$ is the line bundle
of smoothings of the attaching node of the bubble indexed by~$h$.\\

\noindent
Suppose $\T\!=\!(I,\aleph;\under{d})$ is a bubble type such that
$d_i\!=\!0$ for all $i\!\in\!I_0$, 
i.e.~every element in $\U_{\T}(\Pf;J)$ is constant on the principal components.
In this case, the decomposition \e_ref{g1gendecomp_e3} is equivalent to
\begin{equation}\label{g1decomp_e1}\begin{split}
\U_{\T}(\Pf;J)&\approx \Big(\U_{\T_0}(pt)\times
\U_{\bar\T}(\Pf;J)\Big)\big/\hbox{Aut}^*(\T)\\
&\subset\Big(\ov\cM_{1,k}\times
\U_{\bar\T}(\Pf;J)\Big)\big/\hbox{Aut}^*(\T),
\end{split}\end{equation}
where $k\!=\!|I_1|$ and
$$\U_{\bar\T}(\Pf;J)=
\big\{(b_h)_{h\in I_1}\!\in\!\prod_{h\in I_1}\!\!\U_{\T_h}(\Pf;J)\!:
\ev_0(b_{h_1})\!=\!\ev_0(b_{h_2})~\forall h_1,h_2\!\in\!I_1\big\}.$$
Similarly, \e_ref{g1gendecomp_e4a} is equivalent~to
\begin{gather}\label{g1decomp_e3a}
\U_{\T}^{(0)}(\Pf;J)\approx \U_{\T_0}(pt)\times\U_{\bar\T}^{(0)}(\Pf;J)
\subset \ov\cM_{1,k}\times \U_{\bar\T}^{(0)}(\Pf;J),
\qquad\hbox{where}\\
\label{g1decomp_e3b}
\U_{\bar\T}^{(0)}(\Pf;J)=
\big\{(b_h)_{h\in I_1}\!\in\!\prod_{h\in I_1}\!\!\U_{\T_h}^{(0)}(\Pf;J)\!:
\ev_0(b_{h_1})\!=\!\ev_0(b_{h_2})~\forall h_1,h_2\!\in\!I_1\big\}.
\end{gather}
We denote by
$$\pi_P\!:\U_{\T}(\Pf;J)\lra\ov\cM_{1,k}/\hbox{Aut}^*(\T) 
\quad\hbox{and}\quad
\pi_P\!:\U_{\T}^{(0)}(\Pf;J)\lra\ov\cM_{1,k}$$
the projections onto the first component in the decompositions~\e_ref{g1decomp_e1}
and~\e_ref{g1decomp_e3a}.
Let
$$\ev_P\!: \U_{\T}(\Pf;J),\U_{\T}^{(0)}(\Pf;J)\lra \Pf$$
be the map sending each stable map $(\Si,u)$ to its value on the principal 
component $\Si_P$ of $\Si$, i.e~the point~$u(\Si_P)$.\\

\noindent
If $\T\!=\!(I,\aleph;\under{d})$ is as in the previous paragraph, let
\begin{gather*}
\chi(\T)=\big\{i\!\in\!\hat{I}\!:d_i\!\neq\!0;~
d_h\!=\!0~\forall h\!<\!i\big\};\\
\ti\F{\T}=\!\bigoplus_{i\in\chi(\T)}\!\!\ti{\cal F}_{h(i)}\T
\lra\U_{\T}^{(0)}(\Pf;J),
\quad\hbox{where}\quad
h(i)\!=\!\min\{h\!\in\!\hat{I}\!:h\!\le\!i\}\in I_1.
\end{gather*}
The subset $\chi(\T)$ of $I$ indexes the first-level effective
bubbles of every element of~$\U_{\T}^{(0)}(\Pf;J)$.
For each element $b\!=\!(\Si_b,u_b)$ of $\U_{\T}^{(0)}(\Pf;J)$
and $i\!\in\!\chi(\T)$, let
$${\cal D}_ib=\big\{du_b|_{\Si_{b,i}}\big\}\big|_{\i}e_{\i}\in T_{\ev_P(b)}\Pf,
\qquad\hbox{where}\quad e_{\i}=(1,0,0)\in T_{\i}S^2.$$
The complex span of ${\cal D}_ib$ in $T_{\ev_P(b)}\Pf$
is the tangent line to the rational component~$\Si_{b,i}$
at the node of~$\Si_{b,i}$ closest to a principal component of~$\Si_b$.
If the branch corresponding to~$\Si_{b,i}$ has a cusp at this node,
then ${\cal D}_ib\!=\!0$.\\

\noindent
Let $\E\!\lra\!\ov\M_{1,k}$ denote the Hodge line bundle,
i.e.~the line bundle of holomorphic differentials.
For each $i\!\in\!\chi(\T)$, we define the bundle map 
$${\cal D}_{J,i}\!:\ti{\cal F}_{h(i)}\T\lra 
\pi_P^*\E^*\otimes_J \ev_P^*T\Pf,$$
over $\U_{\T}^{(0)}(\Pf;J)$ by
$$\big\{{\cal D}_{J,i}(\ti\ups)\big\}(\psi)=
\psi_{x_{h(i)}(b)}(\ti{v})\cdot_J{\cal D}_ib \in T_{\ev_P(b)}\Pf
\quad\hbox{if}\quad \psi\!\in\!\pi_P^*\E,
~\ti\ups\!=\!(b,\ti{v})\!\in\!\ti{\cal F}_{h(i)}{\T},
~b\!\in\!\U_{\T}^{(0)}(\Pf;J),$$
and $x_{h(i)}(b)\!\in\!\Si_{b,\aleph}$ is the node joining the bubble $\Si_{b,h(i)}$
of $b$ to the principal component $\Si_{b,\aleph}$ of~$\Si_b$.
For each $\ups\!\in\!\ti{\cal F}\T$, we~put
\begin{gather*}
\rho(v)\!=\!\big(b;\rho_i(v)\big)_{i\in\chi(\T)} \in \ti\F\T,
\quad\hbox{where}\quad
\rho_i(v)\!=\!\!\prod_{h\in\hat{I},h\le i}\!\!\!\!\!\!v_h\in\ti{\cal F}_{h(i)}\T,
\qquad\hbox{if}\\
\ups\!=\!\big(b;v_{\aleph},(v_i)_{i\in\hat{I}}\big),~~
b\!\in\!\U_{\T}^{(0)}(\Pf;J),~~
(b,v_{\aleph})\!\in\!\ti{\cal F}_{\aleph}\T,~~
(b,v_h)\!\in\!\ti{\cal F}_h\T~\hbox{if}~h\!\in\!I_1,~~
v_i\!\in\!\C~\hbox{if}~i\!\in\!\hat{I}\!-\!I_1.
\end{gather*}
These definitions are illustrated in Figure~\ref{g1bdstr_fig} on page~\pageref{g1bdstr_fig}.
While the bundle maps $D_{J,i}$ and $\rho$ do not necessarily descend to the vector bundle 
${\cal FT}$ over $\U_{\T}(\Pf;J)$, the~map
\begin{gather*}
{\cal D}_{\T}\!:{\cal FT}\lra \pi_P^*\E^*\!\otimes\!\ev_P^*T\Pf
\big/\hbox{Aut}^*(\T), \qquad
{\cal D}_{\T}(\ups)=\sum_{i\in\chi(\T)}\!\!\!{\cal D}_{J,i}\rho_i(\ups),
\end{gather*}
is well-defined.\\

\noindent
Let $\ti\V_1^d\!\lra\!\U_{\T}^{(0)}(\Pf;J)$ be 
the vector bundle such~that the fiber of $\ti\V_1^d$
over a point $b\!=\!(\Si_b,u_b)$ in $\U_{\T}^{(0)}(\Pf;J)$
is $\ker\bpar_{\na,b}$, where $\na$ is the standard connection 
in line bundle $\L\!=\!\ga^{*\otimes5}$ over~$\Pf$;
see Subsection~\ref{review_subs}, as well as
Subsection~\ref{g1cone-g1conelocalstr_subs2} in~\cite{g1cone}.
If $b\!=\!(\Si_b,u_b)\!\in\!\U_{\T}^{(0)}(\Pf;J)$,
$\xi\!=\!(\xi_h)_{h\in I}\!\in\!\Ga(b;\L)$, and $i\!\in\!\chi({\T})$, let
$$\ti\D_{\T,i}\xi=\na_{e_{\i}}\xi_i\in\L_{\ev_0(b)},$$
as in Subsection~\ref{g1cone-notation1_subs} in~\cite{g1cone}.
We next define the bundle~map
$$\D_{\T}\!: \ti\V_1^d\!\otimes\!\ti\F\T\lra \pi_P^*\E^*\!\otimes\!\ev_P^*\L$$
over $\U_{\T}^{(0)}(\Pf;J)$ by
\begin{gather*}
\big\{\ti\D_{\T}(\xi\!\otimes\!\ti\ups)\big\}(\psi)
=\sum_{i\in\chi(\T)}\!\psi_{x_{h(i)}(b)}(\ti\ups_i)\cdot 
\ti\D_{\T,i}\xi \in\L_{\ev_P(b)} \qquad\hbox{if}\\
\xi\in\ti\V_1^d|_b\subset\Ga(b;\L),\qquad
\ti\ups=(\ti\ups_i)_{i\in\chi(\T)}\in\ti\F\T|_b,
\quad\hbox{and}\quad  \psi\in\E_{\pi_P(b)}.
\end{gather*}
The bundle map $\ti\D_{\T}$ induces a linear bundle~map over~$\U_{\T}(\P;J)$:
$$\D_{\T}\!:\V_1^d\!\otimes\!\F\T
\lra \pi_P^*\E^*\!\otimes\!\ev_P^*\L\big/\hbox{Aut}^*(\T),
 \quad\hbox{where}\quad 
\F\T=\!\Big(\bigoplus_{i\in\chi(\T)}\!\!
\pi_P^*L_{P,h(i)}\!\otimes\!\pi_i^*L_0\Big)\big/\hbox{Aut}^*(\T),$$
$L_{P,h}\!\lra\!\ov\M_{1,k}$ is the universal tangent line bundle at
the marked point $x_h$,
$L_0\!\lra\!\U_{\T'}(\Pf;J)$ is 
the universal tangent line bundle at the special marked point $(i,\i)$
for any bubble type $\T'$ of rational stable maps, and
$$\pi_i\!:\U_{\T}(\Pf;J)\lra\U_{\T_i'}(\Pf;J)$$
is the projection map sending each bubble map $b\!=\!(\Si_b,u_b)$ 
to its restriction to the component~$\Si_{b,i}$.\\

\noindent
Finally, if $\T$ is any bubble type, for genus-zero or genus-one maps,
and $K$ is a subset of $\U_{\T}(\Pf;J)$,
we denote by $K^{(0)}$ the preimage of $K$ under the quotient projection map
$\U_{\T}^{(0)}(\Pf;J)\!\lra\!\U_{\T}(\Pf;J)$.
All vector orbi-bundles we encounter will be assumed to be normed.
Some will come with natural norms; for others, we implicitly choose a norm
once and for~all.
If \hbox{$\pi_{\F}\!:{\frak F}\!\lra\!\X$} is a normed vector bundle
and $\de\!:\X\!\lra\!\R$ is any function, possibly constant,
let
$$\F_{\de}=\big\{\ups\!\in\!\F\!: |\ups|\!<\!\de(\pi_{\F}(\ups))\big\}.$$
If $\Om$ is any subset of $\F$, we take  $\Om_{\de}\!=\!\Om\cap\F_{\de}$.

\section{On Genus-One Gromov-Witten Invariants}
\label{dbar_sec}

\subsection{Setup}
\label{setup_subs}

\noindent
In this section, we prove Propositions~\ref{maincontr_prp} and~\ref{bdcontr_prp}.
We start by clarifying the setup described after Proposition~\ref{bdcontr_prp}.
We also specify the open subsets of admissible perturbations
of the $\bpar_J$-operator
to be used in proving Propositions~\ref{maincontr_prp} and \ref{bdcontr_prp};
see Definition~\ref{contr_dfn}.\\

\noindent
Let $U_s$ be the neighborhood of~$Y$ in~$\Pf$ and 
$\ti{T}Y$ the subbundle of $T\Pf|_{U_s}$ as in Subsection~\ref{outline_subs2}.
We~set
$$\X_s=\big\{[\Si,j;u]\!\in\!\X_1(\Pf,d)\!:u(\Si)\!\subset\!U_s\big\}.$$
Let $\nu$ be a multisection of the bundle $\Ga^{0,1}_1(\Pf,d)$ such~that\\
${}\quad$ ($\nu1$) for every open neighborhood $\U$ of $\ov\M_1(\Pf,d;J)$ in $\X_1(\Pf,d)$,
there exists $\ep_{\nu}(\U)\!>\!0$\\  
${}\qquad~~~$ such that $\{\bpar_J\!+\!t\nu\}^{-1}(0)$ is contained in~$\U$
               for all $t\!\in\!(0,\ep_{\nu}(\U))$;\\
${}\quad$ ($\nu2$) $\nu(b)\!\in\!\Ga(\Si;\La^{0,1}_{J,j}T^*\Si\!\otimes\!u^*\ti{T}Y)/
\hbox{Aut}(b)$ if $b\!=\![\Si,j,u]\!\in\!\X_s$, and
$\nu(b)\!=\!0$ if $b\!\not\in\!\X_s$;\\
${}\quad$ ($\nu3$) for some $\ep_{\nu}\!>\!0$ and for all $t\!\in\!(0,\ep_{\nu})$,
       the multisection $\bar{\partial}_J\!+\!t\nu$ does not vanish on\\
${}\qquad~~~$ $\X_1(Y,d)\!-\!\X_1^0(Y,d)$ and is transversal to 
               the zero set in $\Ga^{0,1}_1(Y,d;J)$ along~$\X_1^0(Y,d)$.\\
The middle condition implies that $\bpar_{\na,u}\{s\!\circ\!u\}\!=\!0$ if 
$[\Si,j,u]\!\in\!\{\bpar_J\!+\!t\nu\}^{-1}(0)$.
It can be shown, by slightly modifying the proof of Corollary~\ref{bdcontr_crl2},
that the finite-dimensional conditions ($\nu3a$)-($\nu3c$) stated below
imply~($\nu3$).\\

\noindent
If $\nu$ is a section of the bundle $\Ga^{0,1}_1(\Pf;d)$ over 
$\X_1(\Pf,d)$ as in~($\nu1$) and~($\nu2$) above,
for all $\ka\!\in\!\S_0(Y;J)$, we define a section of the bundle
$$\W_{\ka,d/d_{\ka}}^{1,1}\lra\ov\M_1^1(\ka,d/d_{\ka})
\qquad\hbox{by}\qquad \pi_{\nu,\ka}^1(b)=[\nu(b)],$$
where $[\nu(b)]$ is the $(0,1)$-cohomology class of $\nu(b)$
and $\W_{\ka,d/d_{\ka}}^{1,1}$ is as in Subsection~\ref{outline_subs2}.
For each $q\!\in\!\Z^+$, we define a section of the bundle
\begin{gather}\label{redcoker_e}
\wt\W_{\ka,d/d_{\ka}}^{1,q}\!\equiv\!
\big(\pi_P^*\E^*\!\otimes\!\pi_B^*\ev_0^*N_Y\ka\oplus
\pi_B^*\W_{\ka,d/d_{\ka}}^{0,q}\big)/S_q,
\lra \ov\M_1^q(\ka,d/d_{\ka})\\
\hbox{by}\qquad
\ti\pi_{\nu,\ka}^q(b)=\pi_{\ka}^{\perp}[\nu(b)],\notag
\end{gather}
where 
$$\pi_{\ka}^{\perp}\!:\W_{\ka,d/d_{\ka}}^{1,q}\lra\wt\W_{\ka,d/d_{\ka}}^{1,q}$$
is the projection map corresponding to the quotient of $\W_{\ka,d/d_{\ka}}^{1,q}$
by $\pi_P^*\E^*\!\otimes\!\pi_B^*\ev_0^*T\ka$; see~\e_ref{cokerkdecomp_e}.
Finally, we define a section of the bundle
$$\W_{\ka,d/d_{\ka}}^{1,0}\lra\ov\M_1^0(\ka,d/d_{\ka})
\qquad\hbox{by}\qquad
\pi_{\nu,\ka}^0(b)=
\begin{cases}
\ti\pi_{\nu,\ka}^q(b),&
\hbox{if}~b\!\in\!\ov\M_1^q(\ka,d/d_{\ka}),~q\!\in\!\Z^+;\\
[\nu(b)],&\hbox{otherwise};
\end{cases}$$
see~\e_ref{coker0decomp_e2}.  This section is well-defined on 
$\ov\M_1^{q_1}(\ka,d/d_{\ka})\cap\ov\M_1^{q_2}(\ka,d/d_{\ka})$.\\

\noindent
We denote by $\ti{\cal A}_1^d(\bar{\partial},J)$ the space of multisections~$\nu$
as in~($\nu1$) and~($\nu2$) such that for all $\ka\!\in\!\S_0(Y;J)$:\\
${}\quad$ ($\nu3a$) the section $\pi_{\nu,\ka}^0$ does not vanish on
$\ov\M_1^0(\ka,d/d_{\ka})\!-\!\M_1^0(\ka,d/d_{\ka})$ and is transversal to\\
${}\qquad\quad~~$ the zero set on~$\M_1^0(\ka,d/d_{\ka})$;\\
${}\quad$ ($\nu3b$) the section $\pi_{\nu,\ka}^1$ does not vanish on
$\ov\M_1^1(\ka,d/d_{\ka})$;\\
${}\quad$ ($\nu3c$) the section $\ti{\pi}_{\nu,\ka}^1$ does not vanish on
$\ov\M_1^1(\ka,d/d_{\ka})\!-\!\M_1^1(\ka,d/d_{\ka})$ and is transversal to\\ 
${}\qquad\quad~~$ the zero set on~$\M_1^1(\ka,d/d_{\ka})$.\\
By \e_ref{spacesdim_e1}, \e_ref{spacesdim_e0}, \e_ref{rankk_e1}, 
and Lemmas~\ref{g1str_lmm1} and~\ref{g1str_lmm2}, 
these conditions are satisfied by a dense open path-connected subset of sections~$\nu$.

\subsection{Proof of Proposition~\ref{maincontr_prp}}
\label{maincontr_subs}

\noindent
We will focus on the last case of Proposition~\ref{maincontr_prp},
which follows from  Proposition~\ref{maincontr_prp2}.
The claim in the first case is clear,
since the single-element set $\M_1^0(\ka,d/d_{\ka})$ consists of 
a transverse zero of the section $\bpar_J$ over~$\X_1(Y,d)$.
The proof of Proposition~\ref{maincontr_prp2} applies to this case as well,
except there is no gluing to be done.\\

\noindent
Let $s$ and $Y$ be as in Proposition~\ref{maincontr_prp}.
For every bubble type $\T$ and every rational $J$-holomorphic curve
$\ka$ in~$Y$, we~put
$$\U_{\T;\ka}=\big\{[\cC,u]\!\in\!\U_{\T}(\Pf;J)\!: u(\C)\!=\!\ka\big\}.$$

\begin{prp}
\label{maincontr_prp2}
Suppose $d$, $Y$, and $J$ are as in Proposition~\ref{maincontr_prp},
$\nu\!\in\!\ti{\cal A}_1^d(\bar{\partial};J)$ is a generic perturbation
of the $\bar{\partial}_J$-operator on $\X_1(\Pf,d)$,
$\ka\!\in\!\S_0(Y;J)$, and
$\T\!=\!(I,\aleph;\under{d})$ is a bubble type such that 
$\sum_{i\in I}\!d_i\!=\!d$ and
$d_i\!\neq\!0$ for some minimal element $i$ of~$I$.
If $|I|\!>\!1$ or $\aleph\!\neq\!\eset$, 
for every compact subset~$K$ of $\U_{\T;\ka}$,
there exist $\ep_{\nu}(K)\!\in\!\R^+$ and 
an open neighborhood~$U(K)$ of~$K$ in $\X_1(Y,d)$ such~that
$$\{\bpar_J\!+\!t\nu\}^{-1}(0)\!\cap\! U(K)=\eset
\qquad\forall t\!\in\!(0,\ep_{\nu}(K)).$$
If $|I|\!=\!1$ and $\aleph\!=\!\eset$, for every compact subset~$K$ of $\U_{\T;\ka}$, 
there exist $\ep_{\nu}(K)\!\in\!\R^+$ and an open neighborhood~$U(K)$ 
of~$K$ in $\X_1(Y,d)$ with the following properties:\\
${}\quad$ (a) the section $\bar{\partial}_J\!+\!t\nu$ is transverse to the zero set
in $\Ga_1^{0,1}(Y,d;J)$ over $U(K)$ for all $t\!\in\!(0,\ep_{\nu}(K))$;\\
${}\quad$ (b) for every open subset $U$ of $\X_1(Y,d)$,
there exists $\ep(U)\!\in\!(0,\ep_{\nu}(K))$ such that
$$^{\pm}\!\big|\{\bpar_J\!+\!t\nu\}^{-1}\!\cap\!U\big|=
\blr{e(\W_{\ka,d/d_{\ka}}^{1,0}),\big[\ov\M_1^0(\ka,d/d_{\ka})\big]}
\quad\hbox{if}\quad
\pi_{\nu,\ka}^{0~-1}(0)\!\subset\!K\!\subset\!U\!\subset\!U(K),
~t\!\in\!(0,\ep(U)).$$
\end{prp}

\noindent
In other words, the contribution from  the main stratum $\M_1^0(\ka,d/d_{\ka})$  of 
$\ov\M_1^0(\ka,d/d_{\ka})$ to the number $N_1(d)$,
as computed via the section $\bpar_J$, is the euler class of the vector bundle 
$\W_{\ka,d/d_{\ka}}^{1,0}$ over~$\ov\M_1^0(\ka,d/d_{\ka})$.
None of the boundary strata of $\M_1^{\{0\}}(\ka,d/d_{\ka})$
contributes to~$N_1(d)$.\\ 

\noindent
We fix a $J$-compatible metric $g_{\Pf}$ on $\Pf$ and proceed as in 
Subsection~\ref{g1comp-reg1_subs1} of~\cite{g1comp}.
For each sufficiently small element $\ups\!=\!(b,v)$ of $\ti{\cal F}{\T}^{\eset}$,
let 
$$b(\ups)=\big(\Si_{\ups},j_{\ups};u_{\ups}\big),
\qquad\hbox{where}\qquad
u_{\ups}=u_b\circ q_{\ups},$$
be the corresponding approximately holomorphic stable map.
Here
$$q_{\ups}\!:\Si_{\ups}\lra\Si_b$$
is the basic gluing map constructed in Subsection~\ref{g1comp-reg1_subs1} of~\cite{g1comp}.
Since $d_i\!\neq\!0$ for some minimal element $i$ of~$I$, 
i.e.~the stable map~$b$ is non-constant on the principal curve of the domain~$\Si_b$ of~$b$,
the linearization~$D_{J,b}$ of the $\bar{\partial}_J$-operator at~$b$ is surjective,
since $\|J\!-\!J_0\|_{C^1}\!\le\!\de(b)$.
Thus, if $\ups$ is sufficiently small, the linearization
$$D_{J,\ups}\!:\Ga(\ups;T\Pf)\!\equiv\!L^p_1(\Si_{\ups};u_{\ups}^*T\Pf)\lra
\Ga^{0,1}(\ups,T\Pf;J)\!\equiv\!
L^p(\Si_{\ups};\La^{0,1}_{J,j}T^*\Si_{\ups}\!\otimes\!u_{\ups}^*T\Pf),$$
of the $\bar{\partial}_J$-operator at~$b(\ups)$,
defined via the $J$-compatible connection $\na^J$ in $T\Pf$
corresponding to the Levi-Civita connection of the metric $g_{\Pf}$, 
 is also surjective.
In particular, we can obtain an orthogonal decomposition
\begin{equation}\label{gadecomp_e1}
\Ga(\ups;T\Pf)=\Ga_-(\ups;T\Pf)\oplus\Ga_+(\ups;T\Pf)
\end{equation}
such that the linear operator 
$$D_{J,\ups}\!:\Ga_+(\ups;T\Pf)\lra\Ga^{0,1}(\ups;T\Pf;J)$$  
is an isomorphism, while 
$$\Ga_-(\ups;T\Pf)=\big\{\xi\circ q_{\ups}\!:  \xi\!\in\!\Ga_-(b;\Pf)\big\},
\qquad\hbox{where}\qquad
\Ga_-(b;T\Pf)=\ker D_{J,b}.$$
The $L^2$-inner product on~$\Ga(\ups;T\Pf)$ used in the orthogonal decomposition
is defined via the metric~$g_{\Pf}$ on~$\Pf$ and 
the metric~$g_{\ups}$ on~$\Si_{\ups}$ induced by the pregluing construction.
The Banach spaces $\Ga(\ups;T\Pf)$ and $\Ga^{0,1}(\ups;T\Pf;J)$ carry the norms 
$\|\cdot\|_{\ups,p,1}$ and $\|\cdot\|_{\ups,p}$, respectively,
which are also defined by the pregluing construction.
These norms are equivalent to the ones used in~\cite{LT}.
In particular, the norms of $D_{J,\ups}$ and of the inverse of its restriction
to $\Ga_+(\ups;T\Pf)$ have fiberwise uniform upper bounds, 
i.e.~dependent only on $[b]\!\in\!\U_{\T}(\Pf;J)$, and
not on $\ups\!\in\!\ti{\cal F}{\T}^{\eset}$.

\begin{lmm}
\label{maincontr_lmm1}
If ${\T}$ is a bubble type and  $\nu$ is an admissible perturbation of 
the $\bar{\partial}_J$-operator on $\X_1(\Pf,d)$  as in Proposition~\ref{maincontr_prp2},
for every precompact open subset $K$ of $\U_{\T}(\Pf;J)$,
there exist $\de_K,\ep_K,C_K\!\in\!\R^+$
and an open neighborhood $U_K$ of $K$ in $\X_1(\Pf,d)$
with the following properties:\\
(1) for all $\ups\!=\!(b,v)\!\in\!\ti{\cal F}{\T}^{\eset}|_{K^{(0)}}$,
\begin{gather*}
\|D_{J,\ups}\xi\|_{\ups,p}\le C_K|\ups|^{1/p}\|\xi\|_{\ups,p,1} 
\quad\forall\xi\!\in\!\Ga_-(\ups;T\Pf)
\qquad\hbox{and}\\
C_K^{-1}\|\xi\|_{\ups,p,1}\le \|D_{J,\ups}\xi\|_{\ups,p} \le C_K\|\xi\|_{\ups,p,1}
\quad\forall  \xi\!\in\!\Ga_+(\ups;T\Pf);
\end{gather*}
(2) for all $\ups\!=\!(b,v)\!\in\!\ti{\cal F}{\T}^{\eset}|_{K^{(0)}}$
and  $t\!\in\![0,\de_K)$, the equation
$$\bar{\partial}_J\exp_{u_{\ups}}\!\xi+t\nu(\exp_{u_{\ups}}\!\xi)=0,
\qquad \xi\!\in\!\Ga_+(\ups;T\Pf),~~\|\xi\|_{\ups,p,1}\!\le\!\ep_K,$$
has a unique solution $\xi_{t\nu}(\ups)$, and 
$\|\xi_{t\nu}(\ups)\|_{C^0}\!\le\!C_K(t\!+\!|\ups|^{1/p})$;\\
(3) there exist a smooth bundle map 
$\ze_{\nu}\!:\ti{\cal F}{\T}^{\eset}\!\lra\!\ti\Ga(T\Pf,d)$ 
over $\U_{\T}^{(0)}(\Pf;J)$ and
a continuous function $\ve_{\nu}\!:{\cal FT}^{\eset}\!\lra\!\R$
such that for all $\ups\!=\!(b,v)\!\in\!\ti{\cal F}{\T}^{\eset}|_{K^{(0)}}$
and $t\!\in\![0,\de_K)$, 
$$\big\|\xi_{t\nu}(\ups)-\xi_0(\ups)-tq_{\ups}^*\ze_{\nu}(b)\big\|_{C^0}
\le C_K\big(t\!+\!\ve_{\nu}(\ups)\big)t
\qquad\hbox{and}\qquad
\lim_{\ups\lra b}\!\!\ve_{\nu}(\ups)=0;$$
(4) the map 
$$\phi_{\T,t\nu}\!: {\cal FT}_{\de_K}^{\eset}\big|_K\!\lra\!\X_1^0(\Pf,d),
\quad [\ups]\!\lra\!\big[\ti{b}_{t\nu}(\ups)\big],
\qquad\hbox{where}\quad
\ti{b}_{t\nu}(\ups)=\big(\Si_{\ups},j_{\ups};\exp_{u_{\ups}}\!\xi_{t\nu}(\ups)\big),$$
is an orientation-preserving diffeomorphism
onto $\big\{\bar{\partial}_J\!+\!t\nu\}^{-1}(0)\cap\X_1^0(\Pf,d)\cap U_K$.
\end{lmm}

\noindent
The first claim of the lemma is a special case of Lemma~\ref{g1comp-reg1_lmm1} in~\cite{g1comp}.
The second statement is obtained by expanding the equation at~$u_{\ups}$
and applying the Contraction Principle; see Subsection~3.6 in~\cite{gluing}.
The uniqueness part means that there is a unique solution
for each branch of the multisection~$\nu$.
In~(3), $\ti\Ga(T\Pf,d)$ denotes the Banach bundle over 
the space $\U_{\T}^{(0)}(\Pf,d;J)$ such that
$$\ti\Ga(T\Pf,d)\big|_{(\Si_b,u_b)}=\Ga(\Si_b;u_b^*T\Pf).$$
Let $P_{\ups}$ and $P_b$ denote the inverses of~$D_{J,\ups}$ on $\Ga_+(\ups;T\Pf)$
and of~$D_{J,b}$ on~$\Ga_+(b;T\Pf)$, respectively.
The Banach space $\Ga_+(b;T\Pf)$ is the orthogonal complement of $\Ga_-(b;T\Pf)$
in 
$$\Ga(b;T\Pf)\equiv L^p_1(\Si_b;u_b^*T\Pf);$$
see Subsection~3.1 in~\cite{gluing}.
Taking the difference of the expansions for the equations
in~(2) describing $\xi_{t\nu}(\ups)$ and $\xi_0(\ups)$ and applying~$P_{\ups}$,
one finds that
$$\big\|\xi_{t\nu}(\ups)-\xi_0(\ups)-tP_{\ups}\nu(u_{\ups})\big\|_{C^0}
\le C_K(t\!+\!|\ups|^{1/p})t.$$
On the other hand, a direct computation shows that 
\begin{equation*}\begin{split}
\|P_{\ups}\nu(u_{\ups})-q_{\ups}^*P_b\nu(u_b)\big\|_{C^0}
&\le C(b)\|\nu(u_{\ups})-D_{J,\ups}q_{\ups}^*P_b\nu(u_b)\big\|_{\ups,p}+
\ti{\ve}_{\nu}(\ups)\\
&\le C(b)\|\nu(u_{\ups})-q_{\ups}^*\nu(u_b)\big\|_{\ups,p}
+\ti{\ve}_{\nu}'(\ups)
\le \ve_{\nu}(\ups);
\end{split}\end{equation*}
see Subsection~4.1 in~\cite{gluing} for a similar computation.
These two bounds imply (3) of Lemma~\ref{maincontr_lmm1}, 
with $\ze_{\nu}(b)\!=\!P_b\nu(u_b)$.
Finally, the proof of (4) is similar to 
Subsections~3.8 and  4.3-4.5 of~\cite{gluing}.

\begin{lmm}
\label{maincontr_lmm2}
Suppose ${\T}$ and $\nu$ are as in Lemma~\ref{maincontr_lmm1}.
For every precompact open subset $K$ of $\U_{\T}(\Pf;J)$,
there exist $\de_K,\ep_K,C_K\!\in\!\R^+$,
an open neighborhood $U_K$ of $K$ in $\X_1(\Pf,d)$, 
and injective vector-bundle homomorphisms 
$$\ti{\phi}_{\T,t\nu}\!:
\pi_{\cal FT}^*\V_1^d\big|_{{\cal F}^{\eset}{\T}_{\de_K}}\lra\Ga(\L;d),$$ 
covering  the maps $\phi_{\T,t\nu}$ of Lemma~\ref{maincontr_lmm1},
with the following properties:\\
(1) requirements (1)-(4) of Lemma~\ref{maincontr_lmm1} are satisfied;\\
(2) $\lim_{(\ups,w)\lra (b,w^*)}\ti{\phi}_{\T,t\nu}(\ups;w)\!=\!w^*$ 
for all $b\!\in\!K$ and $w^*\!\in\!\V_1^d$;\\
(3) $s_1^d(\phi_{\T,t\nu}(\ups))\!\equiv\![s\!\circ\!\exp_{u_{\ups}}\!\xi_{t\nu}]
\in\Im\ti{\phi}_{\T,t\nu}$, and 
for all $[\ups]\!=\![b,v]\!\in\!{\cal FT}_{\de_K}^{\eset}\big|_K$
$$\big|\ti\phi_{\T,t\nu}^{-1}s_1^d\big(\phi_{\T,t\nu}(\ups)\big)
-\ti\phi_{\T,0}^{-1}s_1^d\big(\phi_{\T,0}(\ups)\big)-
t\{\na s\}\ze_{\nu}(b)\big|
\le C_K\big(t\!+\!\ve_{\nu}(\ups)\big)t,$$
where $\ve_{\nu}\!:{\cal FT}^{\eset}\!\lra\!\R$ is
a continuous function 
such that $\lim\limits_{\ups\lra b}\!\ve_{\nu}(\ups)\!=\!0$
for all $b\!\in\!\U_{\T}(\Pf;J)$.
\end{lmm}

\noindent
{\it Proof:} (1)
We need to construct a lift $\ti\phi_{\T,t\nu}$ that has the desired properties.
For each element 
$\ups\!=\!(b,v)$ of $\ti{\cal F}\T^{\eset}_{\de_K}\big|_{K^{(0)}}$,
$t\!\in\![0,\de_K)$, and $\xi\!\in\!\Ga(b;\L)$, let
$$R_{\ups}\xi\in
\Ga(\ups;\L)\!\equiv\!L^p_1\big(\Si_{\ups};u_{\ups}^*\L\big)
\quad\hbox{and}\quad
R_{\ups,t\nu}\xi\in
\Ga(\ti{b}_{t\nu}(\ups);\L)
\!\equiv\!L^p_1\big(\Si_{\ups};\{\exp_{u_{\ups}}\!\xi_{t\nu}(\ups)\}^*\L\big)$$
be defined by
$$\big\{R_{\ups}\xi\big\}(z)=\xi\big(q_{\ups}(z)\big)
\quad\hbox{and}\quad
\big\{R_{\ups,t\nu}\xi\big\}(z)=
\Pi_{\{\xi_{t\nu}(\ups)\}(z)}\big\{R_{\ups}\xi\big\}(z)
\qquad\forall~z\!\in\!\Si_{\ups},$$
where $\Pi_{\{\xi_{t\nu}(\ups)\}(z)}\big\{R_{\ups}\xi\big\}(z)$
is the $\na$-parallel transport of $\big\{R_{\ups}\xi\big\}(z)$
along the $\na^J$-geodesic 
$$\ga_{\{\xi_{t\nu}(\ups)\}(z)}\!:[0,1]\lra\Pf,\qquad
\tau\lra\exp_{u_{\ups}(z)}\tau\{\xi_{t\nu}(\ups)\}(z).$$
We denote the image of 
$$\Ga_-(b;\L)\equiv\ker\bpar_{\na,b}$$
under the linear map $R_{\ups,t\nu}$ by $\ti\Ga_-(\ti{b}_{t\nu}(\ups);\L)$.
If $\de_K$ is sufficiently small, the $L^2$-orthogonal projection
$$\ti\pi_{\ups,t\nu}\!:\Ga(\ti{b}_{t\nu}(\ups);\L)
\lra\ti\Ga_-(\ti{b}_{t\nu}(\ups);\L),$$
defined with respect to the metric~$g_{\ups}$ on~$\Si_{\ups}$,
restricts to an isomorphism~on 
$$\Ga_-\big(\ti{b}_{t\nu}(\ups);\L\big)\equiv
\ker\bpar_{\na,\ti{b}_{t\nu}(\ups)};$$
see Subsection~\ref{g1cone-g1conelocalstr_subs1} in~\cite{g1cone}.
Let $\ti{\pi}_{\ups,t\nu}^{-1}$ be the inverse of this isomorphism.
We set
$$\ti{\phi}_{\T,t\nu}([\ups;\xi])=
[\ti{\pi}_{\ups,t\nu}^{-1}R_{\ups,t\nu}\xi]\qquad
\forall\xi\!\in\!\Ga_-(b;\L).$$
(2) By our assumptions on~$\nu$, 
$$\bpar_{\na,\ti{b}_{t\nu}(\ups)}\big(s\circ\exp_{u_{\ups}}\xi_{t\nu}(\ups)\big)=0
\qquad\Lra\qquad
s_1^d\big(\phi_{\T,t\nu}(\ups)\big)\in\Im\ti{\phi}_{\T,t\nu}.$$
It remains to prove the estimate in part~(3) of the lemma.
If $\ep_K$ is sufficiently small,
$\ups\!\in\!\ti{\cal F}\T_{\de_K}^{\eset}|_{K^{(0)}}$,
$\xi\!\in\!\Ga(\ups;T\Pf)$, and $\|\xi_{\ups}\|_{\ups,p,1}\!<\!\ep_K$,
we define 
$$N^s_{\ups}\xi\!\in\!\Ga\big(\ups;\L\big)
\quad\hbox{by}\quad
\Pi_{\xi(z)}^{-1}s\big(\!\exp_{u_{\ups}(z)}\!\xi(z)\big)
=s(u_{\ups}(z))+\na s\big|_{u_{\ups}(z)}\xi(z)+
\big\{N^s_{\ups}\xi\big\}(z)~~\forall z\!\in\!\Si_{\ups}.$$
The quadratic term $N^s_{\ups}$ varies smoothly with $\ups$, 
$N^s_{\ups}0\!=\!0$, and 
\begin{equation}\label{s2term_e}
\big\|N^s_{\ups}\xi_1-N^s_{\ups}\xi_2\big\|_{C^0}\le
C_s\big(\|\xi_1\|_{C^0}\!+\!\|\xi_2\|_{C^0}\big)
\|\xi_1\!-\!\xi_2\|_{C^0}
\end{equation}
for some $C_s\!\in\!\R^+$ and for all 
$\xi_1,\xi_2\!\in\!\Ga(\ups)$ such that
$\|\xi_1\|_{\ups,p,1},\|\xi_2\|_{\ups,p,1}\!<\!\ep_{\T,\nu}(K)$.
If $\xi\!\in\!\Ga_-(b;\L)$,
$$\bigllrr{s_1^d\big(\phi_{\T,t\nu}(\ups)\big),R_{\ups,t\nu}\xi}
=\bigllrr{\Pi_{\xi_{t\nu}(\ups)}^{-1}s_1^d(\phi_{\T,t\nu}(\ups))
,\xi\circ q_{\ups}}.$$
Thus, the estimate in (3) of Lemma~\ref{maincontr_lmm2}
follows from \e_ref{s2term_e} and the estimate in~(3) of Lemma~\ref{maincontr_lmm1}.

\begin{crl}
\label{maincontr_lmm2crl}
Suppose ${\T}$ and $\nu$ are as in Lemma~\ref{maincontr_lmm1}.
For every precompact open subset $K$ of $\U_{\T}(\Pf;J)$,
there exist $\de_K,\ep_K,C_K\!\in\!\R^+$,
an open neighborhood $U_K$ of $K$ in $\X_1(\Pf,d)$,
and for each $t\!\in\!(0,\ep(K))$ a sign-preserving bijection
$$\{\bar{\partial}_J\!+\!t\nu\}^{-1}(0)\cap\X_1(Y,d)\cap U_K\lra
\big\{u\!\in\!\M_1^0(\Pf,d;J)\cap U_K\!:
\{s_1^d\!+\!t\vt_t\}(u)\!=\!0\big\},$$
where $\vt_t\!\in\!\Ga(\M_1^0(\Pf,d;J)\cap U_K;\V_1^d)$
is a family of smooth sections such that
\begin{gather*}
\lim_{\ups\lra b,t\lra0}\vt_t(\ups)=[\{\na s\}P_b\nu(b)]
\qquad\forall b\!\in\!K\qquad\hbox{and}\\
\big|\na_X\ti{\phi}_{\T,0}^{-1}s_1^d\phi_{\T,0}(\ups)\big|
\le C_K|X|
\qquad\forall b\!\in\!K,~\ups\!\in\!{\cal FT}_{\de_K}^{\eset}\big|_b,~
X\!\in\!\ker D_{J,b},
\end{gather*}
where $\phi_{\T,0}$ and $\ti{\phi}_{\T,0}$
are as in Lemma~\ref{maincontr_lmm2}.
\end{crl}

\noindent 
{\it Proof:} The section $\vt_t$ is given~by
$$t\vt_t\big(\phi_{\T,0}(\ups)\big)
=\ti\phi_{\T,0}\ti\phi_{\T,t\nu}^{-1}
                        \big(s_1^d(\phi_{\T,t\nu}(\ups))\big)
-s_1^d(\phi_{\T,0}(\ups))
\qquad\forall\, [\ups]\!\in\!{\cal FT}^{\eset}_{\de_K}|_{K^{(0)}}.$$
This corollary is immediate from Lemma~\ref{maincontr_lmm2},
with the exception of the last estimate.
This estimate follows from the behavior of the various terms involved in defining~$\vt_t$;
see Subsections~3.4 and~4.2 in~\cite{gluing}.\\

\noindent
For each $\ka\!\in\!\S_0(J;Y)$, $\U_{\T;\ka}$ 
is  a smooth suborbifold of $\U_{\T}(\Pf;J)$.
We denote its normal bundle by~${\cal N}^{\ka}{\T}$.
Its fiber at $[b]\!\in\!\U_{\T;\ka}$ 
is the quotient $\Ga_-(b;T\Pf)/\Ga_-(b;TY)$, where
$$\Ga_-(b;TY)\equiv\Ga_-(b;T\Pf)\cap\Ga(b;TY) =\Ga_-(b;T\ka),$$
by the assumption ($J_Y2$) on $J_Y$.
We identify ${\cal N}^{\ka}{\T}$ with the $L^2$-orthogonal complement
of $\Ga_-(b;TY)$ in $\Ga_-(b;T\Pf)$.
Let
$$\vph_{\T;\ka}\!:
{\cal N}^{\ka}{\T}_{\de_{\ka}}\lra\U_{\T}(\Pf;J)$$
be an orientation-preserving identification of neighborhoods
of $\U_{\T;\ka}$ in ${\cal N}^{\ka}{\T}$
and in $\U_{\T}(\Pf;J)$ and~let
$$\ti{\vph}_{\T;\ka}\!:\pi_{{\cal N}^{\ka}{\T}}^*{\cal FT}
\big|_{{\cal N}^{\ka}{\T}_{\de_{\ka}}}\lra{\cal FT}
\quad\hbox{and}\quad
\ti{\vph}_{\T;\ka}\!:\pi_{{\cal N}^{\ka}{\T}}^*\V_1^d
\big|_{{\cal N}^{\ka}{\T}_{\de_{\ka}}}\lra\V_1^d$$
be lifts of~$\vph_{\T;\ka}$ to vector-bundle isomorphisms
restricting to the identity over~$\U_{\T;\ka}$.\\

\noindent
The section $s_1^d$ is smooth on $\U_{\T}(\Pf;J)$ and its differential along $\U_{\T;\ka}$, 
i.e.~the homomorphism~$j_0$ in the long exact sequence
\begin{equation}\label{les_e}
0 \lra \Ga_-(b;TY) \stackrel{i_0}{\lra} \Ga_-(b;T\Pf)
\stackrel{j_0}{\lra} \Ga_-(b;\L) \stackrel{\d_0}{\lra}  H_J^1(b;TY) \lra 0,
\end{equation}
is injective on~${\cal N}^{\ka}{\T}$.
We denote the image bundle of~$j_0$ by $\V_+\!\subset\!\V_1^d$ and 
its $L^2$-orthogonal complement in $\V_1^d$ by~$\V_-$. 
Let $\pi_+$ and $\pi_-$ be the corresponding projection maps.

\begin{lmm}
\label{maincontr_lmm3}
Suppose ${\T}$ is a bubble type as in Lemma~\ref{maincontr_lmm1}
and $\ka\!\in\!\S_0(Y;J)$.
For every precompact open subset $K$ of $\U_{\T;\ka}$,
there~exist\\
${}\quad$ (a) $\de_K,\de_K'\!\in\!\R^+$
and an open neighborhood~$U_K$ of~$K$ in~$\X_1(\Pf,d)$;\\
${}\quad$ (b) an orientation-preserving diffeomorphism
$\phi_{\T,\ka}\!:{\cal N}^{\ka}{\T}_{\de_K'}\!\times_K\!
{\cal FT}_{\de_K}^{\eset}\!\lra\!\M_1^0(\Pf,d;J)\cap U_K$;\\
${}\quad$ (c) a lift $\ti\phi_{\T,\ka}\!: 
\pi_{{\cal N}^{\ka}{\T}\oplus{\cal FT}}^*\V_1^d\!\lra\!\V_1^d$ 
of $\phi_{\T,\ka}$ to a vector-bundle isomorphism;\\
with the following property.
If $\vt_t\!\in\!\Ga(\M_1^0(\Pf,d;J)\!\cap\!U_K;\V_1^d)$
is a family of smooth sections such that for some $C\!\in\!\R^+$,
\begin{equation*}\begin{split}
|\vt_t(u)|&\le C\\
\big|\ti{\phi}_{\T,\ka}^{-1}\vt_t\phi_{\T,\ka}(X,\ups)-
\ti{\phi}_{\T,\ka}^{-1}\vt_t\phi_{\T,\ka}(X',\ups)\big|
&\le C|X\!-\!X'|
\end{split}
\qquad\forall X,X'\!\in\!{\cal N}^{\ka}{\T}_{\de_K'},~
\ups\!\in\!{\cal FT}_{\de_K}^{\eset}\big|_K,
\end{equation*}
then  there exists $\ep\!\in\!\R^+$ such that for all 
$t\!\in\![0,\ep)$, $b\!\in\!K$, and $\ups\!\in\!{\cal FT}_{\de_K}|_b$,
the equation 
$$\pi_+\ti{\phi}_{\T,\ka}^{-1}
\big(\{s_1^d\!+\!t\vt_t\}\phi_{\T,\ka}(b;X,\ups)\big)=0$$
has a unique solution $X\!=\!X_t(\ups)\!\in\!{\cal N}^{\ka}{\T}_{\de_K'}|_b$.
Furthermore,
$$\lim_{t\lra0,\ups\lra b}t^{-1}\pi_-\ti{\phi}_{\T,\ka}^{-1}
\big(s_1^d\phi_{\T,\ka}(b;X_t(\ups),\ups)\big)=0
\qquad\forall\, b\!\in\!K.$$
\end{lmm}

\noindent
{\it Proof:} 
(1) The desired maps $\phi_{\T,\ka}$ and $\ti{\phi}_{\T,\ka}$
are simply the compositions 
$\phi_{\T,0}\!\circ\!\ti{\vph}_{\T,\ka}$ and 
$\ti{\phi}_{\T,0}\!\circ\!\ti{\vph}_{\T,\ka}$, respectively.
For each $b\!\in\!K$, $X\!\in\!{\cal N}^{\ka}{\T}|_b$, and 
$\ups\!\in\!{\cal FT}^{\eset}|_b$
sufficiently small, we define $\ti{N}_s(X)$ and $\ti{N}_s'(X,\ups)$
in $\V_1^d|_b$ by
\begin{alignat}{1}
\label{maincontr_lmm3e1a}
&\ti{\vph}_{\T,\ka}^{-1}\big(s_1^d\vph_{\T,\ka}(b;X)\big)
=s_1^d(b)+j_0X+\ti{N}_s(X)=j_0X+\ti{N}_s(X);\\
\label{maincontr_lmm3e1b}
&\ti{\phi}_{\T,\ka}^{-1}\big(s_1^d\phi_{\T,\ka}(b;X,\ups)\big)
=\ti{\vph}_{\T,\ka}^{-1}s_1^d\vph_{\T,\ka}(b;X)+\ti{N}_s'(X,\ups).
\end{alignat}
Since $j_0$ is the derivative of $s_1^d$ on $\U_{\T}(\Pf;J)$,
for some $C_K\!\in\!\Bbb{R}^+$,
\begin{equation}\label{maincontr_lmm3e2a}
\ti{N}_s(0)\!=\!0,\quad
\big|\ti{N}_s(X)-\ti{N}_s(X')\big|<
C_K\big(|X|\!+\!|X'|\big)|X\!-\!X'|
\quad\forall X,X'\!\in\!{\cal N}^{\ka}{\T}_{\de_K'}|_K.
\end{equation}
For $\ti{N}_s'(\cdot,\cdot)$, we similarly have
\begin{equation}\label{maincontr_lmm3e2b}\begin{split}
\big|\ti{N}_s'(X,\ups)\big|&\le C_K|\ups|^{1/p},\\
\big|\ti{N}_s'(X,\ups)\!-\!\ti{N}_s'(X',\ups)\big| &<
C_K|\ups|^{1/p}|X\!-\!X'|
\end{split}
\qquad\forall X,X'\!\in\!{\cal N}^{\ka}{\T}_{\de_K'}|_K,
~\ups\!\in\!{\cal FT}_{\de_K}^{\eset}|_K.
\end{equation}
The first estimate above is clear from (2) of Lemma~\ref{maincontr_lmm1}.
The second bound follows from the analogous bound on the behavior 
of the vector field~$\xi_0(\ups)$ of~(2) of Lemma~\ref{maincontr_lmm1};
see Subsection~4.2 in~\cite{gluing}.\\
(2) If $\vt_t$ is a family of smooth sections as in the statement of the lemma,
by \e_ref{maincontr_lmm3e1a} and~\e_ref{maincontr_lmm3e1b},
\begin{gather}\label{maincontr_lmm3e5}
\pi_+\ti{\phi}_{\T,\ka}^{-1}
\big(\{s_1^d\!+\!t\vt_t\}\phi_{\T,\ka}(b;X,\ups)\big)
=j_0X+\pi_+\ti{N}_s(X)+\pi_+\ti{N}_s'(X,\ups)+t\pi_+\ti{\vt}_t(X,\ups),\\
\hbox{where}\qquad
\ti{\vt}_t(X,\ups)=\ti{\phi}_{\T,\ka}^{-1}
\big(\vt_t\phi_{\T,\ka}(b;X,\ups)).\notag
\end{gather}
By \e_ref{maincontr_lmm3e2a}-\e_ref{maincontr_lmm3e5} and the Contraction Principle,
there exist $\de,\de'\!\in\!\R^+$, dependent on~$j_0$ and~$C_K$,
and $\ep,C'\!\in\!\R^+$, dependent on~$j_0$, $C_K$, and~$C$, such that
for all 
$t\!\in\![0,\ep)$, $b\!\in\!K$, and $\ups\!\in\!{\cal FT}_{\de}^{\eset}|_b$,
the equation 
$$\pi_+\ti{\phi}_{\T,\ka}^{-1}
\big(\{s_1^d\!+\!t\vt_t\}\phi_{\T,\ka}(b;X,\ups)\big)=0$$
has a unique solution $X\!=\!X_t(\ups)\!\in\!{\cal N}^{\ka}{\T}_{\de'}|_b$.
Furthermore,
\begin{equation}\label{maincontr_lmm3e7}
\big|X_0(\ups)\big|\le C'|\ups|^{1/p}
\quad\hbox{and}\quad
\big|X_t(\ups)\!-\!X_0(\ups)\big|\le C't.
\end{equation}
(3) By \e_ref{maincontr_lmm3e1a}-\e_ref{maincontr_lmm3e2b} and
\e_ref{maincontr_lmm3e7},
\begin{equation}\label{maincontr_lmm3e9a}
\big|\pi_-\ti{\phi}_{\T,\ka}^{-1}
\big(s_1^d\phi_{\T,\ka}(b;X_t(\ups),\ups)\big)-
\pi_-\ti{\phi}_{\T,\ka}^{-1}
\big(s_1^d\phi_{\T,\ka}(b;X_0(\ups),\ups)\big)\le
C''\big(t\!+\!|\ups|^{1/p}\big)t.
\end{equation}
On the other hand, as can be seen from Lemma~\ref{maincontr_lmm4} below,
\begin{equation}\label{maincontr_lmm3e9b}
\pi_-\ti{\phi}_{\T,\ka}^{-1}
\big(s_1^d\phi_{\T,\ka}(b;X_0(\ups),\ups)\big)=0
\end{equation}
for all $\ups\!\in\!{\cal FT}^{\eset}|_K$ sufficiently small.
The last claim of the lemma follows from \e_ref{maincontr_lmm3e9a}
and~\e_ref{maincontr_lmm3e9b}.\\

\noindent
{\it Proof of Proposition~\ref{maincontr_prp2}:} 
(1) By Corollary~\ref{maincontr_lmm2crl}, if $t\!\in\!\R^+$ 
and $U(K)$ are sufficiently small,
there is a one-to-one correspondence between 
$\{\bpar_J\!+\!t\nu\}^{-1}(0)\!\cap\!U(K)$ and the~set
$$\big\{u\!\in\!\M_1^0(\Pf,d;J)\!\cap\!U(K)\!:
\{s_1^d\!+\!t\vt_t\}(u)\!=\!0\big\},$$
where $\vt_t\!\in\!\Ga(\M_1^0(\Pf,d;J)\!\cap\!U_K;\V_1^d)$
is a family of smooth sections as in Lemma~\ref{maincontr_lmm3}.
In addition,
$$\lim_{\ups\lra b,t\lra0}\vt_t(\ups)=[\{\na s\}P_b\nu(b)]
\qquad\forall\, b\!\in\!K.$$
The homomorphism~$\d_0$ in the long exact sequence~\e_ref{les_e}
restricts to an isomorphism on~$\V_-$ and vanishes on~$\V_+$.
By definition of~$\d_0$ and~$P_b$,
$$\d_0([\{\na s\}P_b\nu(b)]) = \pi_{\ka,\nu}^0(b)
\qquad\forall\, b\!\in\!\U_{\T;\ka}.$$
Thus, by Lemma~\ref{maincontr_lmm3}, 
$$\{\bpar_J\!+\!t\nu\}^{-1}(0)\!\cap\!U(K)=\eset
\qquad\hbox{if}\qquad
\pi_{\ka,\nu}^{0~-1}(0)\!\cap\!K=\eset.$$
The case $|I|\!>\!1$ or $\aleph\!\neq\!\eset$ of Proposition~\ref{maincontr_prp2} 
now follows from the assumption~($\nu3a$).\\
(2) If $|I|\!=\!1$ and $\aleph\!=\!\eset$,  
by the assumption~($\nu3a$) and Lemma~\ref{maincontr_lmm3}, 
the section $\bpar_J\!+\!t\nu$ is transverse to the zero set on $U(K)$ and 
$$^{\pm}\!\big|\{\bpar_J\!+\!t\nu\}^{-1}\!\cap\!U(K)\big|=
~ ^{\pm}\!\big|\pi_{\ka,\nu}^{0~-1}(0)\!\cap\! K\big|.$$
Since $\pi_{\ka,\nu}^{0~-1}(0)\!\subset\!\U_{\T;\ka}=\M_1^0(\ka,d/d_{\ka})$,
$$^{\pm}\!\big|\{\bpar_J\!+\!t\nu\}^{-1}\!\cap\!U(K)\big|
=~^{\pm}\!\big|\pi_{\ka,\nu}^{0~-1}(0)\big|
=\blr{e(\W_{\ka,d/d_{\ka}}^{1,0}),\big[\ov\M_1^0(\ka,d/d_{\ka})\big]},$$
provided $\pi_{\ka,\nu}^{0~-1}(0)\!\subset\!K$ and $t$ and $U(K)$ are sufficiently small.\\

\noindent
We conclude this subsection with Lemma~\ref{maincontr_lmm4},
which was used in Lemma~\ref{maincontr_lmm3}.

\begin{lmm}
\label{maincontr_lmm4}
Suppose $\T$ is a bubble type as in Lemma~\ref{maincontr_lmm1} and $\ka\!\in\!\S_0(Y,J)$.
For every precompact open subset $K$ of $\U_{\T;\ka}$,
there exist $\de\!\in\!\Bbb{R}^+$, an open neighborhood~$U$ of~$K$ in~$\X_1(\Pf,d)$,
and an orientation-preserving diffeomorphism
$$\phi_{\T,\ka}'\!:{\cal FT}_{\de}^{\eset}|_K\lra
\M_1^0(\ka,d/d_{\ka})\!\cap\!U \subset\M_1^0(\Pf,d;J).$$
\end{lmm}

\noindent
{\it Proof:} If $\T\!=\!(I,\aleph;\under{d})$, 
\begin{gather*}
\U_{\T;\ka}= \U_{\T'}(\ka;J_0)\approx
\U_{\T'}(\Bbb{P}^1;J_0)
\quad\hbox{and}\quad
{\cal FT}|_{\U_{\T;\ka}}={\cal FT}'\lra \U_{\T'}(\Bbb{P}^1;J_0),\\
\hbox{where}\qquad
{\T}'\!=\!(I,\aleph;\under{d}')\quad
\quad\hbox{and}\quad d_i'\!=\!d_i/d_{\ka}.
\end{gather*}
Thus, Lemma~\ref{maincontr_lmm4} is the $\Bbb{P}^1$-analogue 
of the $t\!=\!0$ case of (4) of Lemma~\ref{maincontr_lmm1}.

\subsection{Proof of Proposition~\ref{bdcontr_prp}}
\label{bdcontr_subs}

\noindent
Proposition~\ref{bdcontr_prp}
follows immediately from Proposition~\ref{bdcontr_prp2}.

\begin{prp}
\label{bdcontr_prp2}
Suppose $d$, $Y$, and $J$ are as in Proposition~\ref{maincontr_prp},
$\nu\!\in\!\ti{\cal A}_1^d(\bpar;J)$ is a generic perturbation
of the $\bpar_J$-operator on $\X_1(\Pf,d)$,
$\ka\!\in\!\S_0(Y;J)$, and
$\T\!=\!(I,\aleph;\under{d})$ is a bubble type such that 
$\sum_{i\in I}\!d_i\!=\!d$ and $d_i\!=\!0$ for all minimal elements~$i$ of~$I$.
If $|\hat{I}|\!>\!1$ or $\aleph\!\neq\!\eset$, 
for every compact subset~$K$ of $\U_{\T;\ka}$,
there exist $\ep_{\nu}(K)\!\in\!\Bbb{R}^+$ and 
an open neighborhood~$U_{\nu}(K)$ of~$K$ in $\X_1(Y,d)$ such~that
$$\{\bpar_J\!+\!t\nu\}^{-1}(0)\!\cap\! U_{\nu}(K)=\eset
\qquad\forall\, t\!\in\!(0,\ep_{\nu}(K)).$$
If $|\hat{I}|\!=\!1$ and $\aleph\!=\!\eset$, 
for every compact subset~$K$ of $\U_{\T;\ka}$, 
there exist $\ep_{\nu}(K)\!\in\!\R^+$ and an open neighborhood~$U(K)$ 
of~$K$ in $\X_1(Y,d)$ with the following properties:\\
${}\quad$ (a) the section $\bpar_J\!+\!t\nu$ is transverse to the zero set
in $\Ga_1^{0,1}(Y,d;J)$ over $U(K)$ for all $t\!\in\!(0,\ep_{\nu}(K))$;\\
${}\quad$ (b) for every open subset $U$ of $\X_1(Y,d)$,
there exists $\ep(U)\!\in\!(0,\ep_{\nu}(K))$ such that
\begin{gather*}
^{\pm}\!\big|\{\bpar_J\!+\!t\nu\}^{-1}\!\cap\!U\big|=
\frac{d/d_{\ka}}{12}\,
\blr{e(\W_{\ka,d/d_{\ka}}^0),\big[\ov\M_0(\ka,d/d_{\ka})\big]}
\qquad\hbox{if}\\
\ti\pi_{\nu,\ka}^{1~-1}(0)\!\subset\!K\!\subset\!U\!\subset\!U(K),
~t\!\in\!(0,\ep(U)).
\end{gather*}
\end{prp}

\noindent
In simpler words, none of the strata of $\ov\M_1^q(\ka,d/d_{\ka})$ with $q\!\ge\!2$
contributes to the number~$N_1(d)$.
Neither does any of the boundary strata of $\ov\M_1^1(\ka,d/d_{\ka})$.
On the other hand, $\M_1^1(\ka,d/d_{\ka})$ contributes the euler class
of the bundle~$\wt{\cal W}_{\ka,d/d_{\ka}}^{1,1}$;
see Subsection~\ref{setup_subs}.\\

\noindent
We will proceed similarly to Subsection~\ref{maincontr_subs}, but
run the gluing construction in~$Y$, instead of~$\Pf$, and
make use of the assumption~($J_Y2$) from the start.
We will also use the family of metrics on~$\Bbb{P}^1$ provided by 
Lemma~2.1 in~\cite{g2n2and3}, which we now restate:

\begin{lmm}
\label{flat_metrics_lmm}
There exist $r_{\Bbb{P}^1}\!>\!0$ and a smooth family of Kahler metrics
$\{g_{\Bbb{P}^1,q}\!: q\!\in\!\Bbb{P}^1\}$ on~$\Bbb{P}^1$ with the following property.
If $B_q(q',r)\!\subset\!\Bbb{P}^1$ denotes the $g_{\Bbb{P}^1,q}$-geodesic ball about~$q'$,
the triple $(B_q(q,r_{\Bbb{P}^1}),J_0,g_{\Bbb{P}^1,q})$ is isomorphic
to a ball in $\C^1$ for \hbox{all $q\!\in\!\Bbb{P}^1$}.
\end{lmm}

\noindent
In this case, the operators $D_{J,b}|_{\Ga(\Si_b;u_b^*TY)}$ are not surjective 
for  $b\!\in\!\U_{\T;\ka}^{(0)}$,
where $\U_{\T;\ka}^{(0)}$ is the preimage of $\U_{\T;\ka}$
under the quotient projection map 
$\U_{\T}^{(0)}(\Pf;J)\!\lra\!\U_{\T}(\Pf;J)$.
Thus, in contrast to the case of Lemma~\ref{maincontr_lmm1},
we encounter an obstruction bundle in trying to solve 
the $\bar{\partial}_J$-equation near~$\U_{\T;\ka}$,
as in Subsections~3.3-3.5 of~\cite{gluing}.
Subsections~3.3-3.5 in~\cite{g2n2and3} 
describe a special case of an analogous construction in circumstances
similar to the present situation.\\

\noindent
First, we describe a convenient ``exponential" map for $Y$ defined 
on a neighborhood of each smooth curve \hbox{$\ka\!\in\!\S_0(Y;J)$}.
We identify the rational curve $\ka$ with~$\Bbb{P}^1$.
For each $b\!\in\!\U_{\T;\ka}$, let $g_{Y,b}$ be 
a $J$-compatible extension of the metric  $g_{\ka,b}\!\equiv\!g_{\Bbb{P}^1,\ev_P(b)}$ on~$\ka$
provided by Lemma~\ref{flat_metrics_lmm} to a Riemannian metric on
a neighborhood of $\ka$ in~$Y$.
We identify the normal bundle $N_Y\ka$ 
of $\ka$ in~$Y$ with the $g_{Y,b}$-orthogonal complement of $T\ka$ in~$TY|_{\ka}$.
Let
$$\exp_b\!:T\ka\lra\ka \qquad\hbox{and}\qquad
\widetilde{\exp}_b\!:\pi_{T\ka}^*N_Y\ka\lra N_Y\ka$$
be the exponential map with respect to the metric $g_{\ka,b}$ and 
a lift of $\exp_b$ to a vector bundle homomorphism restricting to the identity over~$\ka$.
For example, $\wt{\exp}_b$ can be taken to be the $g_{Y,b}$-parallel transport
along the $g_{\ka,b}$-geodesics.
For each $q\!\in\!\ka$ and $\xi\!\in\!T_qY$ sufficiently small, let
$$\exp_b\xi=\exp_{g_{Y,b},{\exp_{b,q}\xi_-}}\big(\widetilde{\exp}_{b,\xi_-}\xi_+\big)
\qquad\hbox{if}\quad
\xi=\xi_-\!+\!\xi_+\in T\ka\!\oplus\!N_Y\ka=TY,$$
where $\exp_{g_{Y,b}}$ is the exponential map for the metric $g_{Y,b}$.
One useful property of this ``exponential" map is that
$\exp_b\xi\!\in\!\ka$ if $\xi\!\in\!T\ka\!\subset\!TY$.\\

\noindent
For each element $b\!=\!(\Si_b,u_b)$ of $\U_{\T;\ka}^{(0)}$,
we identify the cokernel $H_J^1(b;TY)$ of the operator
$$D_{J,b}\!:\Ga(b;TY)\lra\Ga^{0,1}(b;TY;J)$$
with the space $\Ga_-^{0,1}(b;TY)$ of $(J,j)$-antilinear 
$u_b^*TY$-valued harmonic forms on~$\Si_b$.
The elements of $\Ga_-^{0,1}(b;TY)$ may have simple poles at the nodes of~$\Si_b$
with the residues adding up to zero at each node.
If ${\cal H}_{b,P}$ denotes the one-dimensional vector space of
harmonic antilinear differentials on the principal component(s) $\Si_{b,P}$ of~$\Si_b$,
$$\Ga_-^{0,1}(b;TY)=\Ga_-^{0,1}(b;T\ka)\oplus \Ga_-^{0,1}(b;N_Y\ka)
={\cal H}_{b,P}\!\otimes\!T_{\ev_P(b)}\ka\oplus \Ga_-^{0,1}(b;N_Y\ka).$$
This decomposition is $L^2$-orthogonal. 
Furthermore, $\Ga_-^{0,1}(b;N_Y\ka)$ is isomorphic to the cokernel $H_J^1(b;N_Y\ka)$
of the operator
$$D_{J,b}^{\perp}\!:\Ga(b;N_Y\ka)\lra\Ga^{0,1}(b;N_Y\ka;J)$$
induced by the operator $D_{J,b}$ via the quotient projection map
$$\pi_{\ka}^{\perp}:TY|_{\ka}\lra N_Y\ka=TY|_{\ka}\big/T\ka.$$
We note that if $\aleph\!=\!\eset$ and $|\hat{I}|\!=\!1$,
$\Ga_-^{0,1}(b;TY)$ is a subspace of $\Ga^{0,1}(b;TY;J)$.\\

\noindent
We are now ready to proceed with the pregluing construction.
For each sufficiently small element $\ups\!=\!(b,v)$ of $\ti{\cal F}{\T}^{\eset}$, let 
$$b(\ups)=\big(\Si_{\ups},j_{\ups};u_{\ups}\big)$$
be the corresponding approximately holomorphic stable map,
as in Subsection~\ref{maincontr_subs}.
In the present case, the linearization~$D_{J,b}$ of the $\bar{\partial}_J$-operator 
at~$b$ is not surjective.
Thus, the linearization $D_{J,\ups}$ of the $\bar{\partial}_J$-operator at $b(\ups)$,
defined via the Levi-Civita connection of the metric~$\ti{g}_{Y,b}$,
is not uniformly surjective.
An approximate cokernel of~$D_{J,b}$ is given~by
\begin{equation}\label{cokerbundle_e4}
\Ga_-^{0,1}(\ups;TY)=\Ga_-^{0,1}(\ups;T\ka)\oplus\Ga_-^{0,1}(\ups;N_Y\ka),
\end{equation}
with the vector spaces $\Ga_-^{0,1}(\ups;T\ka)$ and $\Ga_-^{0,1}(\ups;N_Y\ka)$
explicitly describable from $\Ga_-^{0,1}(b;T\ka)$ and $\Ga_-^{0,1}(b;N_Y\ka)$,
respectively, via the basic gluing map $q_{\ups}\!:\Si_{\ups}\!\lra\!\Si_b$.
In fact, we can simply take
\begin{equation}\label{cokerbundle_e5}
\Ga_-^{0,1}(\ups;N_Y\ka)=\big\{q_{\ups}^*\eta\!:\eta\!\in\!\Ga_-^{0,1}(b;N_Y\ka)\big\}.
\end{equation}
While we can define the space $\Ga_-^{0,1}(\ups;T\ka)$ in the same way
from $\Ga_-^{0,1}(b;T\ka)$,
in the $\aleph\!=\!\eset$, $|\hat{I}|\!=\!1$ case
it is more convenient to take
$$\Ga_-^{0,1}(\ups;T\ka)=
\big\{R_{\ups}\eta\!: \eta\!\in\!\Ga_-^{0,1}(b;T\ka)\big\},$$
where $R_{\ups}\eta$ is a smooth extension of $\eta$ such that 
$R_{\ups}\eta$ is harmonic on the neck attaching the only bubble $\Si_{b,h}$
of $\Si_b$ and below a small collar of the neck and vanishes past a slighter larger collar.
For an explicit description of~$R_{\ups}\eta$, see 
the construction at the beginning of Subsection~2.2 in~\cite{g2n2and3}.
We observe that
\begin{equation}\label{cokerbundle_e9}\begin{split}
&\llrr{\eta_{\ka},\ti\eta}_{\ups,2}=0, \qquad
\llrr{\bpar_Ju_{\ups},\ti\eta}_{\ups,2}=0,\qquad
\llrr{D_{J,\ups}\xi,\ti\eta}_{\ups,2}=0,\\
&\qquad\forall~\xi\!\in\!\Ga(\ups;T\ka),~
\eta_{\ka}\!\in\!\Ga_-^{0,1}(\ups;T\ka),~
\ti\eta\!\in\!\Ga_-^{0,1}(\ups;N_Y\ka),
\end{split}\end{equation}
where $\llrr{\cdot,\cdot}_{\ups,2}$ is the $L^2$-inner product
of the metric $g_{Y,b}$ on $Y$.
This inner-product is independent of the choice of a metric on $\Si_{\ups}$
compatible with the complex structure~$j_{\ups}$ on~$\Si_{\ups}$,
though we will always view $\Si_{\ups}$ as carrying the metric~$g_{\ups}$
induced by the pregluing construction.
If $\aleph\!=\!\eset$ and $|\hat{I}|\!=\!1$, $\Ga_-^{0,1}(\ups;TY)$
is a subspace of $\Ga^{0,1}(\ups;TY;J)$, and we denote its $L^2$-orthogonal complement
by~$\Ga_+^{0,1}(\ups;TY)$. Let
$$\pi^{0,1}_{\ups;\ka},\ti\pi^{0,1}_{\ups},\pi^{0,1}_{\ups;+}\!:
\Ga^{0,1}(\ups;TY;J)\lra \Ga^{0,1}_-(\ups;T\ka),\Ga^{0,1}_-(\ups;N_Y\ka),
\Ga^{0,1}_+(\ups;TY)$$
be the $L^2$-projection maps.\\

\noindent
As in Subsection~\ref{maincontr_subs},
if $\ups$ is sufficiently small, we can also obtain a decomposition
\begin{equation}\label{kerbundle_e1}
\Ga(\ups;TY)=\Ga_-(\ups;TY)\oplus\ti\Ga_+(\ups;TY)
\end{equation}
such that the linear operator 
$$D_{J,\ups}\!:\Ga_+(\ups;TY)\lra \Ga^{0,1}(\ups;TY;J)$$ 
is injective,  while 
$$\Ga_-(\ups;TY)=\big\{\xi\circ q_{\ups}\!: \xi\!\in\!\Ga_-(b;TY)\big\}.$$
In this case, $D_{J,\ups}$ denotes the linearization of the $\bar{\partial}_J$-operator
at $b(\ups)$ with the respect to the ``exponential" map chosen above.
In \e_ref{kerbundle_e1}, we can take the space $\ti\Ga_+(\ups;TY)$
to be the $L^2$-orthogonal complement of $\Ga_-(\ups;TY)$, and
we do so unless $\aleph\!=\!\eset$ and $|\hat{I}|\!=\!1$.
If $\aleph\!=\!\eset$ and $|\hat{I}|\!=\!1$,
we can choose $\ti\Ga_+(\ups;TY)$ in such a way that
\begin{equation}\label{kerbundle_e2}
\llrr{D_{J,\ups}\xi,\eta}_{\ups,2}=0
\qquad\forall\, \xi\!\in\!\ti\Ga_+(\ups;TY)\!\cap\!\Ga(\ups;T\ka),
~\eta\!\in\!\Ga^{0,1}_-(\ups;TY),
\end{equation}
the operator 
$$D_{J,\ups}\!:
\ti\Ga_+(\ups;T\ka)\!\equiv\!\ti\Ga_+(\ups;TY)\!\cap\!\Ga(\ups;T\ka)
\lra 
\Ga_+^{0,1}(\ups;T\ka)\!\equiv\!\Ga_+^{0,1}(\ups;TY)\!\cap\!\Ga^{0,1}(\ups;T\ka)$$
is an isomorphism, and the intersection of $\ti\Ga_+(\ups;TY)$ with 
the $L^2$-orthogonal complement of $\Ga_-(\ups;TY)$ has codimension one in both spaces.
The subspace $\ti\Ga_+(\ups;TY)$ of $\Ga(\ups;TY)$
is constructed by restricting the procedure described
in Subsection~2.3 of~\cite{g2n2and3} to 
the line ${\cal H}_{b,\aleph}\!\otimes\!T_{\ev_P(b)}\ka$.\\

\noindent
Similarly, let $\Ga_+(\ups;N_Y{\ka})$ be the $L^2$-orthogonal complement of
$$\Ga_-(\ups;N_Y\ka)=\big\{\xi\circ q_{\ups}\!: \xi\!\in\!\Ga_-(b;N_{\ka}Y)\big\}$$
in $\Ga(\ups;N_Y{\ka})\!\equiv\!L^p_1(\Si_{\ups};u_{\ups}^*N_Y\ka)$.
If $\ups$ is sufficiently small,  the linear operator 
$$D_{J,\ups}^{\perp}\!:\Ga_+(\ups;N_Y\ka)\lra \Ga^{0,1}(\ups;N_Y\ka;J)$$ 
is injective. 
The key properties of this setup are described in Lemma~\ref{bdcontr_lmm1}:

\begin{lmm}
\label{bdcontr_lmm1}
If ${\T}$, $\nu$, and $\ka$ are as in Proposition~\ref{bdcontr_prp2},
for every precompact open subset $K$ of $\U_{\T;\ka}$,
there exist $\de_K,C_K\!\in\!\R^+$ and an open neighborhood $U_K$ of $K$ 
in $\X_1(Y,d)$ with the following properties:\\
(1) for every $[\ti{b}]\!\in\X_1^0(Y,d)\cap U_K$, 
there exist $\ups\!\in\!\ti{\cal F}{\T}_{\de_K}^{\eset}|_{K^{(0)}}$ and 
$\ze\!\in\!\ti\Ga_+(\ups;TY)$ such that $\|\ze\|_{\ups,p,1}\!<\!\de_K$
and $[\exp_{b(\ups)}\!\ze]\!=\![\ti{b}]$, and the pair
$(b,\ze)$ is unique up to the action of the group 
$\hbox{Aut}({\T})\!\propto\!(S^1)^{\hat{I}}$;\\
(2) for all $\ups\!=\!(b,v)\!\in\!\ti{\cal F}{\T}^{\eset}_{\de_K}|_{K^{(0)}}$,
\begin{gather*}
\|\bar{\partial}_Ju_{\ups}\|_{\ups,p}\le  C_K|\ups|^{1/p};\\
C_K^{-1}\|\ze\|_{\ups,p,1}\le \|D_{J,\ups}\ze\|_{\ups,p} 
\le C_K\|\ze\|_{\ups,p,1}   \quad\forall\, \ze\!\in\!\ti\Ga_+(\ups;TY);\\
C_K^{-1}\|\xi\|_{\ups,p,1}\le \|D_{J,\ups}^{\perp}\xi\|_{\ups,p} 
\le C_K\|\xi\|_{\ups,p,1}   \quad\forall\, \xi\!\in\!\Ga_+(\ups;N_{\ka}Y);
\end{gather*}
(3) for all $\ups\!=\!(b,v)\!\in\!\ti{\cal F}{\T}^{\eset}_{\de_K}|_{K^{(0)}}$,
$\xi\!\in\!\Ga(\ups;N_{\ka}Y)$, and $\eta\!\in\!\Ga^{0,1}_-(\ups;N_{\ka}Y)$,
$$\big|\llrr{D_{J,\ups}^{\perp}\xi,\eta}_{\ups,2}\big|
\le C_K|\ups|^{1/p}\|\xi\|_{\ups,p,1}\|\eta\|_{\ups,1}.$$
\end{lmm}

\noindent
In the first claim of Lemma~\ref{bdcontr_lmm1},
$$\exp_{b(\ups)}\!\ze=\big(\Si_{\ups},j_{\ups};\exp_{b,u_{\ups}}\!\ze\big).$$
This statement is a variation on (2) of Lemma~\ref{g1comp-reg1_lmm3} in~\cite{g1comp}
and holds for the same reasons.
The first estimate in~(2) and~(3) of Lemma~\ref{bdcontr_lmm1} can be obtained 
by direct computations.
The two remaining estimates are proved analogously to the corresponding estimates
of Lemma~\ref{maincontr_lmm1}.

\begin{crl}
\label{bdcontr_crl1}
Suppose $\nu$, ${\T}$, and $\ka$ are as in Proposition~\ref{bdcontr_prp2}.
If $q\!\in\!\Z^+$ and $K$ is a compact subset of 
$\U_{\T;\ka}\!\subset\!\ov\M_1^q(\ka,d/d_{\ka})$
such that 
$$\ti{\pi}_{\ka,\nu}^{q~-1}(0)\!\cap\!K=\eset,$$
then there exist $\ep_{\nu}(K)\!\in\!\R^+$
and an open neighborhood $U_{\nu}(K)$ of $K$ in $\X_1(Y,d)$ such~that
$$\{\bpar_J\!+\!t\nu\}^{-1}(0)\!\cap\! U_{\nu}(K)=\eset
\qquad\forall\, t\!\in\!(0,\ep_{\nu}(K)).$$
\end{crl}

\noindent
{\it Proof:} 
(1) As usually, for all $\ze\!\in\!\Ga(\ups;TY)$ sufficiently small,
$$\Pi_{\ze}^{-1} \{\bpar_J\!+\!t\nu\}\exp_{b(\ups)}\!\ze
=\bpar_Ju_{\ups}+D_{J,\ups}\ze+N_{\ups}\ze+tN_{\nu,\ups}\ze+t\nu|_{u_{\ups}},$$
where $\Pi_{\ze}$ denotes the parallel transport with respect to the Levi-Civita
connection of the metric $g_{Y,b}$ along the geodesics of the map~$\exp_b$.
The nonlinear terms satisfy
\begin{equation}\label{bdcontr_crl1e2}\begin{split}
\|N_{\ups}\ze\!-\!N_{\ups}\ze'\|_{\ups,p}&\le
C_K\big(\|\ze\|_{\ups,p,1}\!+\!\|\ze'\|_{\ups,p,1}\big)\|\ze\!-\!\ze'\|_{\ups,p,1}\\
\|N_{\nu,\ups}\ze\!-\!N_{\nu,\ups}\ze'\|_{\ups,p}&\le C_K\|\ze\!-\!\ze'\|_{\ups,p,1}
\end{split}
\qquad  \forall\, \ze,\ze'\!\in\!\Ga(\ups;TY);
\end{equation}
see Subsection~3.6 in~\cite{gluing} for example.
Our choice of the map $\exp_b$ also implies that
\begin{equation}\label{bdcontr_crl1e3}
N_{\ups}\ze\in\Ga^{0,1}(\ups;T\ka)
\qquad\forall\, \ze\!\in\Ga(\ups;T\ka)\subset\Ga(\ups;TY).
\end{equation}
(2) Suppose $\ups\!=\!(b,\ups)\!\in\!\ti{\cal F}{\T}^{\eset}|_{K^{(0)}}$,
$\ze\!\in\!\Ga_+(\ups;TY)$, and 
\begin{gather}
\{\bpar_J\!+\!t\nu\}\exp_{b(\ups)}\!\ze=0
\qquad\Lra\notag\\
\label{bdcontr_crl1e5}
\bpar_Ju_{\ups}+D_{J,\ups}\ze+N_{\ups}\ze+tN_{\nu,\ups}\ze+t\nu|_{u_{\ups}}=0.
\end{gather}
From (2) of Lemma~\ref{bdcontr_lmm1} and \e_ref{bdcontr_crl1e2},
we then obtain
\begin{equation}\label{bdcontr_crl1e6}
\|\ze\|_{\ups,p,1}\le C_K\big(|\ups|^{1/p}\!+\!t\big).
\end{equation}
On the other hand, applying the projection map $\pi_{\ka}^{\perp}$
to both sides of \e_ref{bdcontr_crl1e5}, we~get
\begin{gather}
\label{bdcontr_crl1e7}
D_{J,\ups}^{\perp}\ze^{\perp}
+N_{\ups}^{\perp}\ze+tN_{\nu,\ups}^{\perp}\ze+t\nu^{\perp}|_{u_{\ups}}=0
\in\Ga^{0,1}(\ups;N_Y\ka;J),
\qquad\hbox{if}\\
\ze=\ze^t+\ze^{\perp}\in  \Ga(\ups;T\ka)\oplus\Ga(\ups;N_Y\ka),\quad
N_{\ups}^{\perp}\ze=\pi_{\ka}^{\perp}N_{\ups}\ze,\quad
N_{\nu,\ups}^{\perp}\ze=\pi_{\ka}^{\perp}N_{\nu,\ups}\ze,\quad
\nu^{\perp}=\pi_{\ka}^{\perp}\nu.\notag
\end{gather}
By \e_ref{bdcontr_crl1e2}, \e_ref{bdcontr_crl1e3}, and \e_ref{bdcontr_crl1e6},
$$\big\|N_{\ups}^{\perp}\ze\big\|_{\ups,p}
=\big\|\pi_{\ka}^{\perp}\big(N_{\ups}(\ze^t\!+\!\ze^{\perp})
-N_{\ups}\ze^t\big)\big\|_{\ups,p}
\le  C_K\big(|\ups|^{1/p}\!+\!t\big)\,\|\ze^{\perp}\|_{\ups,p,1}.$$
Thus, by (2) of Lemma~\ref{bdcontr_lmm1}, \e_ref{bdcontr_crl1e2},
and \e_ref{bdcontr_crl1e7},
\begin{equation}\label{bdcontr_crl1e9}
\|\ze^{\perp}\|_{\ups,p,1}\le C_Kt,
\end{equation}
provided $\de_K$ is sufficiently small.
Combining (3) of Lemma~\ref{bdcontr_lmm1}, \e_ref{bdcontr_crl1e2}, \e_ref{bdcontr_crl1e7},
and~\e_ref{bdcontr_crl1e9}, we obtain
$$\big|\llrr{\nu|_{u_{\ups}},\eta}_{\ups,2}\big|
\le C_K\big(|\ups|^{1/p}\!+\!t\big)\|\eta\|_{\ups,1}
\qquad\forall\,\eta\!\in\!\Ga^{0,1}_-(\ups;N_Y\ka).$$
Since the section $\ti{\pi}_{\ka,\nu}^{q~-1}$ of the bundle
$\wt{\cal W}_{\ka,d/d_{\ka}}^{1,q}$ does not vanish over the compact set~$K$,
it follows that
$$\{\bpar_J\!+\!t\nu\}^{-1}(0)\!\cap\! U_{\nu}(K)=\eset$$
if $t$ and  $U_{\nu}(K)$ are sufficiently small.\\

\noindent
By Lemma~\ref{g1str_lmm2}, the spaces $\ov\M_1^q(\ka,d/d_{\ka})$ 
with $q\!\ge\!2$ are contained in $\ov\M_1^0(\ka,d/d_{\ka})$.
In particular, if ${\T}$ is a bubble type as in the first claim
of Proposition~\ref{bdcontr_prp2},
$$\U_{\T;\ka}\subset 
\big(\ov\M_1^0(\ka,d/d_{\ka})\!-\!\M_1^0(\ka,d/d_{\ka})\big)
\cup \big(\ov\M_1^1(\ka,d/d_{\ka})\!-\!\M_1^1(\ka,d/d_{\ka})\big).$$
Thus, Corollary~\ref{bdcontr_crl1}, along with the regularity assumptions
($\nu3a$) and~($\nu3c$), implies the first claim of Proposition~\ref{bdcontr_prp2}.

\begin{crl}
\label{bdcontr_crl2}
Suppose $\nu$, ${\T}$, and $\ka$ are as in Proposition~\ref{bdcontr_prp2}.
If $\aleph\!=\!\eset$ and  $|\hat{I}|\!=\!1$,
for every compact subset~$K$ of $\U_{\T;\ka}$ containing
$\ti{\pi}_{\nu,\ka}^{1~-1}(0)$, 
there exist $\ep_{\nu}(K)\!\in\!\R^+$ and an open neighborhood~$U(K)$ 
of~$K$ in $\X_1(Y,d)$ with the following properties:\\
${}\quad$ (a) the section $\bpar_J\!+\!t\nu$ is transverse to the zero set
in $\Ga_1^{0,1}(Y,d;J)$ over $U(K)$ for all $t\!\in\!(0,\ep_{\nu}(K))$;\\
${}\quad$ (b) for every open subset $U$ of $\X_1(Y,d)$,
there exists $\ep(U)\!\in\!(0,\ep_{\nu}(K))$ such that
$$^{\pm}\!\big|\{\bpar_J\!+\!t\nu\}^{-1}\!\cap\!U\big|=
\frac{d/d_{\ka}}{12}
\blr{ e(\W_{\ka,d/d_{\ka}}^0),\big[\ov\M_0(\ka,d/d_{\ka})\big]}
\quad\hbox{if}\quad
K\!\subset\!U\!\subset\!U(K), ~t\!\in\!(0,\ep(U)).$$
\end{crl}

\noindent
{\it Proof:} (1) By Corollary~\ref{bdcontr_crl1} and the assumption ($\nu3a$) on $\nu$,
it can be assumed that the compact set~$K$ is disjoint from $\ov\M_1^0(\ka,d/d_{\ka})$.
Thus, if $h$ is the unique element of~$\hat{I}$,
\begin{equation}\label{bdcontr_crl2e1}
\big|{\cal D}_{J,h}^{(1)}\ups\big|\ge C_K|\ups|
\qquad\forall\, [\ups]\!=\![b,v]\!\in\!{\cal FT}|_K;
\end{equation}
see Lemma~\ref{g1str_lmm2}.
By Lemma~\ref{bdcontr_lmm1} and the proof of Corollary~\ref{bdcontr_crl1}, 
we need to determine the number of solutions $[\ups,\ze]$ of the equation
\begin{equation}\label{bdcontr_crl2e3}
\bar{\partial}_Ju_{\ups}+D_{J,\ups}\ze+
N_{\ups}\ze+tN_{\nu,\ups}\ze+t\nu=0,
\quad \ups\!\in\!\ti{\cal F}{\T}^{\eset}_{\de_K}|_{K^{(0)}},~
\ze\!\in\!\ti\Ga_+(\ups;TY), ~\|\ze\|_{\ups,p,1}\!\le\ep_K.
\end{equation}
In this case, $\Ga_-^{0,1}(\ups;TY)$ is a subspace of $\Ga^{0,1}(\ups;TY;J)$,
and the middle estimate in (2) of Lemma~\ref{bdcontr_lmm1} implies that
$$C_K^{-1}\|\ze\|_{\ups,p,1}\le \|\pi^{0,1}_{\ups;+}D_{J,\ups}\ze\|_{\ups,p} 
\le C_K\|\ze\|_{\ups,p,1}   \quad\forall\, \ze\!\in\!\ti\Ga_+(\ups;TY).$$
Thus, the linear operator
$$\pi^{0,1}_{\ups;+}D_{J,\ups}\!: \ti\Ga_+(\ups;TY)\lra\Ga^{0,1}_+(\ups;TY)$$
is an isomorphism.
It then follows from the Contraction Principle, the first estimate in (2)
of Lemma~\ref{bdcontr_lmm1}, and \e_ref{bdcontr_crl1e2} that the equation
$$\pi_{\ups;+}^{0,1}\big(\bar{\partial}_Ju_{\ups}\!+\!D_{J,\ups}\ze\!+\!
N_{\ups}\ze\!+\!tN_{\nu,\ups}\ze\!+\!t\nu\big)=0,
\quad \ze\!\in\!\ti\Ga_+(\ups;TY), ~\|\ze\|_{\ups,p,1}\!\le\ep_K,$$
has a unique solution $\ze_{\ups}\!\in\!\ti\Ga_+(\ups;TY)$,
provided $\ups\!\in\!\ti{\cal F}{\T}^{\eset}_{\de_K}|_{K^{(0)}}$ 
is sufficiently small.
Furthermore,
\begin{equation}\label{bdcontr_crl2e5}
\big\|\ze_{\ups}\|_{\ups,p,1}\le C_K\big(|\ups|^{1/p}\!+\!t\big).
\end{equation}
Thus, the number of solutions $[\ups,\ze]$ of \e_ref{bdcontr_crl2e3} is the same
as the number of solutions of 
\begin{equation}\label{bdcontr_crl2e6}
\Psi_{t\nu}(\ups)\!\equiv t^{-1}\cdot
\pi_{\ups;-}^{0,1}\big(\bpar_Ju_{\ups}\!+\!D_{J,\ups}\ze_{\ups}\!+\!
N_{\ups}\ze_{\ups}\!+\!tN_{\nu,\ups}\ze_{\ups}\!+\!t\nu\big)=0,
\qquad [\ups]\!\in\!{\cal FT}^{\eset}_{\de_K}|_{K^{(0)}},
\end{equation}
where
$$\pi^{0,1}_{\ups;-}\!=\!\pi^{0,1}_{\ups;\ka}\!\oplus\!\ti{\pi}^{0,1}_{\ups}\!:
\Ga^{0,1}(\ups;TY;J)\lra\Ga_-^{0,1}(\ups;T\ka)\oplus\Ga_-^{0,1}(\ups;N_Y\ka)$$
is the $L^2$-projection map.\\
(2) With our choice of the space $\Ga_-^{0,1}(\ups;T\ka)$,
\begin{equation}\label{bdcontr_crl2e9}
\pi_{\ups;\ka}^{0,1}\bpar_Ju_{\ups}=R_{\ups}{\cal D}_{J,h}\ups
\in \Ga_-^{0,1}(\ups;T\ka);
\end{equation}
see Subsection~4.1 in~\cite{g2n2and3}.
Furthermore, 
\begin{equation}\label{bdcontr_crl2e10}
\pi_{\ups;\ka}^{0,1}N_{\ups}\ze=0 \qquad\forall\, \ze\!\in\!\Ga(\ups;T\ka),
\end{equation}
since the supports of all elements of $\eta\!\in\!\Ga_-^{0,1}(\ups;T\ka)$
are disjoint from the support of $N_{\ups}\ze$, for $\ze\!\in\!\Ga(\ups;T\ka)$,
due to our choice of the ``exponential" map.
By \e_ref{bdcontr_crl1e2}, \e_ref{bdcontr_crl2e5}, \e_ref{bdcontr_crl2e9}, 
\e_ref{bdcontr_crl2e10}, and the same argument as in the proof 
of Corollary~\ref{bdcontr_crl1},
\begin{equation}\label{bdcontr_crl2e11}\begin{split}
\big\|\ti\pi^{0,1}_{\ups}\Psi_{t\nu}(\ups)
-R_{\ups}\ti\pi_{\ka,\nu}^1(b)\big\|_{\ups,2}
&\le C_K\big(|\ups|^{1/p}\!+\!t\big)+\|\nu(u_{\ups})\!-\!R_{\ups}\nu(b)\|_{\ups,2}\\
\big\|\pi^{0,1}_{\ups;\ka}\Psi_{t\nu}(\ups)
-R_{\ups}\big(\pi_{\ka,\nu;\ka}^1(b)+t^{-1}{\cal D}_{J,h}\ups\big)\big\|_{\ups,2}
&\le C_K\big(|\ups|^{1/p}\!+\!t\big)+\|\nu(u_{\ups})\!-\!R_{\ups}\nu(b)\|_{\ups,2},
\end{split}\end{equation}
where $R_{\ups}\eta\!=\!q_{\ups}^*\eta$ for $\eta\!\in\!\Ga^{0,1}_-(b;N_Y\ka)$
and 
$$\pi_{\ka,\nu;\ka}^1(b)=\pi_{\ka}\pi_{\ka,\nu}^1(b)
\equiv\pi_{\ka,\nu}^1(b)\!-\!\pi_{\ka}^{\perp}\pi_{\ka,\nu}^1(b)
\in\Ga_-^{0,1}(b;T\ka).$$
Since
$$\lim_{t\lra0,|\ups|\lra0}\!
\big(|\ups|^{1/p}\!+\!t\big)+\|\nu(u_{\ups})\!-\!R_{\ups}\nu(b)\|_{\ups,2}=0,$$
by \e_ref{bdcontr_crl2e1}, \e_ref{bdcontr_crl2e11},
and the same cobordism argument as in Subsection~3.1 of~\cite{g2n2and3},
the number of solutions of \e_ref{bdcontr_crl2e6} is the same as the number of 
solutions of the~system
$$\begin{cases}
\ti\pi_{\ka,\nu}^1(b)=0\in\wt\W_{\ka,d/d_{\ka}}^{1,1}\\
\pi_{\ka,\nu;\ka}^1(b)+{\cal D}_{J,h}\ups=0
\in \pi_P^*\E^*\!\otimes\!\pi_B^*\ev_0^*T\ka
\end{cases}
\qquad [\ups]\!=\![b,v]\in{\cal FT}\lra\ov\M_1^1(k,d/d_{\ka}),$$
if the interior of the compact set $K$ contains 
the finite set $\ti{\pi}_{\nu,\ka}^{1~-1}(0)$.
Since $D_{J,h}$ does not vanish on $\ti{\pi}_{\ka,\nu}^{1~-1}(0)$,
the number of solution of this system~is
\begin{equation*}\begin{split}
^{\pm}\!\big|\ti{\pi}_{\nu,\ka}^{1~-1}(0)\big|
&=\blr{ e(\wt\W_{\ka,d/d_{\ka}}^{1,1}),\big[\ov\M_1^1(\ka,d/d_{\ka})\big]}\\
&=\blr{ e\big(\pi_P^*\E^*\!\otimes\!\pi_B^*\ev_0^*N_Y\ka\big)   
e(\ti\pi_B^*\W_{\ka,d/d_{\ka}}^0),
\big[\ov\cM_{1,1}\!\times\!\ov\M_{0,1}(\ka,d/d_{\ka})\big]}\\
&=\Big(\!-\frac{1}{24}\Big)\big(-2(d/d_{\ka})\big)\,
\blr{ e(\W_{\ka,d/d_{\ka}}^0),\big[\ov\M_0(\ka,d/d_{\ka})\big]},
\end{split}\end{equation*}
as claimed in Proposition~\ref{bdcontr_prp2}.

\section{On the Euler Class of the Cone $\V_1^d\lra\ov\M_1^0(\Pf,d;J)$}
\label{cone_sec}

\subsection{The Structure of the Moduli Spaces $\M_1^0(\P,d;J)$}
\label{g1str_subs}

\noindent
In this section, we prove Proposition~\ref{cone_prp2} by constructing 
a perturbation $\vt$ of the section $s_1^d$ of the cone
$\V_1^d$ over $\ov\M_1^0(\Pf,d;J)$ and counting the number of zeros
of the multisection $s_1^d\!+\!t\vt$ for a small $t\!\in\!\R^+$
that lie near each stratum~of 
\begin{equation}\label{conedecomp_e}
s_1^{d~-1}(0)\cap \ov\M_1^0(\Pf,d;J)
=\ov\M_1^0(Y,d;J)
=\bigsqcup_{\ka\in\S_0(Y;J)}\!\!\!\!\!\!\ov\M_1^0(\ka,d/d_{\ka})
~\sqcup  
\bigsqcup_{\ka\in\S_1(Y;J)}\!\!\!\!\!\!\M_1^0(\ka,d/d_{\ka}).
\end{equation}
Since the single-element orbifold $\M_1^0(\ka,d/d_{\ka})$ is a transverse zero
of $s_1^d$, for  $\ka\!\in\!\S_1(Y;J)$, 
\begin{equation}\label{coneg1contr_e}
\cC_{\M_1^0(\ka,d/d_{\ka})}\big(s_1^d\big)
=\, ^{\pm}\big|\big\{s_1^d\!+\!t\vt\}^{-1}(0)\!\cap\! U_{\ka}\big|
=\, ^{\pm}\big|\M_1^0(\ka,d/d_{\ka})\big|,
\end{equation}
if $U_{\ka}$ is a small neighborhood of $\M_1^0(\ka,d/d_{\ka})$
in $\ov\M_1^0(\Pf,d;J)$.
The second equality in~\e_ref{coneg1contr_e} holds for every multisection~$\vt$
of $\V_1^d$ and every $t\!\in\!\R$ sufficiently small.
Thus, the key to proving Proposition~\ref{cone_prp2} is computing the $s^d_1$-contribution 
from each stratum of the moduli space $\ov\M_1^0(\ka,d/d_{\ka})$.
This is achieved by Proposition~\ref{coneg0contr_prp} and Corollary~\ref{bdcorr_crl}.\\

\noindent
In this subsection, we describe the structure of the moduli space $\ov\M_1^0(\P,d;J)$,
with $J$ sufficiently close to~$J_0$.
Lemmas~\ref{g1str_lmm1} and~\ref{g1str_lmm2} are special cases
of Lemmas~\ref{g1cone-g1mainstr_lmm} and~\ref{g1cone-g1bdstr_lmm}, respectively,
in~\cite{g1cone}.
In turn, the latter two lemmas follow immediately from 
Theorems~\ref{g1comp-reg_thm} and~\ref{g1comp-str_thm} in~\cite{g1comp}.

\begin{lmm}
\label{g1str_lmm1}
If $n,d\!\in\!\Z^+$, there exists $\de_n(d)\!\in\!\R^+$ with the following property.
If $J$ is an almost complex structure on~$\P$, such that $\|J\!-\!J_0\|_{C^1}\!<\!\de_n(d)$, 
and $\T\!=\!(I,\aleph;\under{d})$ is a bubble type such that
$\sum_{i\in I}\!d_i\!=\!d$ and $d_i\!\neq\!0$ for some minimal element $i$ of $I$,
then $\U_{\T}(\P;J)$ is a smooth orbifold,
$$\dim\U_{\T}(\P;J)=2\big(d(n\!+\!1)-|\aleph|-|\hat{I}|\big),
\qquad\hbox{and}\qquad
\U_{\T}(\P;J)\subset\ov\M_1^0(\P,d;J).$$
Furthermore, there exist $\de\!\in\!C(\U_{\T}(\P;J);\R^+)$,
an open neighborhood $U_{\T}$ of $\U_{\T}(\P;J)$ in $\X_1(\P,d)$, 
and an orientation-preserving homeomorphism
$$\phi_{\T}\!:{\cal FT}_{\de}\lra \ov\M_1^0(\P,d;J)\!\cap\!U_{\T},$$
which restricts to a diffeomorphism 
${\cal FT}_{\de}^{\eset}\!\lra\!\M_1^0(\P,d;J)\!\cap\!U_{\T}$.
\end{lmm}

\begin{lmm}
\label{g1str_lmm2}
If $n,d\!\in\!\Z^+$, there exists $\de_n(d)\!\in\!\R^+$ with the following property.
If $J$ is an almost complex structure on~$\P$, such that $\|J\!-\!J_0\|_{C^1}\!<\!\de_n(d)$,
and  ${\T}\!=\!(I,\aleph;\under{d})$ is a bubble type such that
$\sum_{i\in I}\!d_i\!=\!d$ and $d_i\!=\!0$ for all minimal elements $i$ of~$I$, then
$\U_{\T}(\P;J)$ is a smooth orbifold,
\begin{gather*}
\dim\U_{\T}(\P;J)=2\big(d(n\!+\!1)-|\aleph|\!-\!|\hat{I}|+n\big),
~~\hbox{and}~~
\ov\M_1^0(\P,d;J)\cap\U_{\T}(\P;J)=
\U_{\T;1}(\P;J),\\
\hbox{where}\qquad
\U_{\T;1}(\P;J)=
\big\{[b]\!\in\!\U_{\T}(\P;J)\!:\dim_{\C}\hbox{Span}_{(\C,J)}
\{{\cal D}_ib\!:i\!\in\!\chi({\T})\}\!<\!|\chi({\T})|\big\}.
\end{gather*}
The space $\U_{\T;1}(\P;J)$ admits a stratification
by smooth suborbifolds of $\U_{\T}(\P;J)$:
\begin{gather*}
\U_{\T;1}(\P;J)=
\bigsqcup_{m=\max(|\chi({\T})|-n,1)}^{m=|\chi({\T})|}
\!\!\!\!\!\!\!\U_{\T;1}^m(\P;J),
\qquad\hbox{where}\\
\U_{\T;1}^m(\P;J)=\big\{
[b]\!\in\!\U_{\T}(\P;J)\!:\dim_{\C}\hbox{Span}_{(\C,J)}
\{{\cal D}_ib\!:i\!\in\!\chi({\T})\}=|\chi({\T})|\!-\!m\big\},\\
\begin{split}
\dim\U_{\T;1}^m(\P;J)
&= 2\big(d(n\!+\!1)-|\aleph|\!-\!|\hat{I}|
+n+\big(|\chi({\T})|\!-\!n\!-\!m)m\big)\\
&\le \dim\M_1^0(\P,d;J)-2.
\end{split}
\end{gather*}
Furthermore, the space
$${\cal F}^1{\T}^{\eset}\equiv
\big\{[b,\ups]\!\in\!{\cal FT}^{\eset}\!: {\cal D}_{\T}(\ups)\!=\!0\big\}$$
is a smooth oriented suborbifold of ${\cal FT}$.
Finally, there exist $\de\!\in\!C(\U_{\T}(\P;J);\R^+)$,
an open neighborhood $U_{\T}$ of $\U_{\T}(\P;J)$ in $\X_1(\P,d)$, 
and an orientation-preserving diffeomorphism
$$\phi_{\T}\!:
{\cal F}^1{\T}^{\eset}_{\de}\lra \M_1^0(\P,d;J)\!\cap\! U_{\T},$$
which extends to a homeomorphism 
$$\phi_{\T}\!:{\cal F}^1{\T}_{\de}\lra\ov\M_1^0(\P,d;J)\!\cap\! U_{\T},$$
where ${\cal F}^1{\T}$ is the closure of ${\cal F}^1{\T}^{\eset}$
in ${\cal FT}$.
\end{lmm}

\begin{figure}
\begin{pspicture}(-1.1,-1.8)(10,1.25)
\psset{unit=.4cm}
\psellipse(5,-1.5)(1.5,2.5)
\psarc[linewidth=.05](3.2,-1.5){2}{-30}{30}\psarc[linewidth=.05](6.8,-1.5){2}{150}{210}
\pscircle[fillstyle=solid,fillcolor=gray](2.5,-1.5){1}\pscircle*(3.5,-1.5){.2}
\pscircle[fillstyle=solid,fillcolor=gray](.5,-1.5){1}\pscircle*(1.5,-1.5){.2}
\pscircle(7.5,-1.5){1}\pscircle*(6.5,-1.5){.2}
\pscircle[fillstyle=solid,fillcolor=gray](8.91,-.09){1}\pscircle*(8.21,-.79){.2}
\pscircle[fillstyle=solid,fillcolor=gray](8.91,-2.91){1}\pscircle*(8.21,-2.21){.2}
\rput(2.5,0){$h_1$}\rput(.5,0){$h_2$}\rput(7.3,0){$h_3$}
\rput(10.5,0.1){$h_4$}\rput(10.5,-2.9){$h_5$}
\rput(5,-5){\small ``tacnode"}
\pnode(5,-5){A1}\pnode(3.5,-1.5){B1}
\ncarc[nodesep=.35,arcangleA=-25,arcangleB=-15,ncurv=1]{->}{A1}{B1}
\pnode(5,-4.65){A2}\pnode(7.3,-1.5){B2}
\ncarc[nodesep=0,arcangleA=40,arcangleB=30,ncurv=1]{-}{A2}{B2}
\pnode(8,-.95){B2a}\pnode(8.02,-2.02){B2b}
\ncarc[nodesep=0,arcangleA=0,arcangleB=10,ncurv=1]{->}{B2}{B2a}
\ncarc[nodesep=0,arcangleA=0,arcangleB=10,ncurv=1]{->}{B2}{B2b}
\rput(25,-1.5){\begin{tabular}{l}
$\chi({\T})\!=\!\{h_1,h_4,h_5\}$,~
$\rho(\ups)\!=\!(\ups_{h_1},\ups_{h_3}\ups_{h_4},\ups_{h_3}\ups_{h_5})$\\
\\
${\cal F}^1{\T}^{\eset}=\big\{[b;v_1,v_2,v_3,v_4,v_5]\!:v_h\!\in\!\C^*$\\
${}\qquad\qquad~~{\cal D}_{J,h_1}v_{h_1}\!+\!{\cal D}_{J,h_4}v_{h_3}v_{h_4}\!+\!
{\cal D}_{J,h_5}v_{h_3}v_{h_5}\!=\!0\big\}$
\end{tabular}}
\end{pspicture}
\caption{An Illustration of Lemma~\ref{g1str_lmm2}} 
\label{g1bdstr_fig}
\end{figure}

\noindent
We illustrate Lemma~\ref{g1str_lmm2} in Figure~\ref{g1bdstr_fig}.
As before, the shaded discs represent the components of the domain
on which every stable map~$[b]$ in $\U_{\T}(\P;J)$ is non-constant.
The element $[\Si_b,u_b]$ of $\U_{\T}(\P;J)$ is 
in the stable-map closure of $\M_1^0(\P,d;J)$ if and only if
the branches of $u_b(\Si_b)$ corresponding to the attaching nodes on 
the first-level effective bubbles of $[\Si_b,u_b]$ form a generalized tacnode.
In the case of Figure~\ref{g1bdstr_fig}, this means that  either\\
${}\quad$ (a) for some $h\!\in\!\{h_1,h_4,h_5\}$, 
the branch of $u_b|_{\Si_{b,h}}$ at the node~$\i$ has a cusp, or\\
${}\quad$ (b) for all $h\!\in\!\{h_1,h_4,h_5\}$,
the branch of $u_b|_{\Si_{b,h}}$ at the node~$\i$ is smooth, but the dimension\\
${}\qquad~~~$ of the span of the three lines tangent to these
branches is less than three.\\

\noindent
If  $\ka\!\in\!\S_0(Y;J)$, we put
$$\U_{\T;\ka;1}=\U_{\T;\ka} \!\cap\! \U_{\T;1}(\Pf;J)
\subset \ov\M_1^0(\ka,d/d_{\ka}),
\qquad
\U_{\T;\ka;1}^m = \U_{\T;\ka} \!\cap\! \U_{\T;1}^m(\Pf;J)
\subset \U_{\T;\ka;1}.$$
By the $n\!=\!1$ case of Lemma~\ref{g1str_lmm2}, 
$$\U_{\T;\ka;1}=\U_{\T;\ka}
=\U_{\T;\ka;1}^{|\chi({\T})|-1}\cup \U_{\T;\ka;1}^{|\chi({\T})|}
\qquad\hbox{if}\quad |\chi({\T})|\ge2.$$
The last space may be empty.
In particular,
$$\ov\M_1^q(\ka,d/d_{\ka})=
\big(\ov\cM_{1,q}\!\times\!\ov\M_0^q(\ka,d/d_{\ka})\big)/S_q
\subset\ov\M_1^0(\ka,d/d_{\ka})
\quad\hbox{if}~~q\!\ge\!2.$$\\ 

\noindent
Let 
$$\M_{0,1;1}^0(\Pf,d;J)=\big\{[\Bbb{P}^1,u]\!\in\!\M_{0,1}^0(\Pf,d;J)\!:
du|_{\i}\!=\!0\big\}.$$
In other words, $\M_{0,1;1}^0(\Pf,d;J)$ is the subset of $\M_{0,1}^0(\Pf,d;J)$
consisting of the elements $[\Bbb{P}^1,u]$ such that the differential of
$u$ vanishes at the marked point of~$\Bbb{P}^1$, which we always take to be~$\i$.
The image of a generic element in $\M_{0,1;1}^0(\Pf,d;J)$ is a rational curve 
$J$-holomorphic curve in $\Pf$ with a cusp at the image of the marked point.
We denote by $\ov\M_{0,1;1}(\Pf,d;J)$ the closure of $\M_{0,1;1}^0(\Pf,d;J)$
in~$\ov\M_{0,1}(\Pf,d;J)$.
If $\ka\!\in\!\S_0(Y;J)$, we~put
\begin{gather*}
\M_{0,1;1}^0(\ka,d/d_{\ka})=
\M_{0,1;1}^0(\Pf,d;J)\cap\ov\M_{0,1}(\ka,d/d_{\ka}),\\
\ov\M_{0,1;1}(\ka,d/d_{\ka})=
\ov\M_{0,1;1}(\Pf,d;J)\cap\ov\M_{0,1}(\ka,d/d_{\ka}),\\
\M_{1;1}^{1;0}(\ka,d/d_{\ka})=\ov\cM_{1,1}\times\M_{0,1;1}^0(\ka,d/d_{\ka}),
\quad\hbox{and}\quad
\ov\M_{1;1}^1(\ka,d/d_{\ka})=\ov\cM_{1,1}\times\ov\M_{0,1;1}(\ka,d/d_{\ka}).
\end{gather*}
By Lemma~\ref{g1str_lmm2},
$$\ov\M_1^1(\ka,d)\cap \ov\M_1^0(\ka,d)=\ov\M_{1;1}^1(\ka,d)
\qquad\forall\, d\!\in\!\Z^+.$$ 
We note that
\begin{equation}\label{dim_e3}
\dim_{\C}\ov\M_{0,1;1}(\ka,d)=2d-2
\quad\hbox{and}\quad
\dim_{\C}\ov\M_{1;1}^1(\ka,d)=2d-1.
\end{equation}

\subsection{The Structure of the Cone $\V_1^d\lra\ov\M_1^0(\Pf,d;J)$}
\label{conestr_subs}

\noindent
We next describe the structure of the cone $\V_1^d$
near each stratum $\U_{\T}(\Pf;J)$ and $\U_{\T;1}^m(\Pf;J)$
of $\ov\M_1^0(\Pf,d;J)$.
We then state several regularity conditions that we will require the perturbation $\vt$
of $s_1^d$ to satisfy.
The first lemma stated is a special case of Lemma~\ref{g1cone-g1conemainstr_lmm}
in~\cite{g1cone}.

\begin{lmm}
\label{g1conestr_lmm1}
If $d$, $\L$, and $\V_1^d$ are as in Proposition~\ref{cone_prp2},
there exists $\de(d)\!\in\!\R^+$  with the following property.
If $J$ is an almost complex structure on~$\Pf$, such that $\|J\!-\!J_0\|_{C^1}\!<\!\de_n(d)$,
and $\T\!=\!(I,\aleph;\under{d})$ is a bubble type such that 
$\sum_{i\in I}\!d_i\!=\!d$ and $d_i\!\neq\!0$ for some minimal element $i$ of $I$,
then the requirements of Lemma~\ref{g1str_lmm1} are satisfied.
Furthermore, the restriction $\V_1^d\!\lra\!\U_{\T}(\Pf;J)$
is a smooth complex vector orbibundle of rank~$5d$.
Finally, there exists a smooth vector-bundle isomorphism
$$\ti\phi_{\T}\!:
\pi_{{\cal FT}_{\de}^{\eset}}^*\big(\V_1^d\big|_{\U_{\T}(\Pf;J)}\big)
\lra \V_1^d\big|_{\ov\M_1^0(\Pf,d;J)\cap U_{\T}},$$
covering the homeomorphism $\phi_{\T}$ of Lemma~\ref{g1str_lmm1},
such that $\ti{\phi}_{\T}$ is the identity over 
$\U_{\T}(\Pf;J)$ and is smooth over~${\cal FT}_{\de}^{\eset}$.
\end{lmm}

\noindent
For every $\ka\!\in\!\S_0(Y;J)$, 
the family of boundary operators $\d_0$ in the long exact sequence~\e_ref{les_e},
with $b\!\in\!\U_{\T;\ka}^{(0)}$ and $\T$ as in Lemma~\ref{g1conestr_lmm1},
induces a surjective bundle homomorphism
$$\d_{\ka,d/d_{\ka}}^{1,0}\!:\V_1^d\lra \W_{\ka,d/d_{\ka}}^{1,0}$$
over $\M_1^{\{0\}}(\ka,d/d_{\ka})$.
The first two regularity conditions on a perturbation $\vt$ of the section 
$s_1^d$ over $\ov\M_1^0(\Pf,d;J)$ are that for every $\ka\!\in\!\S_0(Y;J)$\\
${}\quad$ ($\vt1a$) the section  
$\d_{\ka,d/d_{\ka}}^{1,0}\!\vt|_{\M_1^0(\ka,d/d_{\ka})}$ 
is transverse to the zero set in $\W_{\ka,d/d_{\ka}}^{1,0}$;\\
${}\quad$ ($\vt1b$) the section $\d_{\ka,d/d_{\ka}}^{1,0}\!\vt$
does not vanish on $\M_1^{\{0\}}(\ka,d/d_{\ka}))\!-\!\M_1^0(\ka,d/d_{\ka}))$.\\
By Lemma~\ref{g1conestr_lmm1}, the $n\!=\!1$ case of Lemma~\ref{g1str_lmm1},
and~\e_ref{coker0decomp_e2},
the collection of multisections $\vt$ of $\V_1^d$ that satisfy ($\vt1a$) and ($\vt1b$)
is open and dense in the space of all multisections of~$\V_1^d$.\\

\noindent
The next lemma, which is the analogue of Lemma~\ref{g1conestr_lmm1} for
the strata $\U_{\T;1}^m(\Pf;J)$ of Lemma~\ref{g1str_lmm2},
is a special case of Proposition~\ref{g1cone-g1conebdstr_prp} and 
Lemma~\ref{g1cone-g1conebdstr_lmm} in~\cite{g1cone}. 
For any $b\!\in\!\U_{\T;1}^m(\Pf;J)$, we~put
$$\F^1{\T}_b=\big\{
\ti\ups\!=\!(\ti{\ups}_i)_{i\in\chi({\T})}\in\F{\T}_b\!:
\sum_{i\in\chi({\T})}\!\!{\cal D}_{J,i}\ti{\ups}_i\!=\!0\big\}.$$

\begin{lmm}
\label{g1conestr_lmm2}
If $d$, $\L$, and $\V_1^d$ are as in Proposition~\ref{cone_prp2},
there exists $\de(d)\!\in\!\R^+$  with the following property.
If $J$ is an almost complex structure on~$\Pf$
such that $\|J\!-\!J_0\|_{C^1}\!<\!\de(d)$,
then the requirements of Lemma~\ref{g1str_lmm2} and of Lemma~\ref{g1conestr_lmm1} 
are satisfied for all appropriate bubble types.
Furthermore, if $\T\!=\!(I,\aleph;\under{d})$
is a bubble type such that $\sum_{i\in I}\!d_i\!=\!d$ and 
$d_i\!=\!0$ for all minimal elements $i$ of $I$, then
the restriction \hbox{$\V_1^d\!\lra\!\U_{\T}(\Pf;J)$}
is a smooth complex vector orbibundle of rank~$5d\!+\!1$.
In addition, for every integer
$$m\in \big(\max(|\chi({\T})|\!-\!4,1),|\chi({\T}|\big),$$
there exist a neighborhood $U_{\T}^m$ of $\U_{\T;1}^m(\Pf;J)$ in 
$\X_1(\Pf,d)$ and a topological vector orbibundle 
$$\V_{1;{\T}}^{d;m} \lra \ov\M_1^0(\Pf,d;J)\cap U_{\T}^m$$
such~that
$\V_{1;{\T}}^{d;m}\!\lra\!\M_1^0(\Pf,d;J)\cap U_{\T}^m$
is a smooth complex vector orbibundle contained in $\V_1^d$ and
$$\V_{1;{\T}}^{d;m}|_{\U_{\T;1}^m(\Pf;J)}
=\big\{\xi\!\in\!\V_1^d|_b\!: b\!\in\!\U_{\T;1}^m(\Pf;J);~~
\D_{\T}(\xi\!\otimes\!\ti{\ups})\!=\!0
~\forall \ti{\ups}\!\in\!\F^1{\T}_b\big\}.$$
There also exists a continuous vector-bundle isomorphism
$$\ti\phi_{\T}^m\!:
\pi_{{\cal F}^1{\T}_{\de}}^*\big(\V_{1;\T}^{d;m}
\big|_{\U_{\T;1}(\Pf;J)\cap U_{\T}^m}\big)
\lra \V_{1;{\T}}^{d;m}\big|_{\ov\M_1^0(\Pf,d;J)\cap U_{\T}^m},$$
covering the homeomorphism $\phi_{\T}$ of Lemma~\ref{g1str_lmm1},
such that $\ti\phi_{\T}^m$ is the identity over $\U_{\T;1}^m(\Pf;J)$. 
Finally, if $\T$ and $\T'$ are two bubble types as above and 
$m,m'\!\in\!\Z^+$, then 
$$\V_{1;{\T}'}^{d;m'}
\big|_{\U_{\T}^m(\Pf;J)\cap U_{\T'}^{m'}}
\subset\V_{1;{\T}}^{d;m}
\big|_{\U_{\T;1}^m(\Pf;J)\cap U_{\T'}^{m'}}
\quad\hbox{if}\quad m'\!\ge\!m.$$\\
\end{lmm}

\noindent
If $[b]\!\in\!\U_{\T;\ka;1}^m$, we put
\begin{alignat*}{1}
\Ga_-(b;\L;0)&=\big\{\xi\!\in\!\Ga_-(b;\L)\!:\ti\D_{\T,i}\xi\!=\!0
~\forall i\!\in\!\chi({\T})\big\};\\
\ti\V_{1;{\T}}^d|_b&=
\big\{[\xi]\!\in\!\V_{1;{\T}}^d|_b\!:\D_{\T,i}\xi\!=\!0
~\forall i\!\in\!\chi({\T})\big\}\subset\V_{1;{\T}}^{d;m}.
\end{alignat*}
In this case,
the standard analogue for $b$ of the long exact sequence \e_ref{les_e} has six terms. 
However, replacing the fourth term by the kernel of the outgoing map at the fourth term,
we~get
\begin{equation}\label{les_e2}
0 \lra \Ga_-(b;TY) \stackrel{i_0}{\lra} \Ga_-(b;T\Pf)
\stackrel{j_0}{\lra} \Ga_-(b;\L) \stackrel{\d_0}{\lra}  H_J^1(\pi_B(b);TY) \lra 0.
\end{equation}
By Theorem~\ref{g1comp-reg_thm} in~\cite{g1comp}, the linear operator
$$\ti\D_{\T,i}^{\Pf}\!:\big\{\ze\!\in\!\Ga_-(b;T\Pf)\!:\ze|_{\Si_{b,\aleph}}\!=\!0\big\}
\lra T_{\ev_P(b)}\Pf, \quad
\ze\lra\na^J_{e_{\i}}\big(\ze|_{\Si_{b,i}}\big),$$
is surjective for every $[b]\!\in\!\U_{\T}(\Pf;J)$,
with ${\T}$ as in Lemma~\ref{g1str_lmm2} and $i\!\in\!\chi({\T})$,
if $J$ is sufficiently close to~$J_0$.
It follows that the homomorphism
$$\ti{j}_0\!:\Ga_-(b;T\Pf)\lra \Ga_-(b;\L)\big/\Ga_-(b;\L;0),$$
induced by the map $j_0$ in~\e_ref{les_e2} is surjective for every
$[b]\!\in\!\U_{\T;\ka}$.
Thus, the family of boundary operators $\d_0$ in~\e_ref{les_e2} with 
$[b]\!\in\!\U_{\T;\ka}$ induces a surjective bundle homomorphism
\begin{equation}\label{surjhom_e}
\d_{\ka,d/d_{\ka}}^{1,q}\!:
\ti\V_{1;\T}^d\lra\pi_B^*\W_{\ka,d/d_{\ka}}^{0,q}
\subset\W_{\ka,d/d_{\ka}}^{1,0}
\end{equation}
over $\U_{\T;\ka;1}$, if $\U_{\T;\ka}\!\subset\!\ov\M_1^q(\ka,d/d_{\ka})$.
Furthermore, 
\begin{equation}\label{pulloper_e}
\d_{\ka,d/d_{\ka}}^{1,q}=\pi_B^*\d_{\ka,d/d_{\ka}}^{0,q},
\end{equation}
where $\d_{\ka,d/d_{\ka}}^{0,q}$ is the surjective bundle homomorphism
over $\U_{\bar\T;\ka}\!\subset\!\ov\M_0^q(\ka,d/d_{\ka})$
defined similarly to~$\d_{\ka,d/d_{\ka}}^{1,q}$.\\

\noindent
We now state additional regularity conditions 
on a perturbation $\vt$ of~$s_1^d$.
We will require that for every $\ka\!\in\!{\cal S}_0(Y;J)$:\\
${}\quad$ ($\vt2a$) the sections
$\d_{\ka,d/d_{\ka}}^{1,1}\!\vt|_{\cM_{1,1}\times\M_{0,1;1}^0(\ka,d/d_{\ka})}$
and 
$\d_{\ka,d/d_{\ka}}^{1,1}\!\vt|_{\partial\ov\cM_{1,1}\times
\M_{0,1;1}^0(\ka,d/d_{\ka})}$
are transverse\\
${}\qquad\qquad$ to the zero set in $\pi_B^*\W_{\ka,d/d_{\ka}}^{1,0}$;\\
${}\quad$ ($\vt2b$) the section $\d_{\ka,d/d_{\ka}}^{1,1}\!\vt$ does not vanish
on $\ov\M_{1;1}^1(\ka,d/d_{\ka})\!-\!\M_{1;1}^{1;0}(\ka,d/d_{\ka})$;\\
${}\quad$ ($\vt2c$) for $q\!\ge\!2$,
the section $\d_{\ka,d/d_{\ka}}^{1,q}\!\vt$ does not vanish
on $\ov\M_1^q(\ka,d/d_{\ka})$.\\
By Lemma~\ref{g1conestr_lmm2}, \e_ref{spacesdim_e1}, and~\e_ref{rankk_e},
the collection of multisections $\vt$ of $\V_1^d$ that satisfy ($\vt2a$) and 
($\vt2c$) with $q\!\ge\!4$
is open and dense in the space of all multisections of~$\V_1^d$.
By \e_ref{spacesdim_e0}, \e_ref{rankk_e}, \e_ref{dim_e3}, and~\e_ref{pulloper_e},
the collection of multisections $\vt$ of $\V_1^d$ that satisfy 
one of the three remaining conditions,
i.e.~($\vt2b$), ($\vt2c$) with $q\!=\!2$, or ($\vt2c$) with $q\!=\!3$, 
is nonempty and open in the space of all multisections of~$\V_1^d$,
but not dense.
Nevertheless, by considering the decompositions of the intersections of 
the corresponding subspaces of $\ov\M_1^0(\ka,d/d_{\ka})$ analogous to~\e_ref{mkdecomp_e}, 
it is straightforward to see that the intersection of these three open sets is still nonempty.
Alternatively, note~that
\begin{alignat*}{1}
\dim_{\C}\big(\ov\M_{1;1}^1(\ka,d/d_{\ka})\!-\!\M_{1;1}^{1;0}(\ka,d/d_{\ka})\big)
&=2(d/d_{\ka})-2=\rk{\cal W}_{\ka,d/d_{\ka}}^{0,1};\\
\dim_{\C}\ov\M_1^3(\ka,d/d_{\ka})
&=2(d/d_{\ka})-2=\rk{\cal W}_{\ka,d/d_{\ka}}^{0,3};\\
\dim_{\C}\ov\M_1^2(\ka,d/d_{\ka})\cap\ov\M_1^3(\ka,d/d_{\ka})
&\le 2(d/d_{\ka})-3.
\end{alignat*}
Furthermore, the space $\ov\M_{1;1}^1(\ka,d/d_{\ka})\!-\!\M_{1;1}^{1;0}(\ka,d/d_{\ka})$
has two irreducible components.
One of them is contained in $\ov\M_1^2(\ka,d/d_{\ka})$, while 
the other intersects $\ov\M_1^2(\ka,d/d_{\ka})$ in subvariety of complex dimension
$2(d/d_{\ka})\!-\!3$.
Thus, if $\vt$ is a generic multisection that satisfies ($\vt2c$) with $q\!=\!2$,
its restrictions~to 
\begin{equation}\label{modspaces_e1}
\ov\M_{1;1}^1(\ka,d/d_{\ka})\!-\!\M_{1;1}^{1;0}(\ka,d/d_{\ka})
\qquad\hbox{and}\qquad  \ov\M_1^3(\ka,d/d_{\ka})
\end{equation}
have finite zero sets, divided equally between positive and negative zeros.
These zeros can be removed in pairs by modifying $\vt$ outside of the boundary
strata of~\e_ref{modspaces_e1}.\\

\noindent
It remains to state one more regularity assumption on $\vt$.
If $\hat{I}\!=\!\chi({\T})\!=\!\{h\}$ is a single-element set,
for every $[b]\!\in\!\U_{\T;\ka;1}$, we~put
\begin{alignat*}{1}
\Ga_-(b;T\Pf;TY)&=\big\{\ze\!\in\!\Ga_-(b;T\Pf)\!:
\ti\D^{\Pf}_{\T,h}\ze\in T_{\ev_P(b)}Y\big\}\\
\Ga_-(b;T\Pf;T\ka)&=\big\{\ze\!\in\!\Ga_-(b;T\Pf)\!:
\ti\D^{\Pf}_{\T,h}\ze\in T_{\ev_P(b)}\ka\big\}
\subset\Ga_-(b;T\Pf;TY).
\end{alignat*}
Since $d\{u_b|_{\Si_{b,h}}\}|_{\i}\!=\!0$ by Lemma~\ref{g1str_lmm2},
the subspaces 
$$\Ga_-(b;T\Pf;TY),\Ga_-(b;T\Pf;T\ka)\subset\Ga_-(b;T\Pf)$$
are in fact independent of the choice of connection $\na^J$ in~$T\Pf$.
Furthermore, by Theorem~\ref{g1comp-reg_thm} in~\cite{g1comp},
\begin{equation}\label{isom_e1}
\Ga_-(b;T\Pf;TY)/\Ga_-(b;T\Pf;T\ka) \approx N_Y{\ka}|_{\ev_P(b)}
\quad\hbox{via}\quad
\ze\lra [\ti\D^{\Pf}_{\T,h}\ze].
\end{equation}
By the paragraph following Lemma~\ref{g1conestr_lmm2}
and condition~$(J_Y2)$ of Definition~\ref{rigid_dfn},
\begin{equation}\label{isom_e2}
\Im j_0|_{\Ga_-(b;T\Pf;TY)}=\ker\d_0\cap\Ga_-(b;\L;0)
\quad\hbox{and}\quad
\ker j_0\subset \Ga_-(b;T\Pf;T\ka),
\end{equation}
where $j_0$ and $\d_0$ are as in~\e_ref{les_e2}.
Let
$$\ti{H}_J^1(\pi_B(b);TY)=
\Ga_-(b;\L;0)\big/\Im j_0|_{\Ga_-(b;T\Pf;T\ka)}.$$
The vector spaces $\ti{H}^1_J(\pi_B(b);TY)$ and the quotient projection maps
induce a vector bundle
over $\M_{1;1}^{1;0}(\ka,d/d_{\ka})$, which we denote by~$Q_{\ka,d/d_{\ka}}^{1,1}$,
and a surjective bundle homomorphism
$$\ti\d_{\ka,d/d_{\ka}}^{1,1}\!:
\ti\V_{\T;1}^d\!=\!\V_{\T;1}^{d;1}
\lra Q_{\ka,d/d_{\ka}}^{1,1}.$$
On the other hand, the boundary operators $\d_0$ in~\e_ref{les_e2} induce
a surjective bundle homomorphism
$$\pi_{\ka}^+\!:Q_{\ka,d/d_{\ka}}^{1,1}\lra\pi_B^*\W_{\ka,d/d_{\ka}}^{0,1}$$
over $\M_{1;1}^{1;0}(\ka,d/d_{\ka})$.
By \e_ref{isom_e1} and \e_ref{isom_e2}, 
$$\ker\pi_{\ka}^+ \approx \pi_B^*(L_0^*\!\otimes\!\ev_0^*N_Y\ka).$$
We also have a surjective bundle homomorphism
$$\pi_{\ka}^-\!:Q_{\ka,d/d_{\ka}}^{1,1}\lra\pi_B^*(L_0^*\!\otimes\!\ev_0^*N_Y\ka).$$
It is induced by the map
\begin{equation}\label{splittinghom_e}
j_0\ze\lra (-2\pi J)\big[\na_{e_{\i}}^J\ze|_{\Si_{b,h}}\big]\in N_Y\ka
\qquad\hbox{if}\quad
\ze\!\in\!\Ga(b;T\Pf),~j_0\ze\!\in\!\Ga_-(b;\L;0).
\end{equation}
Thus, we obtain a splitting of $Q_{\ka,d/d_{\ka}}^{1,1}$:
\begin{equation}\label{bundlesplit_e}
\pi_{\ka}^-\oplus\pi_{\ka}^+\!: Q_{\ka,d/d_{\ka}}^{1,1} \lra
\pi_B^*(L_0^*\!\otimes\!\ev_0^*N_Y\ka) \oplus \pi_B^*\W_{\ka,d/d_{\ka}}^{0,1}
\end{equation}
over $\M_{1;1}^{1;0}(\ka,d/d_{\ka})$.
We note that
\begin{equation}\label{rank_e3}
\rk Q_{\ka,d}^{1,1}=2d
\end{equation}
for all $d\!\in\!\Z^+$.\\

\noindent
Our final regularity condition on $\vt$ is that for every $\ka\!\in\!{\cal S}_0(Y;J)$:\\
${}\quad$ ($\vt3$) the section $\ti\d_{\ka,d/d\ka}^{1,1}\vt$
does not vanish over $\M_{1;1}^{1;0}(\ka,d/d_{\ka})$.\\
By Lemma~\ref{g1conestr_lmm2}, \e_ref{dim_e3}, and~\e_ref{rank_e3},
the collection of multisections $\vt$ of $\V_1^d$ that satisfy ($\vt3$)
is open and dense in the space of all multisections of~$\V_1^d$.
We denote by $\ti{\cal A}_1^d(s;J)$ the collection of all multivalued perturbations
of the section $s_1^d$ of $\V_1^d$ over $\ov\M_1^0(\Pf,d;J)$
that satisfy the regularity conditions ($\vt1$)-($\vt3$).
By the above, $\ti{\cal A}_1^d(s;J)$ is a nonempty open, but not dense, subset 
of the space of all multisections of~$\V_1^d$.\\

\noindent
It is possible to use a dense open collection of perturbations 
in the statement of Proposition~\ref{coneg0contr_prp} below.
However, using such a collection would needlessly complicate
its proof by enlarging the zero set of the sections 
$\d_{\ka,d/d_{\ka}}^{1,0}\!\vt$ by homologically trivial subspaces of
$\ov\M_1^0(\ka,d/d_{\ka})$.
This would also require stating the analogue of~($\vt3$)
for the $q\!=\!2,3$ cases of~($\vt2c$).

\subsection{Proof of Proposition~\ref{cone_prp2}}
\label{cone_subs2}

\noindent
In this subsection, we finally prove Proposition~\ref{cone_prp2}.
It follows immediately from \e_ref{conedecomp_e}, \e_ref{coneg1contr_e}, 
Proposition~\ref{coneg0contr_prp}, and Corollary~\ref{bdcorr_crl}.

\begin{prp}
\label{coneg0contr_prp}
Suppose $J$, $d$, $\L$, and $\V_1^d$ are as in Proposition~\ref{cone_prp2},
$\vt\!\in\!\ti{\cal A}_1^d(s;J)$ is a regular perturbation
of the section $s_1^d$ of $\V_1^d$ on $\ov\M_1^0(\Pf,d;J)$,
$\ka\!\in\!\S_0(Y;J)$, and
${\T}\!=\!(I,\aleph;\under{d})$ is a bubble type such that 
$\sum_{i\in I}\!d_i\!=\!d$. 
If $|I|\!>\!1$ or $\aleph\!\neq\!\eset$, 
for every compact subset~$K$ of $\U_{\T;\ka}$,
there exist $\ep_{\vt}(K)\!\in\!\R^+$ and 
an open neighborhood~$U(K)$ of~$K$ in $\ov\M_1^0(\Pf,d;J)$ such~that
$$\{s_1^d\!+\!t\vt\}^{-1}(0)\!\cap\! U(K)=\eset
\qquad\forall\, t\!\in\!(0,\ep_{\vt}(K)).$$
If $|I|\!=\!1$ and $\aleph\!=\!\eset$, for every compact subset~$K$ of $\U_{\T;\ka}$, 
there exist $\ep_{\vt}(K)\!\in\!\R^+$ and an open neighborhood~$U(K)$ 
of~$K$ in $\ov\M_1^0(\Pf,d;J)$ with the following properties:\\
${}\quad$ (a) the section $s_1^d\!+\!t\vt$ is transverse to the zero set
in $\V_1^d$ over $U(K)$ for all $t\!\in\!(0,\ep_{\vt}(K))$;\\
${}\quad$ (b) for every open subset $U$ of $\ov\M_1^0(\Pf,d;J)$,
there exists $\ep(U)\!\in\!(0,\ep_{\vt}(K))$ such that
\begin{gather*}
^{\pm}\!\big|\{s_1^d\!+\!t\vt\}^{-1}\!\cap\! U\big|=
\blr{e(\W_{\ka,d/d_{\ka}}^{1,0}),\big[\ov\M_1^0(\ka,d/d_{\ka})\big]}
-\cC_{\M_{1;1}^{1;0}(\ka,d/d_{\ka})}\big(\d_{\ka,d/d_{\ka}}^{1,0}\!\vt\big)\\
\hbox{if}\qquad
\big\{\d_{\ka,d/d_{\ka}}^{1,0}\!\vt\big\}^{-1}(0)\!-\!\M_{1;1}^1(\ka,d/d_{\ka})
\!\subset\!K\!\subset\!U\!\subset\!U(K),~t\!\in\!(0,\ep(U)).
\end{gather*}
\end{prp}

\noindent
In other words, the $s_1^d$-contribution from  the main stratum $\M_1^0(\ka,d/d_{\ka})$  
of the space $\ov\M_1^0(\ka,d/d_{\ka})$ to the number 
$$\blr{e(\V_1^d),\big[\ov\M_1^0(\Pf,d;J)\big]},$$
as computed via a perturbation from the open collection $\ti{\cal A}_1^d(s;J)$,
is the euler class of the vector bundle 
$\W_{\ka,d/d_{\ka}}^{1,0}$ over~$\ov\M_1^0(\ka,d/d_{\ka})$
minus the $\d_{\ka,d/d_{\ka}}^{1,0}\!\vt$-contribution
to the latter euler class from the the zeros of $\d_{\ka,d/d_{\ka}}^{1,0}\!\vt$
that lie in $\partial\ov\M_1^0(\ka,d/d_{\ka})$.
Since $\d_{\ka,d/d_{\ka}}^{1,0}\!\vt|_{\M_1^0(\ka,d/d_{\ka})}$
is transverse to the zero set in $\M_1^0(\ka,d/d_{\ka})$,
$$^{\pm}\big|\big\{\d_{\ka,d/d_{\ka}}^{1,0}\!\vt\big\}^{-1}(0)
\!\cap\! \M_1^0(\ka,d/d_{\ka})\big|
=\blr{e(\W_{\ka,d/d_{\ka}}^{1,0}),\big[\ov\M_1^0(\ka,d/d_{\ka})\big]}
-\cC_{\partial\ov\M_1^0(\ka,d/d_{\ka})}\big(\d_{\ka,d/d_{\ka}}^{1,0}\!\vt\big),$$
by Definition~\ref{contr_dfn}.
None of the boundary strata of $\ov\M_1^0(\ka,d/d_{\ka})$
contributes to the euler class of~$\V_1^d$.\\ 

\noindent
{\it Proof:} (1) If ${\T}$ is a bubble type such that $d_i\!\neq\!0$
for some minimal element $i\!\in\!I$, the conclusion of Proposition~\ref{coneg0contr_prp}
follows by the same argument as in the proof of Lemma~\ref{maincontr_lmm3}
and at the end of Subsection~\ref{maincontr_subs}.
The key difference in the $|I|\!=\!1$, $\aleph\!=\!\eset$ case is that
the section $\d_{\ka,d/d_{\ka}}^{1,0}\!\vt$ may vanish on $\partial\ov\M_1^0(\ka,d/d_{\ka})$.
In addition, by the regularity assumptions, ($\vt1b$), ($\vt2b$), and ($\vt2c$), 
\begin{gather*}
\big\{\d_{\ka,d/d_{\ka}}^{1,0}\!\vt\big\}^{-1}(0)
-\M_1^0(\ka,d/d_{\ka})\subset \M_{1;1}^{1;0}(\ka,d/d_{\ka})
\qquad\Lra\\
^{\pm}\big|\big\{\d_{\ka,d/d_{\ka}}^{1,0}\!\vt\big\}^{-1}(0)
\!\cap\! \M_1^0(\ka,d/d_{\ka})\big|
=\blr{e(\W_{\ka,d/d_{\ka}}^{1,0}),\big[\ov\M_1^0(\ka,d/d_{\ka})\big]}
-\cC_{\M_{1;1}^{1;0}(\ka,d/d_{\ka})}\big(\d_{\ka,d/d_{\ka}}^{1,0}\!\vt\big).
\end{gather*}
(2) If ${\T}$ is a bubble type such that $d_i\!=\!0$
for all minimal elements $i\!\in\!I$ and $|\hat{I}|\!=\!1$,
or more generally $|\chi(\T)|\!=\!1$, nearly the same argument still applies.
In this case, ${\cal F}^1\T\!=\!{\cal FT}$ and the conclusions of
Lemmas~\ref{maincontr_lmm3} and~\ref{maincontr_lmm4} are still valid.
The key difference is that the normal bundle of 
$\U_{\T;\ka;1}$ in $\U_{\T;1}(\Pf;J)$ 
is not given by the cokernels of the homomorphisms~$i_0$ in the long exact sequence~\e_ref{les_e}.
Instead, up to the action of the automorphism group of~$b$,
the fiber of ${\cal N}^{\ka}{\T}$ at $[b]\!\in\!\U_{\T;\ka;1}$  is
\begin{gather*}
{\cal N}^{\ka}{\T}=\Ga_-(b;T\Pf;0)\big/\Ga_-(b;TY;0),
\qquad\hbox{where}\\
\begin{aligned}
\Ga_-(b;T\Pf;0)&=\big\{\ze\!\in\!\Ga_-(b;T\Pf)\!: \ti\D^{\Pf}_{\T,h}\ze=0\big\}\\
\Ga_-(b;TY;0)&=\big\{\ze\!\in\!\Ga_-(b;TY)\!: \ti\D^{\Pf}_{\T,h}\ze=0\big\}
\subset \Ga_-(b;T\ka),
\end{aligned}
\end{gather*}
if $h$ is the unique element of $\chi(\T)$.
The reason for this is that 
$$\U_{\T;1}(\Pf;J)
=\big\{b\!\in\!\U_{\T}(\Pf;J)\!:d\{u_b|_{\Si_{b,h}}\}\!=\!0\big\},$$
by Lemma~\ref{g1str_lmm2}.
Since the linear operator 
$$\D_{\T,h}^{\ka}\!:\big\{\ze\!\in\!\Ga_-(b;T\ka)\!:\ze|_{\Si_{b,\aleph}}\!=\!0\big\}
\lra T_{\ev_P(b)}\ka, \quad
\ze\lra\na^{\ka}_{e_{\i}}\big(\ze|_{\Si_{b,h}}\big),$$
is surjective, the image of the homomorphism $j_0$ in the long exact sequence~\e_ref{les_e2}
on $\Ga_-(b;T\Pf;0)$ is the same as on~$\Ga_-(b;T\Pf;T\ka)$.
Thus, the analogue of the bundle $\V_-$ of Subsection~\ref{maincontr_subs}
in this case is the bundle $Q_{\ka,d/d_{\ka}}^{1,1}$ described in the previous
subsection.
The section of $\V_-$ induced by $\vt$ is~$\ti\d_{\ka,d/d_{\ka}}^{1,1}\!\vt$.
Its composition with the map $\pi_{\ka}^+$ in~\e_ref{bundlesplit_e}
is~$\d_{\ka,d/d_{\ka}}^{1,1}\!\vt$.
Thus, Proposition~\ref{coneg0contr_prp} in this case follows from
the regularity assumptions ($\vt2b$) and~($\vt3$),
by the same argument as in Subsection~\ref{maincontr_subs}.\\
(3) Finally, suppose $\T\!=\!(I,\aleph;\under{d})$ is a bubble type such that
$d_i\!=\!0$ for all minimal elements $i\!\in\!I$ and $|\hat{I}|\!\ge\!2$.
In this case, the dimension of the fibers of ${\cal F}^1{\T}$ may not be constant
over $\U_{\T;1}(\Pf;J)$.
Thus, we modify the setup of the second part of Subsection~\ref{maincontr_subs}
by working directly with the normal bundle to the smooth submanifold 
${\cal F}^1\T^{\eset}|_{\U_{\T;\ka;1}}$
in ${\cal F}^1\T^{\eset}|_{\U_{\T;1}(\Pf;J)}$.
Let $\ga\!\lra\!\Bbb{P}\F\T$ be the tautological line bundle
and
$$V=\pi_{\Bbb{P}\F\T}^*\big(\pi_P^*\E^*\!\otimes\!\ev_P^*T\Pf) \lra \Bbb{P}\F\T.$$
We define the section $\al_{\T}$ of $\ga^*\!\otimes\!V$ over $\Bbb{P}\F\T$ by
\begin{gather*}
\big\{\al_{\T}\big(b,(\ti{v}_i)_{i\in\chi(\T)}\big)\big\}(b,\psi)
=\!\sum_{i\in\chi(\T)}\!\!\!\!{\cal D}_{\T,i}\big(b,\psi_{x_{h(i)}}\ti{v}_i\big)
\in T_{\ev_P(b)}\Pf,\\
\hbox{if}\qquad \big(b,(\ti{v}_i)_{i\in\chi(\T)}\big)\!\in\!\ga, 
\quad (b,\psi)\!\in\!\E_{\pi_P(b)}.
\end{gather*}
With our assumptions on $J$, this section is transverse to the zero set
and thus
$$\U_1(\Pf;J)\equiv \al_{\T}^{-1}(0)$$
is a smooth suborbifold of $\Bbb{P}\F{\T}$.
For a similar reason, so is
$$\U_{1;\ka}\equiv \U_1\cap \Bbb{P}\F{\T}\big|_{\U_{\T;\ka}}.$$
Let ${\cal N}^{\ka}{\T}$ denote the normal bundle of $\U_{1;\ka}$ 
in $\U_1(\Pf;J)$.
Up to the action of the automorphism group of $[b,[\ti{\ups}]]\!\in\!\U_{1;\ka}$,
\begin{gather*}
{\cal N}^{\ka}{\T}|_{\big[b,[\ti{\ups}]\big]}=
\Ga_-(b;T\Pf;\ti{\ups})\big/\Ga_-(b;TY;\ti{\ups}),
\qquad\hbox{where}\\
\begin{aligned}
\Ga_-(b;T\Pf;\ti{\ups})
&= \big\{\ze\!\in\!\Ga_-(b;T\Pf)\!:
\sum_{i\in\chi(\T)}\!\!\!\!(\psi_{x_{h(i)}}\ti{v}_i)
\ti\D^{\Pf}_{\T,i}\ze=0~\forall\psi\!\in\!\E_{\pi_P(b)}\big\},\\
\Ga_-(b;TY;\ti{\ups})
&= \big\{\ze\!\in\!\Ga_-(b;TY)\!:
\sum_{i\in\chi(\T)}\!\!\!\!(\psi_{x_{h(i)}}\ti{v}_i)
\ti\D^{\Pf}_{\T,i}\ze=0~\forall\psi\!\in\!\E_{\pi_P(b)}\big\}
\subset \Ga_-(b;T\ka).
\end{aligned}
\end{gather*}
Thus, there is a natural surjective bundle homomorphism
$$\V_-\! \equiv\! \pi_{\Bbb{P}\F\T}^*\V_1^d \big/j_{\ka}({\cal N}^{\ka}\T)
\lra \pi_B^*\W_{\ka,d/d_{\ka}}^{0,q}
\qquad\hbox{if}\quad \U_{\T}(\Pf;J)\subset \ov\M_1^q(\Pf,d;J),$$
where $j_{\ka}\!:{\cal N}^{\ka}{\T}\!\lra\!\pi_{\Bbb{P}\F\T}^*\V_1^d$ 
is the injective bundle homomorphism induced by the maps~$j_0$ in~\e_ref{les_e2}.
We~put
\begin{gather*}
{\cal F}=\pi_{\Bbb{P}\F\T}^*{\cal FT}\lra\Bbb{P}\F\T;\\
{\cal F}^{1\, ,\eset}=\big\{\big(b,[\ti\ups];\ups\big)\!\in\!{\cal F}^{\eset}\!:
\big(b,[\ti\ups]\big)\!\in\!\U_1(\Pf;J);~[\rho(\ups)]\!=\![\ti\ups]\big\}.
\end{gather*}
The smooth orbifold ${\cal F}^{1,\,\eset}$ is diffeomorphic to 
${\cal F}^1{\T}^{\eset}$ by the projection map
$$\big(b,[\ti{\ups}];\ups\big)\lra(b;\ups).$$
Furthermore, ${\cal F}^{1\, ,\eset}\!\lra\!\U_1(\Pf;J)$
is a fiber bundle of smooth varieties.  
We can thus apply the same argument as in the proof of Lemma~\ref{maincontr_lmm3} and 
the end of Subsection~\ref{maincontr_subs}, along with the regularity assumption~($\vt2b$),
to show that 
$$\{s_1^d\!+\!t\vt\}^{-1}(0)\!\cap\! U(K)=\eset
\qquad\forall\, t\!\in\!(0,\ep_{\vt}(K))$$
if $U(K)$ and $t$ are sufficiently small.\\

\noindent
It remains to compute the $\d_{\ka,d/d_{\ka}}^{1,0}\!\vt$-contribution
to the euler class of the bundle $\W_{\ka,d/d_{\ka}}^{1,0}$
over $\ov\M_1^0(\ka,d/d_{\ka})$ from the~set
$$\cZ_{\ka,\vt}\equiv
\big\{\d_{\ka,d/d_{\ka}}^{1,0}\!\vt\big\}^{-1}(0)
-\M_1^0(\ka,d/d_{\ka})
=\big\{\d_{\ka,d/d_{\ka}}^{1,1}\!\vt\big\}^{-1}(0)
\subset \M_{1;1}^{1,0}(\ka,d/d_{\ka}).$$
This contribution is computed by counting the zeros of the section 
$\d_{\ka,d/d_{\ka}}^{1,0}\!\vt\!+\!t\nu$, for a generic section $\nu$
of $\W_{\ka,d/d_{\ka}}^{1,0}$, that lie near~$\cZ_{\ka,\vt}$.
First, let 
$$\pi_{-;\ka}^{\perp}\!:\wt\W_{\ka,d/d_{\ka}}^{1,1}
\lra\pi_P^*\E^*\!\otimes\!\pi_B^*\ev_0^*N_Y\ka$$
denote the (quotient) projection map; see~\e_ref{redcoker_e}.
Our regularity assumptions on~$\nu$ will be that the affine map
\begin{equation}\label{affinemap_e}
\psi_{\vt,\nu}\!:\pi_P^*L_{P,1}\!\otimes\!\pi_B^*L_0
  \lra \pi_P^*\E^*\!\otimes\!\pi_B^*\ev_0^*N_Y\ka, \quad
\psi_{\vt,\nu}(b;\ups)=\pi_{-;\ka}^{\perp}\nu(b)+
\{\ti\d_{\ka,d/d_{\ka}}^{1,1}\vt\}|_b\ups,
\end{equation}
over $\cZ_{\ka,\vt}$ is transverse to the zero set
and all zeros of $\psi_{\vt,\nu}$ lie over 
$$\cZ_{\ka,\vt}^0\equiv
\cZ_{\ka,\vt}\cap\big(\cM_{1,1}\!\times\!\M_{0,1;1}^0(\ka,d/d_{\ka})\big).$$
Since the set $\psi_{\vt,\nu}^{-1}(0)$ is finite, it follows that
it lies over a compact subset $K_{\vt,\nu}$ of~${\cal Z}_{\ka,\vt}^0$.
By the regularity assumption ($\vt2a$), these conditions are satisfied 
by sections $\nu$ in a dense open subset of the space of all sections 
of $\W_{\ka,d/d_{\ka}}^{1,0}$.
We~put
$$\partial\cZ_{\ka,\vt} =
\cZ_{\ka,\vt} \cap\big(\partial\ov\cM_{1,1}\!\times\!\M_{0,1;1}^0(\ka,d/d_{\ka})\big).$$

\begin{lmm}
\label{bdcontr_lmm}
Suppose $J$, $d$, $\L$, $\V_1^d$, $\vt\!\in\!\ti{\cal A}_1^d(s;J)$,
and $\ka\!\in\!\S_0(Y;J)$ are as in Proposition~\ref{coneg0contr_prp}
and in Lemma~\ref{g1conestr_lmm2}.
If $\T\!=\!(I,\aleph;\under{d})$ is a bubble type such that
$d_i\!=\!0$ for all minimal elements $i$ of $I$ and
$\hat{I}\!=\!\{h\}$ is a single-element set, then
there exist $\de\in\!C(\U_{\T;\ka;1};\Bbb{R}^+)$,
$U_{\T}^1$, and $\phi_{\T}^1$ as in the $n\!=\!1$ case
of Lemma~\ref{g1str_lmm2}, $\ve\!\in\!C({\cal FT};\R)$,
and a vector bundle isomorphism
$$\Phi_{\T}\!:\pi_{{\cal FT}_{\de}}^*
\big(\wt\W_{\ka,d/d_{\ka}}^{1,1}|_{\U_{\T;\ka;1}}\big)
\lra \W_{\ka,d/d_{\ka}}^{0,1}\big|_{\ov\M_1^0(\ka,d/d_{\ka})\cap U_{\T}^1},$$
covering the homeomorphism $\phi_{\T}^1$ and restricting to the identity over 
$\U_{\T;\ka;1}$ such~that
$$\big|\pi_{-;\ka}^{\perp}\Phi_{\T}^{-1}
\big(\{\d_{\ka,d/d_{\ka}}^{1,0}\!\vt\}(\phi_{\T}^1(\ups)\big)
-\big\{\pi_{\ka}^-\ti{\d}_{\ka,d/d_{\ka}}^{1,1}\!\vt|_b\big\}\rho(\ups)\big|
\le \ve(\ups)\big|\rho(\ups)\big|
\qquad\forall\, \ups\!=\!(b,v)\!\in\!{\cal FT}_{\de}^{\eset},$$
and $\lim_{|\ups|\lra0}\!\ve(\ups)=\!0$.
\end{lmm}

\noindent
In this case, $\aleph\!=\!\eset$ or $\aleph$ contains one element, and
$$\U_{\T;\ka;1}
=\begin{cases}
\cM_{1,1}\!\times\!\cM_{0,1;1}^0(\ka,d/d_{\ka}),&
\hbox{if}~\aleph\!=\!\eset;\\
\partial\ov\cM_{1,1}\!\times\!\M_{0,1;1}^0(\ka,d/d_{\ka}),&
\hbox{otherwise}.
\end{cases}$$
In either case, by Lemma~\ref{g1str_lmm2}, the normal bundle ${\cal F}^1{\T}$
of $\U_{\T;\ka}$ in $\ov\M_1^0(\ka,d/d_{\ka})$ is~${\cal FT}$.
If $\aleph\!=\!\eset$,
$${\cal FT}=\pi_P^*L_{P,1}\!\otimes\!\pi_B^*L_0
\qquad\hbox{and}\qquad  \rho(\ups)=\ups.$$
Otherwise, ${\cal FT}$ is the direct sum of $\pi_P^*L_{P,1}\!\otimes\!\pi_B^*L_0$
with the line of smoothings of the node of $\Si_{b,\aleph}$,
which in this case is a sphere with two points identified.
If $\ups\!\in\!{\cal FT}$, $\rho(\ups)$ is 
the $\pi_P^*L_{P,1}\!\otimes\!\pi_B^*L_0$-component of~$\ups$.\\

\noindent 
Lemma~\ref{bdcontr_lmm} follows fairly easily from constructions in~\cite{g1comp} and~\cite{g1cone}.
However, its proof is notationally involved, and we postpone it until the next subsection.

\begin{crl}
\label{bdcorr_crl}
Suppose $d$, $\L$, $\V_1^d$, $J$, $\vt\!\in\!\ti{\cal A}_1^d(s;J)$,
and $\ka\!\in\!\S_0(Y;J)$ are as in Proposition~\ref{coneg0contr_prp}.
If $\nu$ is a generic perturbation of the section 
$\d_{\ka,d/d_{\ka}}^{1,0}\!\vt$ of $\W_{\ka,d/d_{\ka}}^{1,0}$
over $\ov\M_1^0(\ka,d/d_{\ka})$,  there exist $\ep_{\nu}\!\in\!\R^+$ and 
an open neighborhood $U$ of $\partial\cZ_{\ka,\vt}$
in $\ov\M_1^0(\ka,d/d_{\ka})$ such~that
$$\big\{\d_{\ka,d/d_{\ka}}^{1,0}\!\vt\!+\!t\nu\big\}^{-1}(0)\!\cap\! U=\eset
\qquad\forall\, t\!\in\!(0,\ep_{\nu}).$$
Furthermore, for every compact subset~$K$ of $\cZ_{\ka,\vt}^0$,
there exist $\ep_{\nu}(K)\!\in\!\R^+$ and an open neighborhood~$U(K)$ 
of~$K$ in $\ov\M_1^0(\ka,d/d_{\ka})$ with the following properties:\\
${}\quad$ (a) the section $\d_{\ka,d/d_{\ka}}^{1,0}\!\vt\!+\!t\nu$ is transverse 
to the zero set in $\W_{\ka,d/d_{\ka}}^{1,0}$ over 
$U(K)$ for all $t\!\in\!(0,\ep_{\nu}(K))$;\\
${}\quad$ (b) for every open subset $U$ of $\ov\M_1^0(\ka,d/d_{\ka})$,
there exists $\ep(U)\!\in\!(0,\ep_{\nu}(K))$ such that
\begin{gather*}
^{\pm}\!\big|\{\d_{\ka,d/d_{\ka}}^{1,0}\!\vt\!+\!t\nu\}^{-1}\!\cap\! U\big|=
-\frac{d/d_{\ka}-1}{12} \,
\blr{e(\W_{\ka,d/d_{\ka}}^0),\big[\ov\M_0(\ka,d/d_{\ka})\big]}\\
\hbox{if}\qquad
K_{\vt,\nu}\!\subset\!K\!\subset\!U\!\subset\!U(K),
~t\!\in\!(0,\ep(U)).
\end{gather*}\\
\end{crl}

\noindent
{\it Proof:} (1) Let $\T\!=\!(I,\aleph;\under{d})$ be a bubble type such that
$\sum_{i\in I}\!d_i\!=\!d$ and $d_i\!=\!0$ for all minimal elements $i$ of~$I$.
By the regularity assumption ($\vt2b$), if
$$\cZ_{\vt,{\T}}\equiv\{\d_{\ka,d/d_{\ka}}^{1,0}\!\vt\}^{-1}(0)
\!\cap\!\U_{\T;\ka;1}\neq\eset,$$
then $\hat{I}\!=\!\{h\}$ is a single-element set, while $|\aleph|\!\in\!0,1$.\\
(2)  We denote by ${\cal N}^{\vt}{\T}$ the normal bundle of 
$\cZ_{\vt,{\T}}$ in~$\U_{\T;\ka;1}$.
Similarly to the proof of Lemma~\ref{maincontr_lmm3},
using the homeomorphism $\phi_{\T}^1$ of Lemma~\ref{g1str_lmm2}
and the bundle isomorphism $\Phi_{\T}$ of Lemma~\ref{bdcontr_lmm}, 
we can obtain an identification of neighborhoods of ${\cal Z}_{\vt,{\T}}$ in 
${\cal N}^{\vt}{\T}\!\oplus\!{\cal FT}$ and in $\ov\M_1^0(\ka,d/d_{\ka})$,
$$\phi_{\T;\vt}\!: 
{\cal N}^{\vt}{\T}_{\de}\!\times_{\cZ_{\vt,\T}}\!{\cal FT}_{\de}
\lra \ov\M_1^0(\ka,d/d_{\ka}),$$
and a lift of $\phi_{\T;\vt}$ to a bundle isomorphism,
$$\Phi_{\T;\vt}\!: 
\pi_{{\cal N}^{\vt}{\T}\!\times_{{\cal Z}_{\vt,{\T}}}{\cal FT}_{\de}}^*
\big(\wt{\cal W}_{\ka,d/d_{\ka}}^{1,1}|_{{\cal Z}_{\vt,{\T}}}\big)
\lra {\cal W}_{\ka,d/d_{\ka}}^{0,1}.$$
For $(b;X,\ups)\!\in\!
       {\cal N}^{\vt}{\T}\!\times_{{\cal Z}_{\vt,{\T}}}\!{\cal FT}_{\de}$,
we~put
\begin{alignat*}{1}
s(b;X,\ups)&=\Phi_{\T;\vt}^{-1}
  \big(\{\d_{\ka,d/d_{\ka}}^{1,0}\!\vt\}\phi_{\T}^1(b;X,\ups)\big)
\in \wt\W_{\ka,d/d_{\ka}}^{1,1}\big|_b;\\
\ti{\nu}(b;X,\ups)&=\Phi_{\T;\vt}^{-1}\big(\nu(\phi_{\T}^1(b;X,\ups))\big)
\in \wt\W_{\ka,d/d_{\ka}}^{1,1}\big|_b.
\end{alignat*}
We define $N_s(X)$ and $N_s'(X,\ups)$ in $\wt{\cal W}_{\ka,d/d_{\ka}}^{1,1}|_b$
by
\begin{alignat}{1}
\label{bdcontr_e3a}
s(b;X,0)&=s(b;0,0)+j_bX+N_sX =j|_bX+N_sX; \\
\label{bdcontr_e3b}
s(b;X,\ups)&=s(b;X,0)+N_s'(X,\ups),
\end{alignat}
where $j_b\!:{\cal N}^{\vt}{\T}\!\lra\!\pi_B^*\W_{\ka,d/d_{\ka}}^{0,1}|_b
\subset\wt\W_{\ka,d/d_{\ka}}^{1,1}|_b$
is the derivative of~$s$.
Thus,
\begin{equation}\label{bdcontr_e4a}
N_s(0)\!=\!0,\quad
\big|N_s(X)-N_s(X')\big|<C_K\big(|X|\!+\!|X'|\big)|X\!-\!X'|
\quad\forall\, X,X'\!\in\!{\cal N}^{\vt}{\T}_{\de}.
\end{equation}
By the continuity $\vt$, 
\begin{equation}\label{bdcontr_e4b}
\big|N_s'(X,\ups)\big|\le\ve(\ups),\quad
\big|N_s'(X,\ups)\!-\!N_s'(X',\ups)\big|\le\ve(\ups)\big|X\!-\!X'|
\quad\forall\, X,X'\!\in\!{\cal N}^{\vt}{\T}_{\de},\,
\ups\!\in\!{\cal FT}_{\de},
\end{equation}
for some $\ve\!\in\!C({\cal FT};\Bbb{R}^+)$
such that $\lim_{|\ups|\lra0}\!\ve(\ups)\!=\!0$.
We also have
\begin{equation}\label{bdcontr_e4c}
\begin{split}
\big|\ti\nu(b;X,\ups)\big|&\le C(b)\\
\big|\ti\nu(b;X,\ups)\!-\!\ti{\nu}(b;X',\ups)\big|&
\le C(b)|X\!-\!X'|
\end{split} \qquad
\forall\, b\!\in\!{\cal Z}_{\vt,{\T}},
\, X,X'\!\in\!{\cal N}^{\vt}\T_{\de}, \,\ups\!\in\!{\cal FT}_{\de},
\end{equation}
for some $C\!\in\!C({\cal Z}_{\vt,{\T}};\Bbb{R})$.\\
(3) Let $\pi_{-;B}\!:\wt\W_{\ka,d/d_{\ka}}^{1,1}\!\lra\!\pi_B^*\W_{\ka,d/d_{\ka}}^{0,1}$
be the natural projection map; see Subsection~\ref{gluing1_subs}.
Let $K$ be a precompact open subset of~$\cZ_{\vt,\T}$.
Since the homomorphism
$$j_b\!:{\cal N}^{\vt}{\T}\lra\pi_B^*\W_{\ka,d/d_{\ka}}^{0,1}$$
is an isomorphism by the regularity assumption~($\nu2a$),
by \e_ref{bdcontr_e3a}-\e_ref{bdcontr_e4c} and the Contraction Principle, the equation
$$\pi_{-;B}\big(s(b;X,\ups)+t\ti{\nu}(b;X,\nu)\big)=0$$
has a unique small solution $X\!=\!X_t(\ups)\!\in\!{\cal N}^{\vt}{\T}_b$ for all
$t\!\in\![0,\de_K)$, $\ups\!\in\!{\cal FT}_{\de_K}|_b$, and $b\!\in\!K$.
Furthermore, 
\begin{equation}\label{bdcontr_e5}
\big|X_t(\ups)\big|\le C_K\big(t\!+\!\ve(\ups)\big).
\end{equation}
(4) By the above, the number of zeros of $\d_{\ka,d/d_{\ka}}^{1,0}\!\vt\!+\!t\nu$,
for $t\!\in\!(0,\de_K)$, in a small neighborhood $U_K$ of~$K$ 
in $\ov\M_1^0(\ka,d/d_{\ka})$ is the number of solutions of the equation
\begin{equation*}\begin{split}
\Psi_t(b;\ups) &\equiv 
t^{-1}\pi_{-;\ka}^{\perp}\big(s(b;X_t(\ups);\ups)\!+\!t\ti\nu(b;X_t(\ups);\ups)\big)\\
&=t^{-1}\pi_{-;\ka}^{\perp}N_s'\big(\ups;X_t(\ups)\big)
+\pi_{-;\ka}^{\perp}\ti\nu\big(b;X_t(\ups);\ups\big),
\end{split}\end{equation*}
since $\d_{\ka,d/d_{\ka}}^{1,0}\!\vt|_{\U_{\T;\ka;1}}$ is
a section of $\pi_B^*\W_{\ka,d/d_{\ka}}^{0,1}$ and thus
$\pi_{-;\ka}^{\perp}N_sX\!=\!0$ for all $X\!\in\!{\cal N}^{\vt}{\T}_{\de}$.
By the estimate of Lemma~\ref{bdcontr_lmm}, \e_ref{bdcontr_e5}, and
the smoothness of~$\nu$
\begin{equation}\label{bdcontr_e9}
\big|\Psi_t(b;\ups)-
\big(\pi_{-;\ka}^{\perp}\nu\big(b)+
t^{-1}\{\pi_{\ka}^-\ti\d_{\ka,d/d_{\ka}}^{1,1}\!\vt\}|_b\rho(\ups)
\big)\big|\le C_K\big(t\!+\!\ve(\ups)\big)|\rho(\ups)|.
\end{equation}
By the regularity assumption ($\nu3$), the section 
$\pi_{\ka}^-\ti\d_{\ka,d/d_{\ka}}^{1,1}\!\vt$ does not vanish over $\cZ_{\vt,\T}$.
Suppose $\T$ is a bubble type as in the $|\aleph|\!=\!1$ case  in~(1) above.
By our assumptions on~$\nu$, the affine~map
\begin{equation}\label{bdcontr_e11}
{\cal FT}\lra\pi_P^*\E^*\!\otimes\!\pi_B^*\ev_0^*N_Y\ka,
\qquad \ups\lra \pi_{-;\ka}^{\perp}\nu\big(b)+
\{\pi_{\ka}^-\ti\d_{\ka,d/d_{\ka}}^{1,1}\!\vt\}|_b\rho(\ups),
\end{equation}
which factors through $\psi_{\vt,\nu}$, does not vanish over 
the compact set~$\partial{\cal Z}_{\ka,\vt}$.
Thus, $\Psi_t(b;\ups)\!\neq\!0$ for all $t\!\in\!\Bbb{R}^+$ and 
$\ups\!\in\!{\cal FT}|_{\partial{\cal Z}_{\ka,\vt}}$ sufficiently small.
This concludes the proof of the first statement of Corollary~\ref{bdcorr_crl}.\\
(5) Finally, suppose $\T$ is a bubble type as in the $|\aleph|\!=\!0$ case 
in~(1). If $K$ is a compact subset of $\cZ_{\ka,\vt}^0$ containing
$K_{\vt,\nu}$, by \e_ref{bdcontr_e9}, the regularity assumption~($\nu3$),
and the same cobordism argument as in Subsection~3.1 of~\cite{g2n2and3},
the number of solutions of $\Psi_t(b;\ups)\!=\!0$ with $t\!\in\!\R^+$
and $\ups\!\in\!{\cal FT}|_K$ sufficiently small is 
the number of zeros of the affine map in~\e_ref{bdcontr_e11},
i.e.~$^{\pm}|\psi_{\vt,\nu}^{-1}(0)|$.
Since $\pi_{\ka}^-\ti\d_{\ka,d/d_{\ka}}^{1,1}\!\vt$ does not vanish over~$\cZ_{\ka,\vt}$,
\begin{equation*}\begin{split}
^{\pm}\big|\psi_{\vt,\nu}^{-1}(0)\big|
&=\blr{e(\pi_P^*\E^*\!\otimes\!\pi_B^*\ev_0^*N_Y\ka)
e(\pi_P^*L_{P,1}\!\otimes\!\pi_B^*L_0)^{-1},[\cZ_{\ka,\vt}]}\\
&=\blr{\pi_P^*(-2\la\!+\!\psi_{P,1})+\pi_B^*(c_1(\ev_0^*N_Y\ka)\!+\!\psi_0),[\cZ_{\ka,\vt}]},
\end{split}\end{equation*}
where $\la$ and $\psi_{P,1}$ are the usual tautological classes on $\ov\cM_{1,1}$.
The space $\cZ_{\ka,\vt}$ is the zero set of the section $\d_{\ka,d/d_{\ka}}^{1,1}\!\vt$ 
of the bundle $\pi_B^*\W_{\ka,d/d_{\ka}}^{0,1}$ over $\ov\M_{1;1}^1(\ka,d/d_{\ka})$.
Since this section is transverse to the zero set by~($\nu2a$),
\begin{equation*}\begin{split}
^{\pm}\big|\psi_{\vt,\nu}^{-1}(0)\big|
&=\blr{\pi_P^*(-2\la\!+\!\psi_{P,1})\cdot\pi_B^*e(\W_{\ka,d/d_{\ka}}^{0,1}),
\big[\ov\cM_{1,1}\!\times\!\ov\M_{0,1;1}(\ka,d/d_{\ka})\big]}\\
&= -\frac{1}{24} \, 
\blr{e(\W_{\ka,d/d_{\ka}}^{0,1}),\big[\ov\M_{0,1;1}(\ka,d/d_{\ka})\big]}.
\end{split}\end{equation*}
We note that a generic fiber of the forgetful map 
$$\ti\pi\!: \ov\M_{0,1;1}(\ka,d/d_{\ka}) \lra \ov\M_0(\ka,d/d_{\ka})$$
consists of $2(d/d_{\ka})\!-\!2$ points, corresponding to the branch points
a degree-$d/d_{\ka}$ cover \hbox{$\Bbb{P}^1\!\lra\!\ka$}.
We conclude that for every compact subset $K$ of ${\cal Z}_{\ka,\vt}$ containing~$K_{\vt,\nu}$ 
\begin{equation*}\begin{split}
^{\pm}\!\big|\{\d_{\ka,d/d_{\ka}}^{1,0}\!\vt\!+\!t\nu\}^{-1}\!\cap\! U\big|
&=\, ^{\pm}\big|\psi_{\vt,\nu}^{-1}(0)\big|
= -\frac{1}{24} \, 
\blr{\ti{\pi}^*e({\cal W}_{\ka,d/d_{\ka}}^0),\big[\ov\M_{0,1;1}(\ka,d/d_{\ka})\big]}\\
&=-\frac{1}{24}\big(2(d/d_{\ka})\!-\!2\big) \,
\blr{e({\cal W}_{\ka,d/d_{\ka}}^0),\big[\ov\M_0(\ka,d/d_{\ka})\big]},
\end{split}\end{equation*}
provided that $U$ is a sufficiently small neighborhood
of $K$ in $\ov\M_1^0(\ka,d/d_{\ka})$ and $t\!\in\!(0,\de(U))$.

\subsection{A Genus-One Gluing Procedure}
\label{gluing1_subs}

\noindent
In this subsection, we prove Lemma~\ref{bdcontr_lmm}.
We review the genus-one gluing procedure of 
Subsection~\ref{g1comp-reg1_subs2} in~\cite{g1comp} and
its extensions to the spaces $\Ga_-(b;T\Pf)$ and $\Ga_-(b;\L)$.
As a result, we will be able to describe the behavior of the boundary operator $\d_0$
in the long exact sequence~\e_ref{les_e} for 
$[\ti{b}(\ups)]\!\in\!\M_1^0(\ka,d/d_{\ka})$
with $\ups\!\in\!\ti{\cal F}{\T}^{\eset}$ sufficiently small.\\

\noindent
Let ${\T}\!=\!(I,\aleph;\under{d})$ be a bubble type as in the statement
of Lemma~\ref{bdcontr_lmm}.
If $\ups\!=\!(b,v)\!\in\!\ti{\cal F}{\T}_{\de}^{\eset}$ is small gluing
parameter, let
$$b(\ups)=(\Si_{\ups},j_{\ups};u_{\ups}),
\quad\hbox{where}\quad u_{\ups}=u_b\circ\ti{q}_{\ups},$$
be the (second-stage) approximately $J$-holomorphic map.
Here
$$\ti{q}_{\ups}\!=\!\ti{q}_{\ups_0;2}\!: \Si_{\ups}\lra\Si_b$$
is the basic gluing map constructed in Subsection~\ref{g1comp-reg1_subs2} of~\cite{g1comp}.
In the present case, there is no first stage in this usually two-stage gluing construction,
as there is only one level of bubbles (in fact, only one bubble) to attach.
The key advantage of this gluing construction is that the map $\ti{q}_{\ups}$
is closer to being holomorphic than in the gluing construction used in  Section~\ref{dbar_sec}.
In particular,
$$\big\|\bar{\partial}_Ju_{\ups}\big\|_{\ups,p}\le C(b)\big|\rho(\ups)\big|.$$
If $b\!\in\!\U_{\T;\ka;1}$, then $du_{b,h}|_{\i}\!=\!0$ 
and this estimate improves~to
\begin{equation}\label{dbarest_e}
\big\|\bar{\partial}_Ju_{\ups}\big\|_{\ups,p}\le C(b)\big|\rho(\ups)\big|^2.
\end{equation}
This is immediate from the definition of the map $\ti{q}_{\ups}$.\\

\noindent
We extend the metric $g_{Y,b}$ described in Subsection~\ref{bdcontr_subs}
to a metric $g_{\Pf,b}$ on the bundle~$T\Pf$. Let $\na^J$ be the $J$-compatible connection
corresponding to the Levi-Civita of the metric~$g_{\Pf,b}$.\\

\noindent
Similarly to Subsection~\ref{bdcontr_subs}, let 
$$\Ga_-^{0,1}(b;T\Pf)={\cal H}_{b,P}^{0,1}\otimes T_{\ev_P(b)}\Pf$$
be the space of $u_b^*T\Pf$-valued harmonic $(0,1)$-forms on~$\Si_b$.
If $\ups\!=\!(b,v)$ and $b\!\in\!\U_{\T;\ka;1}$ are as above, we~put
$$\Ga_-^{0,1}(\ups;T\Pf)=\big\{ R_{\ups}\eta\!: \eta\!\in\!\Ga_-^{0,1}(b;T\Pf)\big\}
 \subset\Ga^{0,1}(\ups;T\Pf),$$
where $R_{\ups}\eta$ is a smooth extension of $\eta$ such that 
$R_{\ups}\eta$ is nearly harmonic on the neck attaching the only bubble $\Si_{b,h}$
of $\Si_b$ and below a small collar of the neck and vanishes past a slighter larger collar;
see Subsection~\ref{g1comp-reg1_subs2} in~\cite{g1comp}.
Let
$$\pi_{\ups;-}^{0,1}\!: \Ga^{0,1}(\ups;T\Pf)\lra\Ga_-^{0,1}(b;T\Pf)$$
be the $L^2$-projection map.
We denote its kernel by $\Ga_+^{0,1}(\ups;T\Pf)$. 
By the same argument as in Subsection~2.3 in~\cite{g2n2and3}, we have a decomposition
\begin{gather*}
\Ga(\ups;T\Pf)=\Ga_-(\ups;T\Pf) \oplus \ti\Ga_+(\ups;T\Pf),\\
\hbox{where}\qquad 
\Ga_-(\ups;T\Pf)=\big\{R_{\ups}\ze\!\equiv\!\xi\circ\ti{q}_{\ups}\!:
\ze\!\in\!\Ga_-(b;T\Pf)\big\},
\end{gather*}
such that 
$$D_{J,\ups}\!:\ti\Ga_+(\ups;T\Pf)\lra\Ga_+^{0,1}(\ups;T\Pf)$$
is an isomorphism with fiber-uniformly bounded inverse and
$$\llrr{D_{J,\ups}\ze,\eta}_{\ups,2}=0 \quad\forall\,
\ze\!\in\!\ti\Ga_+(\ups;T\ka)\!\equiv\!\ti\Ga_+(\ups;T\Pf)\cap\Ga(\ups;T\ka),
\, \eta\!\in\!\Ga_-^{0,1}(\ups;T\Pf).$$
Analogously to~\e_ref{dbarest_e}, we also~have
\begin{equation}\label{kerprojest_e1}
\big\|D_{J,\ups}\ze\big\|_{\ups,p}\le C(b)\big|\rho(\ups)\big| \|\ze\|_{\ups,p,1}
\qquad\forall\, \ze\!\in\!\Ga_-(\ups;T\Pf).
\end{equation}
Furthermore, 
\begin{equation}\label{kerprojest_e2}
\big|\pi_{\ups;-}^{0,1}D_{J,\ups}R_{\ups}\ze
+2\pi\rho(\ups)J\,R_{\ups}\ti\D_{\T,h}^{\Pf}\ze\big|
\le C(b)\big|\rho(\ups)\big|^2\|\ze\|_{b,p,1}
\quad\forall\, \ze\!\in\!\Ga_-(b;T\Pf);
\end{equation}
see~(5) of Lemma~\ref{g1comp-reg1_lmm3} in~\cite{g1comp}.
Due to the assumption that $du_{b,h}|_{\i}\!=\!0$, we do not need to require that
$\ze|_{\Si_{b,\aleph}}\!=\!0$.
We also get a slightly sharper bound, though this is not essential.
The estimate~\e_ref{kerprojest_e2} is the fundamental fact behind
the estimate of Lemma~\ref{bdcontr_lmm}.\\

\noindent
Similarly to Subsection~\ref{bdcontr_subs},
the restriction of the homeomorphism $\phi_{\T}^1$ of Lemma~\ref{g1str_lmm2}
can be taken to be of the~form
\begin{gather}
\phi_{\T}^1([\ups])=\big([\ti{b}(\ups)]\big),
\qquad\hbox{where}\quad
\ti{b}(\ups)=(\Si_{\ups},j_{\ups};\ti{u}_{\ups}),
\quad \ti{u}_{\ups}=\exp_{b,u_{\ups}}\!\ze_{\ups}\notag\\
\label{pertest_e}
\ze_{\ups}\!\in\!\ti\Ga_+(\ups;T\ka),\qquad 
\|\ze_{\ups}\|_{\ups,p,1}\le C(b)|\rho(\ups)|^2.
\end{gather}
The last estimate follows from \e_ref{dbarest_e} by the usual argument.\\

\noindent
We denote~by
\begin{gather*}
\Pi_{\ups}^J\!:
\Ga(\ups;T\Pf)\lra\ti\Ga(\ups;T\Pf)\!\equiv\!L^p_1(\ti{b}(\ups);T\Pf)
\qquad\hbox{and}\\
\Pi_{\ups}\!: \Ga(\ups;\L)\!\equiv\!L^p_1(b(\ups);\L)
\lra\ti\Ga(\ups;\L)\!\equiv\!L^p_1(\ti{b}(\ups);\L)
\end{gather*}
the $\na^J$-parallel transport in $T\Pf$ and the $\na$-parallel transport
in $\L$ along the geodesics $\ga_{\ze_{\ups}}$ in $\ka$ of the metric~$g_{\ka,b}$. 
By \e_ref{kerprojest_e1}-\e_ref{pertest_e} and the same argument as 
in Subsection~\ref{g1cone-gluing_subs3} of~\cite{g1cone}, 
there exists an isomorphism
\begin{gather}
\ti{R}_{\ups}\!: \Ga_-(b;T\Pf;0) \lra 
\ti\Ga_-(\ups;T\Pf)\!\equiv\!\ker D_{J,\ti{b}(\ups)}
\qquad\hbox{s.t.} \notag\\
\label{corr_e1}
\big\|\ti{R}_{\ups}\ze-\Pi_{\ups}^JR_{\ups}\ze\big\|_{\ups,p,1}
\le C(b)\big|\rho(\ups)\big|^2\|\ze\|_{b,p,1}
\qquad\forall\, \ze\!\in\!\Ga_-(b;T\Pf;0).
\end{gather}
Similarly, there exists an isomorphism
\begin{gather}
\ti{R}_{\ups}\!: \Ga_-(b;\L;0) \lra 
\ti\Ga_-(\ups;\L)\!\equiv\!\ker\bpar_{\na,\ti{b}(\ups)}
\qquad\hbox{s.t.} \notag\\
\label{corr_e2}
\big\|\ti{R}_{\ups}\xi-\Pi_{\ups}R_{\ups}\xi\big\|_{\ups,p,1}
\le C(b)\big|\rho(\ups)\big|^2\|\xi\|_{b,p,1}
\qquad\forall\, \xi\!\in\!\Ga_-(b;\L;0),
\end{gather}
where again $R_{\ups}\xi\!=\!\xi\circ\ti{q}_{\ups}$.\\

\noindent
We next  describe a convenient family of finite-dimensional spaces 
$$\Ga(b;T\Pf;\L)\subset\Ga(b;T\Pf),$$
parameterized by $b\!\in\!\U_{\T;\ka;1}^{(0)}$,
such that the~homomorphism
$$j_0\!:\Ga(b;T\Pf;\L)\lra\Ga_-(b;\L;0)$$
is an isomorphism.
For every $b\!\in\!\U_{\T;\ka;1}^{(0)}$, let 
$$\Ga_{-,+}(b;T\Pf;N_Y\ka)\approx \ev_P^*N_Y\ka
\qquad\hbox{and}\qquad \Ga_{-,+}(b;T\Pf;\ka)$$ 
be the $L^2$-orthogonal complements of 
$\Ga_-(b;T\Pf;T\ka)$ in $\Ga_-(b;T\Pf;TY)$
and of $\Ga_-(b;T\ka;0)$ in $\Ga_-(b;T\Pf;0)$, respectively.
The map in~\e_ref{splittinghom_e} induces a surjective homomorphism
$$\ti\pi_{\ka}^-\!: \Ga_-(b;\L;0)\lra \ev_P^*N_Y\ka,$$
which restricts to an isomorphism on $j_0(\Ga_{-,+}(b;T\Pf;N_Y\ka))$
and vanishes on $j_0(\Ga_{-,+}(b;T\Pf;\ka))$, where $j_0$ is as in \e_ref{les_e2}.
Let $\Ga_{-,+}(b;\L;0)$ be the $L^2$-orthogonal complement
of $j_0(\Ga_{-,+}(b;T\Pf;\ka))$ in $\ker\ti\pi_{\ka}^-$.
We~set
$$\Ga_{+,-}(b;T\Pf;N_{\Pf}Y)
=\big\{\ze\!\in\!\Ga(b;N_{\Pf}Y)\!:\pi_Y^{\perp}\!\circ\ze\!\in\!\Ga_{-,+}(b;\L;0)\big\},$$
where the normal bundle $N_{\Pf}Y$ of $Y$ in $\Pf$ is identified
with the $g_{\Pf,b}$-orthogonal complement of $TY$ in~$T\Pf$
and
$$\pi_Y^{\perp}\!:T\Pf\lra \L\approx N_{\Pf}Y$$
is the quotient projection.
Then, the map
$$\pi_Y^{\perp}\!: 
\Ga(b;T\Pf;\L)\!\equiv\! \Ga_{-,+}(b;T\Pf;\ka)
\!\oplus\!\Ga_{-,+}(b;T\Pf;N_Y\ka) \!\oplus\!\Ga_{+,-}(b;T\Pf;N_{\Pf}Y)
\lra \Ga_-(b;\L;0)$$
is an isomorphism.
Furthermore, since 
$$\Ga_{-,+}(b;T\Pf;\ka)\oplus\Ga_{-,+}(b;T\Pf;N_Y\ka)\subset\Ga_-(b;T\Pf),$$
the map 
$$\d_b\!:\Ga_{+,-}(b;T\Pf;N_{\Pf}Y) \lra H_J^1(\pi_B(b);N_Y\ka),
\qquad \ze\lra[\pi_{\ka}^{\perp}D_{J,b}\ze],$$
is also an isomorphism, by the definition of the boundary operator $\d_0$ in~\e_ref{les_e2}.\\

\noindent
We now use the subspace $\Ga(b;T\Pf;\L)$ of $\Ga(b;T\Pf)$ to construct 
an analogous subspace $\ti\Ga(\ups;T\Pf;\L)$ of $\ti\Ga(\ups;T\Pf)$
for $\ups\!\in\!\ti{\cal F}{\T}_{\de}^{\eset}$ sufficiently small.
For every
$$\ze\in\Ga_{-,+}(b;T\Pf;N_Y\ka)\!\oplus\!\Ga_{+,-}(b;T\Pf;N_{\Pf}Y),$$
we define $\ti{R}_{\ups}\ze\!\in\!\ti\Ga(\ups;T\Pf)$ by
$$\pi_Y^{\perp}\ti{R}_{\ups}\ze=\ti{R}_{\ups}\pi_Y^{\perp}\ze
\in\ti\Ga_-(\ups;\L)
\quad\hbox{and}\quad
\ti{R}_{\ups}\ze-\Pi_{\ups}^JR_{\ups}\ze\in
\Ga(\Si_{\ups};\ti{u}_{\ups}^*N_{\Pf}Y),$$
where again $R_{\ups}\ze\!=\!\ze\circ\ti{q}_{\ups}$.
By \e_ref{corr_e1} and~\e_ref{corr_e2},
\begin{equation}\label{corr_e3}
\big\|\ti{R}_{\ups}\ze-\Pi_{\ups}R_{\ups}\ze\big\|_{\ups,p,1}
\le C(b)\big|\rho(\ups)\big|^2\|\ze\|_{b,p,1}
\quad\forall\, \ze\!\in\!\Ga(b;\Pf;\L).
\end{equation}
Let $\ti\Ga_{-,+}(\ups;T\Pf;\ka)$, $\ti\Ga_{-,+}(\ups;T\Pf;N_Y\ka)$, 
$\ti\Ga_{+,-}(\ups;T\Pf;N_{\Pf}Y)$, and $\ti\Ga(\ups;T\Pf;\L)$
denote the images of $\Ga_{-,+}(b;T\Pf;\ka)$,  $\Ga_{-,+}(b;T\Pf;N_Y\ka)$, 
$\Ga_{+,-}(b;T\Pf;N_{\Pf}Y)$, and $\Ga(b;T\Pf;\L)$ under~$\ti{R}_{\ups}$.
By~\e_ref{corr_e3}, the~map
$$\pi_Y^{\perp}\!: \ti\Ga(\ups;T\Pf;\L)\lra\ti\Ga_-(\ups;\L)$$
is injective and thus an isomorphism.
Furthermore, since 
$$\ti\Ga_{-,+}(\ups;T\Pf;\ka)\subset\ti\Ga_-(\ups;T\Pf),$$
by the definition of the boundary operator $\d_0$ in~\e_ref{les_e} the map 
$$\d_{\ups}\!: \ti\Ga_{-,+}(\ups;T\Pf;N_Y\ka)
\!\oplus\!\ti\Ga_{+,-}(\ups;T\Pf;N_{\Pf}Y)
\lra H_J^1(\ti{b}(\ups);N_Y\ka), \qquad
\ze\lra[\pi_{\ka}^{\perp}D_{J,\ti{b}(\ups)}\ze],$$
is surjective and thus an isomorphism.
We set
\begin{alignat*}{1}
\ti\Ga_{-,+}^{0,1}(\ups;T\Pf;N_Y\ka)
&=\big\{\pi_{\ka}^{\perp}D_{J,\ti{b}(\ups)}\ze\!: 
\ze\!\in\!\ti\Ga_{-,+}(\ups;T\Pf;N_Y\ka)\big\}
\subset\ti\Ga^{0,1}(\ups;N_Y\ka)\!\equiv\! L^p(\ti{b}(\ups);N_Y\ka);\\
\ti\Ga_{+,-}^{0,1}(\ups;T\Pf;N_{\Pf}Y)
&=\big\{\pi_{\ka}^{\perp}D_{J,\ti{b}(\ups)}\ze\!: 
\ze\!\in\!\ti\Ga_{+,-}(\ups;T\Pf;N_{\Pf}Y)\big\}\subset
\ti\Ga^{0,1}(\ups;N_Y\ka).
\end{alignat*}
It follows from above that the projection map
$$\ti\pi^{0,1}_{\ups}\!:  \ti\Ga_{-,+}^{0,1}(\ups;T\Pf;N_Y\ka)
\!\oplus\!\ti\Ga_{+,-}^{0,1}(\ups;T\Pf;N_{\Pf}Y)
\lra H_J^1(\ti{b}(\ups);N_Y\ka)$$
is an isomorphism.\\

\noindent
The space $\Ga^{0,1}_-(b;N_Y\ka)$ of 
$u_b^*N_Y\ka$-valued harmonic $(0,1)$-forms on~$\Si_b$ splits~as
$$\Ga_-^{0,1}(b;N_Y\ka)=\Ga_{-;P}^{0,1}(b;N_Y\ka)\oplus\Ga_{-;B}^{0,1}(b;N_Y\ka)
={\cal H}_{b,P}\!\otimes\!\ev_P^*N_Y\ka \oplus \Ga_{-;B}^{0,1}(b;N_Y\ka).$$
Here $\Ga_{-;P}^{0,1}(b;N_Y\ka)$ and $\Ga_{-;B}^{0,1}(b;N_Y\ka)$ are the subspaces of 
$\Ga_-^{0,1}(b;N_Y\ka)$ consisting of the differentials supported
on the main components $\Si_{b,\aleph}$ of $\Si_b$ and on 
the only bubble component $\Si_{b,h}$ of~$\Si_b$, respectively.
Let 
$$\pi_{b;P}^{0,1},\pi_{b;B}^{0,1}\!: \Ga_-^{0,1}(b;N_Y\ka) \lra
\Ga_{-;P}^{0,1}(b;N_Y\ka), \, \Ga_{-;B}^{0,1}(b;N_Y\ka)$$
denote the projection maps.
If $\eta\!\in\!{\cal H}_{b,P}\!\otimes\!\ev_P^*N_Y\ka$, we define
$R_{\ups}\eta\!\in\!\Ga^{0,1}(\ups;N_Y\ka)$ as above by
identifying $N_Y\ka$ with the $g_{Y,b}$-orthogonal complement of $T\ka$ 
in $TY\!\subset\!T\Pf$.
If $\eta\!\in\!\Ga_{-;B}^{0,1}(b;N_Y\ka)$, let $R_{\ups}\eta\!=\!\ti{q}_{\ups}^*\eta$.
We denote~by
$$\ti\Ga_{-;P}^{0,1}(\ups;N_Y\ka),\,   \ti\Ga_{-;B}^{0,1}(\ups;N_Y\ka)
\subset \ti\Ga^{0,1}(\ups;N_Y\ka)$$
the images of $\Ga_{-;P}^{0,1}(b;N_Y\ka)$ and $\Ga_{-;B}^{0,1}(b;N_Y\ka)$
under the map $\ti{R}_{\ups}\!\equiv\!\Pi_{\ups}^JR_{\ups}$.
Let
$$\ti\pi^{0,1}_{\ups;P}\!:
 \ti\Ga^{0,1}(\ups;N_Y\ka)\lra\ti\Ga_{-;P}^{0,1}(\ups;N_Y\ka)
\quad\hbox{and}\quad
\ti\pi^{0,1}_{\ups;B}\!:
 \ti\Ga^{0,1}(\ups;N_Y\ka)\lra\ti\Ga_{-;B}^{0,1}(\ups;N_Y\ka)$$
be the $L^2$-projection maps.\\

\noindent
By \e_ref{kerprojest_e2}, \e_ref{pertest_e}, and~\e_ref{corr_e3},
\begin{equation}\label{pi0proj_e}
\big|\ti\pi_{\ups;P}^{0,1}D_{J,\ti{b}(\ups)}\ti{R}_{\ups}\ze
+2\pi\rho(\ups)J\, \ti{R}_{\ups}\pi_{\ka}^{\perp}\ti\D_{\T,h}^{\Pf}\ze\big|
\le C(b)\big|\rho(\ups)\big|^2\|\ze\|_{b,p,1}
~~~\forall\,\ze\!\in\!\Ga_{-,+}^{0,1}(\ups;T\Pf;N_Y\ka).
\end{equation}
In particular, the projection map
$$\ti\pi^{0,1}_{\ups;P}\!: \ti\Ga_{-,+}^{0,1}(\ups;T\Pf;N_Y\ka)
\lra\ti\Ga_{-;P}^{0,1}(\ups;N_Y\ka)$$
is an isomorphism. We denote its inverse by $S_{\ups;P}$.
The projection map
$$\ti\pi^{0,1}_{\ups;B}\!: \ti\Ga_{+,-}^{0,1}(\ups;T\Pf;N_{\Pf}Y)
\lra\ti\Ga_{-;1}^{B,1}(\ups;N_Y\ka).$$
is also an isomorphism, since the map $\d_b$ is.
We denote its inverse by $S_{\ups;B}$.\\

\noindent
Finally, let 
$$T_{\ups}=\ti\pi_{\ups}^{0,1}\circ (S_{\ups;P}\!\oplus\!S_{\ups;B})
\circ \ti{R}_{\ups}\!:
\Ga^{0,1}_-(b;N_Y\ka)\lra H_J^1(\ti{b}(\ups);N_Y\ka).$$
The maps $T_{\ups}$ with $\ups\!\in\!\ti{\cal F}{\T}_{\de}^{\eset}$ induce 
a bundle isomorphism
$$\Phi_{\T}\!:\pi_{{\cal FT}_{\de}^{\eset}}^*
\big(\wt\W_{\ka,d/d_{\ka}}^{1,1}|_{\U_{\T;\ka;1}}\big)
\lra \W_{\ka,d/d_{\ka}}^{0,1}\big|_{\M_1^0(\ka,d/d_{\ka})\cap U_{\T}^1},$$
covering $\phi_{\T}^1|_{{\cal FT}_{\de}^{\eset}}$. 
This isomorphism extends continuously over ${\cal FT}_{\de}\!-\!{\cal FT}_{\de}^{\eset}$,
as can be seen directly from the definition.\\

\noindent
If $\vt(\ups)\!\in\!\ti\Ga_-(\ups;\L)$, we can find a unique
\begin{equation*}\begin{split}
\ze_{\vt}(\ups)=&\, \ze_{\vt}^-(\ups)\oplus \ze_{\vt}^0(\ups)\oplus \ze_{\vt}^+(\ups)\\
&\in \ti\Ga_{-,+}(\ups;T\Pf;\ka)\oplus \ti\Ga_{-,+}(\ups;T\Pf;N_Y\ka)
\oplus \ti\Ga_{+,-}(\ups;T\Pf;N_{\Pf}Y)
\end{split}\end{equation*}
such that $\pi_Y^{\perp}\ze_{\vt}(\ups)\!=\!\vt(\ups)$.
By~\e_ref{pi0proj_e},
\begin{equation}\label{pi0proj_e2}
\big|\pi_{b;P}^{0,1}T_{\ups}^{-1}\d_0\vt(\ups)
+2\pi\rho(\ups)J\,\pi_{\ka}^{\perp}\ti\D_{\T,h}^{\Pf}\ti{R}_{\ups}^{-1}\ze_{\vt}^0(\ups)
\big| \le C_{\vt}(b)\big|\rho(\ups)\big|^2.
\end{equation}
On the other hand, by the definition of the map $\pi_k^-$ in Subsection~\ref{conestr_subs},
\begin{equation}\label{pi0proj_e3}
\pi_k^-\ti{\d}_{\ka,d/d_{\ka}}^{1,1}\vt(b)=
-2\pi J\, \pi_{\ka}^{\perp}\ti\D_{\T,h}^{\Pf}\ze_{\vt}^0(b).
\end{equation}
The estimate of Lemma~\ref{bdcontr_lmm} follows from 
\e_ref{corr_e3}, \e_ref{pi0proj_e2}, \e_ref{pi0proj_e3}, and the continuity of
the section~$\vt$.

\end{document}